\documentclass[reqno,11pt]{amsart}
\usepackage[foot]{amsaddr}
\usepackage{amssymb,amsmath,amsthm,amsfonts}
\usepackage{mathrsfs,dsfont,comment,mathscinet,mathtools}
\usepackage{regexpatch}
\usepackage{graphicx,pgfplots}

\usepackage[utf8]{inputenc}
\usepackage[normalem]{ulem}
\usepackage[left=2.5cm,right=2.5cm,top=2.5cm,bottom=2.5cm]{geometry}

\usepackage{pifont}
\usepackage{tkz-fct}

\usepackage{graphicx,tikz,color}
\usetikzlibrary{quotes,angles}

\usepackage[format=hang,labelfont=bf]{caption}

\usepackage{enumerate,esint}

\usepackage{ wasysym }

\usepackage{epstopdf}
\usepackage[colorlinks=true, pdfstartview=FitV, linkcolor=blue, 
            citecolor=blue, urlcolor=blue]{hyperref}
            
\allowdisplaybreaks
\numberwithin{equation}{section}
\theoremstyle{plain}
\newtheorem{theorem}{Theorem}[section]
\newtheorem{proposition}[theorem]{Proposition}
\newtheorem{lemma}[theorem]{Lemma}
\newtheorem{example}[theorem]{Example}

\newtheorem{definition}[theorem]{Definition}
\theoremstyle{definition}

\newtheorem*{theorem*}{Theorem}

\definecolor{mblue}{HTML}{13439b}

\newcommand{\N}{\mathbb{N}}
\newcommand{\R}{\mathbb{R}}
\newcommand{\RN}{\mathbb{R}^N}
\newcommand{\RNN}{\mathbb{R}^{N\times N}}

\renewcommand{\H}{\mathscr{H}}

\newcommand{\Tan}{\mathrm{T}}

\renewcommand{\O}{\mathrm{O}}
\renewcommand{\d}{\mathrm{d}}
\newcommand{\restr}[1]{|_{#1}}

\renewcommand{\dag}{\text{\textdied}}

\def\XXint#1#2#3{{\setbox0=\hbox{$#1{#2#3}{\int}$ }
\vcenter{\hbox{$#2#3$ }}\kern-.6\wd0}}

\makeatletter
\newcommand{\subalign}[1]{%
	\vcenter{%
		\Let@ \restore@math@cr \default@tag
		\baselineskip\fontdimen10 \scriptfont\tw@
		\advance\baselineskip\fontdimen12 \scriptfont\tw@
		\lineskip\thr@@\fontdimen8 \scriptfont\thr@@
		\lineskiplimit\lineskip
		\ialign{\hfil$\m@th\scriptstyle##$&$\m@th\scriptstyle{}##$\hfil\crcr
			#1\crcr
		}%
	}%
}
\makeatother

\newcommand{\mres}{\mathbin{\vrule height 1.6ex depth 0pt width
		0.13ex\vrule height 0.13ex depth 0pt width 1.3ex}}

\makeatletter
\newcommand{\subsetcong}{\mathrel{\mathpalette\subseteq@cong\relax}}
\newcommand{\subsetsim}{\mathrel{\mathpalette\subset@sim\relax}}

\newcommand{\subseteq@cong}[2]{%
	\vbox{\offinterlineskip\m@th
		\ialign{\hfil$#1##$\hfil\cr
			\sim\cr\subseteq\cr
		}%
	}%
}

\newcommand{\subset@sim}[2]{%
	\vbox{\offinterlineskip\m@th
		\ialign{\hfil$#1##$\hfil\cr
			\sim\cr\subset\cr
		}%
	}%
}

\makeatother

\newcommand{\cof}{\mathrm{cof}}
\newcommand{\adj}{\mathrm{adj}}

\newcommand{\per}{\mathrm{Per}}
\renewcommand{\div}{\mathrm{div}}
\newcommand{\dist}{\mathrm{dist}}
\newcommand\wk{\rightharpoonup}
\newcommand{\wks}{\overset{\ast}{\rightharpoonup}}

\newcommand*\closure[1]{\overline{#1}}

\newcommand{\leb}{\mathscr{L}^N}
\newcommand{\haus}{\mathscr{H}^{N-1}}

\newcommand{\imt}{\mathrm{im}_{\rm T}}
\newcommand{\img}{\mathrm{im}_{\rm G}}
\newcommand{\domg}{\mathrm{dom}_{\rm G}}
\newcommand{\Det}{\mathrm{Det}}

\newcommand\EEE{\color{black}}

\newcommand{\MMM}{\color{magenta}} 

\makeatletter
\xpatchcmd{\proof}{\itshape}{\bfseries}{}{}
\makeatother

\title[Core-radius approximation of singular minimizers in nonlinear elasticity]{Core-radius approximation of singular minimizers in nonlinear elasticity}
\author[M. Bresciani]{Marco Bresciani${}^{*}$}
\address{* Department of Mathematics, Friedrich-Alexander-Universit\"{a}t Erlangen-N\"{u}rnberg, Cauerstrasse 11, 91058 Erlangen, Germany}
\email{marco.bresciani@fau.de}
\author[M. Friedrich]{Manuel Friedrich${}^{*}$}
\email{manuel.friedrich@fau.de}

\date{\today}
\begin{document}

\setlength\parindent{0pt}

\vskip .2truecm
\begin{abstract}
  We study a variational model  in  nonlinear elasticity allowing for cavitation which penalizes both the volume and the perimeter of the cavities. Specifically, we investigate the approximation (in the sense of $\Gamma$-convergence) of the  energy  by means of functionals defined on perforated domains. Perforations are introduced  at  flaw points where singularities are expected  and, hence,   the corresponding deformations  do not exhibit cavitation.   Notably, those points are not prescribed but  rather selected by the variational principle. 
Our analysis is motivated by the numerical simulation of cavitation and extends previous results  on models which solely accounted for  elastic energy but neglected contributions related to the formation of cavities.    
\end{abstract}
\maketitle

\section{Introduction}

\subsection{Motivation}
In the context of nonlinear elasticity, the term cavitation commonly refers to the abrupt formation of voids inside of solids in response to tensile stresses \cite{ball.1982}. To date, this phenomenon has been experimentally observed in different  types  of materials ranging from rubber-like solids \cite{gent.lindley,gent.wang} to ductile metals \cite{gao,petrinic}. Specifically, numerous studies advance the hypothesis that cavitation might play a key role in   fracture initiation by means of growth and coalescence of cavities \cite{tvergaard},  especially   for soft materials \cite{gent.wang,lopezpamies3}.  

Mechanically, cavitation has been interpreted as an elastic instability determined by the sudden expansion of pre-existing micro-voids \cite{fond}. According to this interpretation, for  a  hyperelastic solid with reference configuration $\Omega\subset \RN$ and strain-energy density $W\colon \RNN_+ \to [0,+\infty)$, equilibrium configurations correspond to deformations $\boldsymbol{y}\colon \Omega \to \RN$ minimizing the elastic energy
\begin{equation}
	\label{eq:intro-elastic}
	\boldsymbol{y}\mapsto \int_\Omega W(D\boldsymbol{y}(\boldsymbol{x}))\,\d\boldsymbol{x}.
\end{equation}
Even though this variational principle has been capable of explaining various empirical  measurements   \cite{ball.1982}, some more refined experiments   (e.g., the renowned tests  on ``poker chips'' \cite{gent.lindley})   could not be completely justified within this framework \cite{cho.gent,gent.park,lopezpamies1}.  Therefore, an alternative modeling approach  has been proposed in recent years \cite{lopezpamies2}. Looking at cavitation  as a phenomena of creation of new surface, this approach is based  on the minimization of the energy
\begin{equation}
	\label{eq:intro-elastic+perimeter}
	\boldsymbol{y}\mapsto \int_\Omega W(D\boldsymbol{y}(\boldsymbol{x}))\,\d\boldsymbol{x} + \sum_{\boldsymbol{a}\in C_{\boldsymbol{y}}} \per \left( \imt(\boldsymbol{y},\boldsymbol{a})\right).
\end{equation} 
In the previous equation, $C_{\boldsymbol{y}}$ denotes the set of points $\boldsymbol{a}\in \Omega$ in correspondence of which the deformation $\boldsymbol{y}$ opens the cavity $\imt(\boldsymbol{y},\boldsymbol{a})$, and $\per \left(\imt(\boldsymbol{y},\boldsymbol{a})\right)$ stands for its perimeter, i.e., the area of its boundary.  

From an analytical perspective,   the existence of minimizers for the energy \eqref{eq:intro-elastic+perimeter} has been established in \cite{henao.moracorral.lusin} by building upon previous works for a similar model \cite{mueller.spector}.

From the computational point of view, the simulation of cavitation faces substantial challenges. A major difficulty comes from the fact that minimization problems of nonlinear elasticity often exhibit the Lavrentiev phenomenon \cite{almi.kroemer.molchanova,ball.compu,foss.hrusa.mizel}, which  makes traditional finite-element methods ineffective. Although various   remedies have been proposed to  overcome this issue (see, e.g., \cite{ball.compu,ball.knowles,knowles,li}), these remain mostly inadequate for the numerical simulation of cavitation. The reasons are mainly their inability to effectively resolve the singularities and their proneness towards mesh degeneracy \cite[Subsection 3.2]{henao.xu}. As a practicable way, we mention the works \cite{henao.moracorral.xu,henao.moracorral.xu2} dealing with a phase-field approximation of the energy  \eqref{eq:intro-elastic+perimeter}.  

In this work, we propose a novel approach for the numerical minimization of the energy \eqref{eq:intro-elastic+perimeter} which we analyze from a theoretical perspective. Inspired by previous works in the  framework of pure elasticity   (i.e., dealing with  the minimization of \eqref{eq:intro-elastic}, see the brief review of the literature provided in Subsection \ref{subsec:literature}),   our approach hinges upon the idea of regularizing the energy \eqref{eq:intro-elastic+perimeter} by introducing some perforations inside the reference domain in correspondence of the points where deformations exhibit singularities. Borrowing the terminology usually employed in other fields such as the theories of dislocations and superconductors, we refer to this approach as core-radius approximation.  The method is  detailed in Subsection \ref{subsec:contribution}.

\subsection{Relevant literature}\label{subsec:literature}
We summarize briefly the previous contributions from the literature on the core-radius approximation of the energy \eqref{eq:intro-elastic}. The starting point was the model by Sivaloganathan and Spector \cite{sivaloganathan.spector}. In that work,  the variational approach in \eqref{eq:intro-elastic} is  adopted  under the  assumption that cavitation may only occur  at a prescribed (finite) set of flaw points in the reference configuration. Denoting by $\boldsymbol{a}_1,\dots,\boldsymbol{a}_M\in \Omega$ with $M\in\N$ the flaw points, this requirement is enforced by considering admissible deformations $\boldsymbol{y}\in W^{1,p}(\Omega;\RN)$ with $p>N-1$ whose distributional determinant $\Det D \boldsymbol{y}$ satisfies 
\begin{equation}
	\label{eq:intro-Det}
	\Det D \boldsymbol{y}=(\det D\boldsymbol{y})\leb + \sum_{\boldsymbol{i}=1}^M \kappa_i^{\boldsymbol{y}}\delta_{\boldsymbol{a}_i} \quad  \quad \text{in $\mathcal{M}_{\rm b}(\Omega)$}
\end{equation} 
for some $\kappa_1^{\boldsymbol{y}},\dots,\kappa_M^{\boldsymbol{y}}\geq 0$. Here,  $\kappa_i^{\boldsymbol{y}}$ gives the volume of the cavity opened by $\boldsymbol{y}$ at the point $\boldsymbol{a}_i$.
Note that the representation formula \eqref{eq:intro-Det} is natural and can actually be inferred from the boundedness of suitable surface energies, see \cite[Theorem 4.6]{henao.moracorral.lusin} and \cite[Theorem 8.4]{mueller.spector}. The existence of minimizers of the energy  \eqref{eq:intro-elastic} among all orientation-preserving deformations satisfying \eqref{eq:intro-Det}  is  established in \cite[Theorem 4.1]{sivaloganathan.spector}. 
  Specifically,    admissible deformations are required to fulfill condition (INV) by M\"{u}ller and Spector \cite{mueller.spector}, which is a technical condition closely related to invertibility that excludes unphysical behaviors.

%Precisely, the boundary condition is imposed by requiring the existence of a suitable extension to a larger domain which does not create additional cavities in order to exclude cavitation at the boundary \cite[p. 54--55]{mueller.spector}.  
 
Subsequently, the following approximation for minimizers of \eqref{eq:intro-elastic} under \eqref{eq:intro-Det} has been suggested by  Sivaloganathan,  Spector, and Tilakraj \cite{sivaloganathan.spector.tilakraj}: given $0<\varepsilon \ll 1$, one looks for minimizers of the energy 
\begin{equation}
	\label{eq:intro-elastic-eps}
	\boldsymbol{y}\mapsto \int_{\Omega_\varepsilon} W(D\boldsymbol{y}(\boldsymbol{x}))\,\d \boldsymbol{x}, 
\end{equation} 
where $\Omega_\varepsilon\coloneqq \Omega \setminus \bigcup_{i=1}^M \closure{B}(\boldsymbol{a}_i,\varepsilon)$,   among all orientation-preserving deformations $\boldsymbol{y}\in W^{1,p}(\Omega_\varepsilon;\RN)$ that  satisfy (INV) and do  not exhibit cavitation in the sense that
\begin{equation}
	\label{eq:intro-Det-eps}
	\Det D \boldsymbol{y}=(\det D \boldsymbol{y})\leb \quad \text{in $\mathcal{M}_{\rm b}(\Omega_\varepsilon)$.}
\end{equation}
For $M=1$, it is indeed proved in \cite[Theorem 4.1]{sivaloganathan.spector.tilakraj} that minimizers of the energy  \eqref{eq:intro-elastic-eps} under \eqref{eq:intro-Det-eps}, whose existence is ensured by \cite[Theorem 4.1]{sivaloganathan.spector}, converge, as $\varepsilon\to 0^+$, to minimizers of the energy  \eqref{eq:intro-elastic} under \eqref{eq:intro-Det}. 

A prime difficulty  related to condition (INV)  arose  in the case $M>1$,  and  was  later \EEE  overcome by Henao \cite{henao}.  To do so, he introduced \EEE  a strengthened version of condition (INV), named condition (INV${}^\prime$),  to be imposed on competitors for \eqref{eq:intro-elastic-eps}--\eqref{eq:intro-Det-eps}.  In this way, he was  able to extend  the convergence result  of  \cite{sivaloganathan.spector.tilakraj} to the case  $M>1$,  see   \cite[Theorem 5.1]{henao}.   Subsequently,   this result was further generalized to the setting where  the set of flaw points is not prescribed but  where a minimization  among all possible sets of $M$ flaw points in $\Omega$  is performed, see   \cite[Theorem 6.3]{henao.thesis}. 

Concerning numerical implementation, finite-element analyses of the regularized energy \eqref{eq:intro-elastic-eps} for fixed flaw points have been carried out by Henao and Xu \cite{henao.xu}, and  by  Lian and Li \cite{lian.li,lian.li2}.

\subsection{Contributions of the paper}\label{subsec:contribution}
This study aims at investigating the theoretical viability of the core-radius approximation  in \cite{henao,henao.thesis,sivaloganathan.spector.tilakraj} for the minimization of the energy \eqref{eq:intro-elastic+perimeter}. Precisely, following  \cite{sivaloganathan.spector}, we augment the energy by an additional term penalizing the volume of the cavities.

As in \cite{henao.thesis}, the set of flaw points is not prescribed but rather constitutes one of the unknowns. Thus, for given $\lambda_{\rm v},\lambda_{\rm p}\geq 0$, we are interested in the minimization of the functional
\begin{equation}
	\label{eq:intro-E}
	(A,\boldsymbol{y})\mapsto \int_\Omega W(D\boldsymbol{y}(\boldsymbol{x}))\,\d\boldsymbol{x}+\lambda_{\rm v} \sum_{\boldsymbol{a}\in A} \leb(\imt(\boldsymbol{y},\boldsymbol{a}))+\lambda_{\rm p} \sum_{\boldsymbol{a}\in A} \per \left(\imt(\boldsymbol{y},\boldsymbol{a})\right),
\end{equation} 
where  $A$ denotes the set of flaw points. This set has to be contained in a given compact subset $H\subset \Omega$ and its cardinality cannot exceed a given integer $M\in\N$.  In \cite{henao,henao.thesis,sivaloganathan.spector.tilakraj}, Dirichlet boundary conditions on deformations were formulated by requiring the existence of a suitable extension in order to exclude cavitation at the boundary \cite[pp.~54--55]{mueller.spector}.  
The physical confinement of flaw points assumed here provides an alternative way to  rule out this phenomena, so that boundary conditions can be imposed on a portion of the boundary in the usual sense of traces. \EEE 
Similarly to \eqref{eq:intro-Det}, admissible deformation satisfy 
\begin{equation}
	\label{eq:intro-D}
	\Det D \boldsymbol{y}=(\det D  \boldsymbol{y})\leb+\sum_{\boldsymbol{a}\in A} \kappa_{\boldsymbol{a}}^{\boldsymbol{y}} \delta_{\boldsymbol{a}} \quad \text{in $\mathcal{M}_{\rm b}(\Omega)$}
\end{equation}
for some $\{ \kappa_{\boldsymbol{a}}^{\boldsymbol{y}}:\:\boldsymbol{a}\in A\}\subset (0,+\infty]$, and  fulfill condition (INV).

For the approximation, we introduce a core-radius $0<\varepsilon\ll 1$. We restrict ourselves to sets $A$ of flaws points that are well-separated in terms of $\varepsilon$ and, setting $\Omega^A_\varepsilon\coloneqq \Omega \setminus \bigcup_{\boldsymbol{a}\in A} \closure{B}(\boldsymbol{a},\varepsilon)$, we consider deformations $\boldsymbol{y}\in W^{1,p}(\Omega^A_\varepsilon;\RN)$.  Then,  we look for minimizers of the energy
\begin{equation}
	\label{eq:intro-E-eps}
	(A,\boldsymbol{y})\mapsto \int_{\Omega^A_\varepsilon} W(D\boldsymbol{y}(\boldsymbol{x}))\,\d\boldsymbol{x}+\lambda_{\rm v} \sum_{\boldsymbol{a}\in A} \leb \left( \imt(\boldsymbol{y},B(\boldsymbol{a},\varepsilon))\right)+ \lambda_{\rm p} \sum_{\boldsymbol{a}\in A} \per \left( \imt(\boldsymbol{y},B(\boldsymbol{a},\varepsilon))\right). 
\end{equation}
In the previous equation, $\imt(\boldsymbol{y},B(\boldsymbol{a},\varepsilon))$ denotes the topological image of $B(\boldsymbol{a},\varepsilon)$ under $\boldsymbol{y}$. Roughly speaking, this set corresponds to the region enclosed by $\boldsymbol{y}(S (\boldsymbol{a},\varepsilon))$,  where $S(\boldsymbol{a},\varepsilon)$ denotes the sphere centered at $\boldsymbol{a}$ with radius $\varepsilon$.  \EEE
As in \cite{henao,henao.thesis}, approximating deformations are required to satisfy condition  (INV${}^\prime$). Also, they have to fulfill Dirichlet boundary conditions together with some technical assumptions ensuring that the energy \eqref{eq:intro-E-eps} is well defined. 
 Condition   \eqref{eq:intro-Det-eps}  is strengthened by requiring,  with $\widetilde{\Omega}^A_\varepsilon\coloneqq \Omega  \setminus  \bigcup_{\boldsymbol{a}\in A} B(\boldsymbol{a},\varepsilon)$, that 
\begin{equation}
	\label{eq:intro-D-eps}
	\Det^A_\varepsilon D \boldsymbol{y}=(\det D\boldsymbol{y})\leb \quad \text{in $\mathcal{M}_{\rm b}(\widetilde{\Omega}^A_\varepsilon)$,}
\end{equation}
 where  $\Det^A_\varepsilon D \boldsymbol{y}$ denotes the extended distributional determinant,   i.e., a generalization of the distributional determinant for perforated domains.  
This is the right object to measure the volume of topological images without accounting for the volume enclosed by the perforations. Thus, \eqref{eq:intro-D-eps} ensures that cavitation does not occur at the boundary of perforations. 
  \EEE
 
We observe that, as  a consequence of the technical assumptions of approximating deformations, the last two terms in \eqref{eq:intro-E-eps} can be rewritten using boundary integrals. Precisely,
 \begin{equation}
 	\label{eq:v}
 	 \leb(\imt(\boldsymbol{y},B(\boldsymbol{a},\varepsilon)))=\int_{S(\boldsymbol{a},\varepsilon)}\boldsymbol{y}(\boldsymbol{x})\cdot \left( (\cof D^{\rm t}\boldsymbol{y}(\boldsymbol{x}))\boldsymbol{\nu}_{B(\boldsymbol{a},\varepsilon)}(\boldsymbol{x}) \right)\,\d\haus(\boldsymbol{x})
 \end{equation}
 and
 \begin{equation}
 	\label{eq:p}
 	 \per \left( \imt(\boldsymbol{y},B(\boldsymbol{a},\varepsilon))\right)=\frac{1}{N} \int_{S(\boldsymbol{a},\varepsilon)} \left | (\cof D^{\rm t}\boldsymbol{y}(\boldsymbol{x}))\boldsymbol{\nu}_{B(\boldsymbol{a},\varepsilon)}(\boldsymbol{x}) \right |\,\d\haus(\boldsymbol{x})
 \end{equation}
 for all $\boldsymbol{a}\in A$. In the previous two equations, $D^{\rm t}\boldsymbol{y}$ denotes the weak tangential gradient of $\boldsymbol{y}$. 
 This observation is relevant for  applications since,   compared with those in \eqref{eq:intro-E-eps}, the expressions in \eqref{eq:v}--\eqref{eq:p} seem more manageable from the computational point of view because   they  do  not involve the topological image.   

Our main findings concern the convergence  of the minimization problem   determined by \eqref{eq:intro-E-eps}--\eqref{eq:intro-D-eps}, as $\varepsilon\to 0^+$, to \EEE that given by \EEE \eqref{eq:intro-E}--\eqref{eq:intro-D}  in the sense of $\Gamma$-convergence \cite{dalmaso}. First, in Theorem \ref{thm:equi-coercivity}, we establish the relevant compactness result   by   adapting the arguments in \cite[Theorem 6.3]{henao.thesis} to our setting. Next, in Theorem \ref{thm:lb} we show  that the energy \eqref{eq:intro-E} under \eqref{eq:intro-D} provides a lower bound for the energies \eqref{eq:intro-E-eps} under \eqref{eq:intro-D-eps}, as $\varepsilon \to 0^+$. Eventually, the optimality of the lower bound is proved in Theorem \ref{thm:ub} under an additional regularity assumption. Specifically, the existence of recovery sequences is established for any admissible state $(A,\boldsymbol{y})$ which, roughly speaking, satisfies
\begin{equation}
	\label{eq:intro:P}
	\liminf_{r \to 0^+} \per \left( \imt(\boldsymbol{y},B(\boldsymbol{a},r))  \right)=\per \left( \imt(\boldsymbol{y},\boldsymbol{a})\right) \quad \text{for all $\boldsymbol{a}\in A$.}  
\end{equation}
Some comments on this assumption are provided in Subsection \ref{subsec:overview}.
  
We mention that, for $\lambda_{\rm p}=0$, the previous assumption can be disregarded so that we have a full $\Gamma$-convergence result. This  already provides an extension to \cite[Theorem 6.3]{henao.thesis} by including an energy term related to the formation of cavities. Additionally, in this case, standard $\Gamma$-convergence arguments \cite[Theorem 7.4]{dalmaso} yield the existence of minimizers for the $\Gamma$-limit. 

Note that the same conclusion cannot be inferred for  $\lambda_{\rm p}>0$   because of the restriction \eqref{eq:intro:P}. Still, as we show in Theorem \ref{thm:existence}, the minimization problem   determined by  \eqref{eq:intro-E}--\eqref{eq:intro-D} admits solutions.

 From \MMM a \EEE technical point of view, the proofs of our main results  rely on the adaptation of various lemmas and propositions from \cite{henao.moracorral.lusin,mueller.spector} concerning the relationship between geometric and topological image to the setting of deformations on perforated domains (see Section \ref{sec:perforated}). In particular, we introduce the extended distributional determinant which represents an essential tool within our arguments. Despite its definition being natural, this object seems (at least to our knowledge) new, although it appears closely connected to the extended surface energy defined in \cite{henao.moracorral.xu}. \EEE

\subsection{Discussion and outlook}
\label{subsec:overview}
 
Finding general criteria for the validity of \eqref{eq:intro:P} does not seem feasible.  In particular, this condition is not  stable with respect to weak convergence.
By means of explicit examples,  in Section \ref{sec:example} \EEE we   examine  the validity of \eqref{eq:intro:P}  in relation to  the shape of the cavities (Example \ref{ex:radial}), the choice of the reference configuration (Example \ref{ex:cor}), and the composition with piecewise-smooth transformations (Example~\ref{ex:superpos}). 
Deformations violating \eqref{eq:intro:P}   typically  exhibit the following behavior, as $r\to 0^+$, for at least one $\boldsymbol{a}\in A$: either the (reduced) boundary of the set $\imt(\boldsymbol{y},B(\boldsymbol{a},r))$ oscillates  wildly with oscillations of vanishing amplitude or it produces an exterior kink with two sides that squash onto each other in the limit. To our view, both behaviors can be regarded as pathological.
While the first option seems to be ruled out by  (almost everywhere) injectivity,  the second one can actually occur within our class of admissible deformation. Such an instance is discussed in Example \ref{ex:counterexample}. This example provides a deformation on a planar domain that opens a single cavity formed by the union of a closed ball with a segment connecting one point of its boundary with a point of its complement. Thus, the reduced boundary coincides with the sphere enclosing the volume of the cavity. In particular, the outer segment is not detected by the energy \eqref{eq:intro-elastic+perimeter}. This observation shows that the impossibility of enforcing condition \eqref{eq:intro:P} on admissible deformations is somehow related to the limitations of the model based on \eqref{eq:intro-elastic+perimeter}. Even replacing the perimeter in \eqref{eq:intro-elastic+perimeter} with the area of the topological boundary (leaving aside the question of the existence of minimizers) would not solve the issue as \eqref{eq:intro:P}  requires  to count the   length  of the exterior segment twice. This situation partly resembles that of deformations creating invisible surface in \cite[Theorem 3]{henao.moracorral.fracture}.
Anyways, having \eqref{eq:intro-elastic+perimeter} in mind, the violation of \eqref{eq:intro:P} does not seem energetically favorable. Hence, at least on an intuitive level, this requirement should not constitute a severe restriction on competitors. 

Concerning the well posedness of the model, it must be noted that the existence of minimizers of \eqref{eq:intro-E} under \eqref{eq:intro-D} achieved in Theorem \ref{thm:existence} can   alternatively be deduced from \cite[Theorem 8.5]{henao.moracorral.lusin} by using the representation formula for the surface energy functional \cite[Theorem 4.6]{henao.moracorral.lusin}. Here, we provide a more self-contained proof that does  not  involve this notion.  

In contrast, the existence of minimizers  of  \eqref{eq:intro-E-eps} under \eqref{eq:intro-D-eps} seems out of reach,  at least for $\lambda_{\rm p} > 0 $.   Indeed, the class of approximating deformations is not closed with respect to weak convergence. This comes from the fact that the set flaw points is not prescribed as well as the nature of the technical requirements defining approximating deformations.  From our perspective, this shortcoming does not compromise the viability of the core-radius approach for numerical simulations. Indeed, the implementation of this approach  would rely on an  approximation of the energy \eqref{eq:intro-E-eps} by finite elements for which hopefully existence might be proved.  \EEE

 The relaxation of the energy \eqref{eq:intro-E-eps} is expected to be related to \v{S}ilhav\'{y}'s theory of phase transitions with interfacial energies \cite{silhavy,silhavy2}. 
In the relaxed energy, the last term in \eqref{eq:intro-E-eps}, rewritten as in \eqref{eq:p}, should correspond to the total variation of some  measure   concentrated   on the boundaries of the perforations that is related to the deformation by suitable integral identities. This object should generalize the measure $\sum_{\boldsymbol{a}\in {A}}(\cof D^{\rm t}\boldsymbol{y})\boldsymbol{\nu}_{B(\boldsymbol{a},\varepsilon)}\haus\mres S(\boldsymbol{a},\varepsilon)$.   
 \EEE

 Following the strategy in \cite[Lemma 8.1]{mueller.spector}, we  show  that the extended distributional determinant \EEE is a Radon measure on the union of the perforated domain with the boundaries of the perforations.  It would be interesting to establish whether, analogously to \cite[Theorem 4.6]{henao.moracorral.lusin} and \cite[Theorem 8.4]{mueller.spector}, the boundedness of the energy \eqref{eq:intro-E-eps} yields a representation formula for this object. In principle, this measure might charge the boundary of the perforations. 
Clearly, the extended distributional determinant can  be introduced for general domains by  considering the whole boundary in place of  the boundaries of the perforations. In this setting, this tool might turn our to be useful for characterizing cavitation at the boundary. See \cite{henao.moracorral.oliva} for some  initial work in this direction. \EEE

As we said in the beginning, our analysis  is  motivated by the numerical minimization of the energy \eqref{eq:intro-elastic+perimeter}. Therefore, the next natural step would require to address the discretizaion of the approximating energy \eqref{eq:intro-E-eps} by finite elements. This task is left to future investigations.

\subsection{Structure of the paper}
The paper is organized  as follows.  The setting of this study and our main findings are presented in Section \ref{sec:setting}. Preliminary results from the literature are collected in Section \ref{sec:prelim}. In Section \ref{sec:existence}, we prove Theorem \ref{thm:existence}. Results for deformations on perforated domains are established in Section \ref{sec:perforated}. The proofs of our main results regarding the core-radius approximation, i.e., Theorem \ref{thm:equi-coercivity}, Theorem \ref{thm:lb}, and Theorem \ref{thm:ub}, are presented in Section \ref{sec:proof}. In Section \ref{sec:example}, we discuss some examples of deformations in connection  with assumption \eqref{eq:intro:P}.  Eventually, in the final Appendix, we revise some known facts about tangential differentiability.   

\EEE

\subsection{Notation}
In this work, $N$ denotes an integer with $N\geq 2$.  By $\RNN_+$ we indicate the set of matrices in $\RNN$ with positive determinant. \EEE  The symbols $B(\boldsymbol{a},r)$ and $S(\boldsymbol{a},r)$  are used for   the open ball and the sphere in $\RN$, respectively, centered at $\boldsymbol{a}\in\RN$ with radius $r>0$. We denote the closure of $B(\boldsymbol{a},r)$ by $\closure{B}(\boldsymbol{a},r)$.
 Also, for $R>r$, we define the open annulus $A(\boldsymbol{a},r,R)\coloneqq B(\boldsymbol{a},R)\setminus \closure{B}(\boldsymbol{a},r)$ and we write $\closure{A}(\boldsymbol{a},r,R)$ for its closure.  
 Given $E\subset \RN$, we use the notation $E^\circ$, $\closure{E}$, and $\partial E$ for its interior, closure, and boundary, respectively. 

The Lebesgue measure and the $(N-1)$-dimensional Hausdorff measure on $\RN$ are indicated by $\leb$ and $\haus$.  Unless differently stated, all expressions like ``almost everywhere'' or ``almost every'' as well as all measurability properties refer to the Lebesgue measure.  
Given $E,F\subset \RN$, we write 
$E\cong F$ and $E\simeq F$ whenever $\leb(E \triangle E)=0$ or $\haus(E\triangle F)=0$, respectively. Also, given $\boldsymbol{u},\boldsymbol{v}\colon \RN \to \RN$, we use the notations $\boldsymbol{u}\cong \boldsymbol{v}$ in $E$  and $\boldsymbol{u}\simeq \boldsymbol{v}$ on $E$  to indicate  that  $\leb(\{ \boldsymbol{x}\in E: \boldsymbol{u}(\boldsymbol{x})=\boldsymbol{v}(\boldsymbol{x}) \})=0$ and $\haus(\{ \boldsymbol{x}\in E: \boldsymbol{u}(\boldsymbol{x})=\boldsymbol{v}(\boldsymbol{x}) \})=0$, respectively.

We use standard notation for spaces of continuous and differentiable functions, Lebesgue spaces, and Sobolev spaces. 
 For an open set $G\subset \RN$, the symbol $\mathcal{D}(G)$ is used for the class of smooth and compactly supported functions on $G$.  The symbol  $\mathcal{M}_{\rm b}(G)$ stands for the class of bounded Radon measures defined on the Borel subsets of $G$ with $\delta_{\boldsymbol{a}}$ denoting the Dirac delta centered at $\boldsymbol{a}\in G$. \EEE 
\emph{Whenever $G$ is a bounded Lipschitz domain and $\boldsymbol{y}\in W^{1,p}(G;\RN)$ with $p>N-1$, we adopt the convention of denoting the trace of $\boldsymbol{y}$ on $\partial G$ by $\boldsymbol{y}^*$}. As explained, e.g, in \cite[Remark 4.4.5]{ziemer.wdf}, this notation is compatible with that of precise representatives.

 If $E\subset \RN$ is a set of  finite perimeter \cite[Chapter 3]{ambrosio.fusco.pallara},  then we denote by $\per E$ and $\partial^* E$ its perimeter and reduced boundary, respectively. Also,  $\boldsymbol{\nu}_E\colon \partial^* E \to \RN$ stands for its outer unit normal.

\section{Setting and main results}
\label{sec:setting}

Throughout the paper, we assume that
\begin{equation*}
	\text{$N\in\N$ satisfies $N\geq 2$,\qquad   $\Omega \subset \RN$ is a bounded simply connected  domain}, 
\end{equation*}
\begin{equation*}
	\text{$H\subset \Omega$ is a compact set, \qquad $\Gamma \subset \partial \Omega$ is a rectifiable subset with $\haus(\Gamma)>0$.}
\end{equation*}
Moreover, we fix
\begin{equation*}
	M\in \N, \qquad p>N-1, \qquad \boldsymbol{d}\in W^{1,p}(\Omega;\RN).
\end{equation*}
The set $\Omega$ represents the reference configuration of  a  hyperelastic body. The exponent $p$ specifies the integrability of deformations. These must fulfill a Dirichlet boundary condition on $\Gamma$ with datum $\boldsymbol{d}$ and  are  allowed to open at most $M$ cavities in correspondence of points confined in the set  $H$.

\subsection{Variational model and its regularization}  
We denote by $\mathcal{A}$ the collection of all sets $A\subset H$ with $\mathscr{H}^0(A)\leq M$. 

\begin{definition}[Deformations]\label{def:deformation}
Given $A\in\mathcal{A}$, we define  $\mathcal{Y}^A$ as the class of all  $\boldsymbol{y}\in W^{1,p}(\Omega;\RN)$ satisfying the following properties:
\begin{enumerate}[\rm (i)]
	\item $\det D \boldsymbol{y}\in L^1(\Omega)$ and $\det D \boldsymbol{y}(\boldsymbol{x})>0$ for almost all  $\boldsymbol{x}\in \Omega$;
	\item $\boldsymbol{y}$ satisfies condition  {\rm (INV)};
	\item $\boldsymbol{y}^*\simeq \boldsymbol{d}^*$ on $\Gamma$;
	\item there exist $\{ \kappa_{\boldsymbol{a}}^{\boldsymbol{y}}:\:\boldsymbol{a}\in A \}\subset (0,+\infty)$ such that 
	\begin{equation*}
		\Det D {\boldsymbol{y}}=(\det D {\boldsymbol{y}})\leb + \sum_{\boldsymbol{a}\in A} \kappa_{\boldsymbol{a}}^{\boldsymbol{y}} \delta_{\boldsymbol{a}} \quad \text{in $\mathcal{M}_{\rm b}(\Omega)$.}
	\end{equation*}
\end{enumerate}
\end{definition}

 Item (i) restricts our class to orientation-preserving deformations. 
Item (ii) refers to condition (INV)  recalled in Definition \ref{def:INV}. This is a topological condition which excludes pathological behaviors and entails the almost everywhere injectivity \cite[Lemma~3.4]{mueller.spector}.  Item (iii) imposes the Dirichlet boundary condition in the sense of traces.  In (iv), we denote by $\Det D \boldsymbol{y}$ the distributional determinant of $\boldsymbol{y}$ recalled in Definition \ref{def:distributional-determinant}. In view of  item (iv) and Lemma \ref{lem:Det} below,   $\boldsymbol{y}$ opens a cavity with volume $\kappa_{\boldsymbol{a}}^{\boldsymbol{y}}>0$ in correspondence of each $\boldsymbol{a}\in A$. %More formally, we have  
%\begin{equation}
%	\label{eq:A-cavity}
%	\kappa_{\boldsymbol{a}}=\leb(\imt(\boldsymbol{y},\boldsymbol{a})) \quad \text{for all $\boldsymbol{a}\in A$,} \qquad A=C_{\boldsymbol{y}},
%\end{equation}
%where $\imt(\boldsymbol{y},\boldsymbol{a})$ and $C_{\boldsymbol{y}}$ denote the topological image of $\boldsymbol{a}$ under $\boldsymbol{y}$ and the set of cavitation points of $\boldsymbol{y}$, respectively, both recalled in Definition \ref{def:imt-point}. 

 We set 
\begin{equation*}
	\mathcal{Q}\coloneqq \left\{(A,\boldsymbol{y}):\:A\in\mathcal{A},\:\boldsymbol{y}\in\mathcal{Y}^A  \right\}.
\end{equation*}
Henceforth, we will tacitly assume that $\boldsymbol{d}$ is chosen in such a way that this class is nonempty.
For given parameters $\lambda_{\rm v},\lambda_{\rm p}\geq 0$, we define the  energy functional $\mathcal{E}\colon \mathcal{Q}\to [0,+\infty]$
as
\begin{equation}
	\label{eq:E}
	\mathcal{E}(A,\boldsymbol{y})\coloneqq \mathcal{W}(\boldsymbol{y})+\lambda_{\rm v} \mathcal{V}(A,\boldsymbol{y})+\lambda_{\rm p} \mathcal{P}(A,\boldsymbol{y}),
\end{equation} 
where  
\begin{align*}
	%\label{eq:W}
	\mathcal{W}(\boldsymbol{y})&\coloneqq \int_\Omega W(D\boldsymbol{y}(\boldsymbol{x}))\,\d\boldsymbol{x},\\
	 %\label{eq:V}
	 \mathcal{V}(A,\boldsymbol{y})&\coloneqq \sum_{\boldsymbol{a}\in A} \leb (\imt(\boldsymbol{y},\boldsymbol{a})),\\
	 %\label{eq:P}
	 \mathcal{P}(A,\boldsymbol{y})&\coloneqq \sum_{\boldsymbol{a}\in A} \per (\imt(\boldsymbol{y},\boldsymbol{a})).
\end{align*}
In the previous equation, $W\colon \RNN \to [0,+\infty]$ is the elastic energy density satisfying standard assumptions listed in Subsection \ref{subsec:main} below. The notation  $\imt(\boldsymbol{y},\boldsymbol{a})$ stands for the topological image of a point, see Definition~\ref{def:imt-point} below, which can be interpreted as the  cavity opened by $\boldsymbol{y}$ at the point $\boldsymbol{a}\in A$.    Note that the energy $\mathcal{E}$ possibly penalizes both the volume and the perimeter of the cavities.
% From \eqref{eq:A-cavity}, we note that the dependence of $\mathcal{V}$ and $\mathcal{P}$ on $A$ is fictitious as we can write
%\begin{equation*}
%	\mathcal{V}(A,\boldsymbol{y})=\sum_{\boldsymbol{a}\in C_{\boldsymbol{y}}} \per (\imt(\boldsymbol{y},\boldsymbol{a})), \qquad \mathcal{P}(A,\boldsymbol{y})= \sum_{\boldsymbol{a}\in C_{\boldsymbol{y}}} \leb (\imt(\boldsymbol{y},\boldsymbol{a})).
%\end{equation*}

 We will investigate the following regularization of the previously described model.  
Let $\varepsilon \in \big(0,\bar{\varepsilon}\big)$ with $\bar{\varepsilon}\coloneqq \dist(H;\partial \Omega)$. We denote by $\mathcal{A}_\varepsilon$ the collection of all  sets $A \subset H$ with $\mathscr{H}^0(A)\leq M$ such that $|\boldsymbol{a}-\boldsymbol{b}|\geq  3  \varepsilon$ for all  $\boldsymbol{a},\boldsymbol{b}\in A$ with $\boldsymbol{a}\neq \boldsymbol{b}$. 
For each $A\in\mathcal{A}_\varepsilon$, we {introduce the two sets} 
\begin{equation*}
	\Omega_\varepsilon^A\coloneqq \Omega \setminus \bigcup_{\boldsymbol{a}\in A} \closure{B}(\boldsymbol{a},\varepsilon), \qquad \widetilde{\Omega}^A_\varepsilon\coloneqq \Omega  \setminus  \bigcup_{\boldsymbol{a}\in A} B(\boldsymbol{a},\varepsilon). \EEE 
\end{equation*}

\begin{definition}[Approximating deformations\EEE]\label{def:deformation-perforated}
	Given $\varepsilon\in (0,\bar{\varepsilon})$ and $A\in\mathcal{A}_\varepsilon$, we define $\mathcal{Y}_\varepsilon^A$ as the class of all maps $\boldsymbol{y}\in W^{1,p}(\Omega_\varepsilon^A;\RN)$ satisfying the following properties:
	\begin{enumerate}[\label=\rm(i)]
		\item $\det D \boldsymbol{y}\in L^1(\Omega_\varepsilon^A)$ and $\det D \boldsymbol{y}(\boldsymbol{x})>0$ for almost all $\boldsymbol{x}\in \Omega^A_\varepsilon$;
		\item $\boldsymbol{y}$ satisfies condition {\rm (INV${}^\prime$)}; 
		\item  $\boldsymbol{y}^* \simeq \boldsymbol{d}^*$  on $\Gamma$;
		\item For all $\boldsymbol{a}\in A$, the following  hold:
		\begin{enumerate}[\indent\indent]
			\item[\rm (iv.1)] $\boldsymbol{y}^*\in W^{1,p}(S(\boldsymbol{a},\varepsilon);\RN)$;
			\item[\rm (iv.2)] $\boldsymbol{y}(\boldsymbol{x})\notin \imt(\boldsymbol{y},B(\boldsymbol{a},\varepsilon))$ for almost all $\boldsymbol{x}\in \Omega^A_\varepsilon$;
			\item[\rm (iv.3)] $\deg(\boldsymbol{y},B(\boldsymbol{a},\varepsilon),\boldsymbol{\xi})\geq 0$ for almost all $\boldsymbol{\xi}\in \RN$;
			\item[\rm (iv.4)]  $\boldsymbol{y}^*$  is $\haus$-almost everywhere injective on $S(\boldsymbol{a},\varepsilon)$;
			\item[\rm (iv.5)]  $\haus(S(\boldsymbol{a},\varepsilon)\setminus \domg(\boldsymbol{y}^*,S(\boldsymbol{a},\varepsilon)))=0$ and $\haus(\img(\boldsymbol{y},\Omega^A_\varepsilon)\cap \img(\boldsymbol{y}^*,S(\boldsymbol{a},\varepsilon)))=0$;
		\end{enumerate}
		\item There holds
		\begin{equation*}
			\Det^A_\varepsilon D \boldsymbol{y}=(\det D \boldsymbol{y}) \leb \quad \text{in $\mathcal{M}_{\rm b}(\widetilde{\Omega}^A_\varepsilon)$}.
		\end{equation*}
	\end{enumerate}
\end{definition}

Items (i) and (iii) have the same meaning as in Definition \ref{def:deformation}. 
Condition (INV${}^\prime$) mentioned in (ii) is a generalization of condition (INV) for perforated domains formulated in \cite{henao} and here recalled in Definition \ref{def:INV-prime}.  In (iv),  $\boldsymbol{y}^*$ indicates the trace of $\boldsymbol{y}$ on $S(\boldsymbol{a},\varepsilon)$.  The meaning of (iv.1) and (iv.4) is clear.  In (iv.2)--(iv.3),  $\imt(\boldsymbol{y},B(\boldsymbol{a},\varepsilon))$ and $\deg(\boldsymbol{y},B(\boldsymbol{a},\varepsilon),\cdot)$ denote the topological image of $B(\boldsymbol{a},\varepsilon)$ under $\boldsymbol{y}^*$ and the topological degree of $\boldsymbol{y}^*$ over $B(\boldsymbol{a},\varepsilon)$, respectively, as given in Definition \ref{def:top-deg}.   Roughly speaking, the two properties in (iv.2)--(iv.3) play the role of condition (INV) for the trace $\boldsymbol{y}^*$.   In (iv.5),  $\img(\boldsymbol{y},\Omega^A_\varepsilon)$ is the geometric image of $\Omega^A_\varepsilon$ under $\boldsymbol{y}$ given in Definition \ref{def:geom}, while  $\domg(\boldsymbol{y}^*,S(\boldsymbol{a},\varepsilon))$ and $\img(\boldsymbol{y}^*,S(\boldsymbol{a},\varepsilon))$ denote the geometric domain and image of $\boldsymbol{y}^*$ introduced in Definition \ref{def:geom-surface} below. 
  Eventually, $\Det^A_\varepsilon D\boldsymbol{y}$ in (v) is the extended distributional determinant introduced in Definition \ref{def:extended-distributional-determinant}. In particular, the identity in (v) ensures that $\boldsymbol{y}$ does not open any cavity. 
  
We define
\begin{equation*}
	\begin{split}
		\mathcal{Q}_\varepsilon\coloneqq \big\{ (A,\boldsymbol{y}):\:A\in\mathcal{A}_\varepsilon,\:
		\boldsymbol{y}\in\mathcal{Y}_\varepsilon^A \big\}.
	\end{split}
\end{equation*}
Under the assumption that $\mathcal{Q}$ is nonempty, it can be shown that, for almost every $0< \varepsilon \ll 1$, the same holds for $\mathcal{Q}_\varepsilon$. 
The regularized energy $\mathcal{E}_{\varepsilon}\colon \displaystyle \mathcal{Q}_\varepsilon \to [0,+\infty]$ is defined as
\begin{equation*}
	\mathcal{E}_\varepsilon(A,\boldsymbol{y})\coloneqq \mathcal{W}_\varepsilon(A,\boldsymbol{y})+\lambda_{\rm v} \mathcal{V}_\varepsilon(A,\boldsymbol{y})+\lambda_{\rm p} \mathcal{P}_\varepsilon(A,\boldsymbol{y}),
\end{equation*}
where 
	\begin{align*}
		%\label{eq:We}
		\mathcal{W}_\varepsilon(A,\boldsymbol{y})&\coloneqq 
		 \int_{\Omega_\varepsilon^A} W(D\boldsymbol{y}(\boldsymbol{x}))\,\d\boldsymbol{x}, \quad  \\ 
		 %\label{eq:Ve}
		  \mathcal{V}_\varepsilon(A,\boldsymbol{y})&\coloneqq \sum_{\boldsymbol{a}\in A} \leb(\imt(\boldsymbol{y},B(\boldsymbol{a},\varepsilon))),\\
		   %\label{eq:Pe}
		    \mathcal{P}_\varepsilon(A,\boldsymbol{y})&\coloneqq \sum_{\boldsymbol{a}\in A} \per(\imt(\boldsymbol{y},B(\boldsymbol{a},\varepsilon))).
	\end{align*}
We will see  in Lemma \ref{lem:degree-perimeter-imt} below that $\mathcal{V}_\varepsilon$ and $\mathcal{P}_\varepsilon$ can be equivalently expressed  as boundary integrals.
%\begin{align}
%	\label{eq:VV}
%	\mathcal{V}_\varepsilon(A,\boldsymbol{y})&=\frac{1}{N}\sum_{\boldsymbol{a}\in A} \int_{S(\boldsymbol{a},\varepsilon)} \boldsymbol{y}^*(\boldsymbol{x}) \cdot \left( (\cof D^{\rm t} \boldsymbol{y}^*(\boldsymbol{x}))\boldsymbol{\nu}_{B(\boldsymbol{a},\varepsilon)}(\boldsymbol{x})  \right)\,\d\haus(\boldsymbol{x}),\\
%	\label{eq:PP}
%	  	\mathcal{P}_\varepsilon(A,\boldsymbol{y})&= \sum_{\boldsymbol{a}\in A} \int_{S(\boldsymbol{a},\varepsilon)} |(\cof D^{\rm t} \boldsymbol{y}^*(\boldsymbol{x}))\boldsymbol{\nu}_{B(\boldsymbol{a},\varepsilon)}(\boldsymbol{x})|\,\d\haus(\boldsymbol{x}).
%\end{align}
%In \eqref{eq:VV}--\eqref{eq:PP}, we write $D^{\rm t}\boldsymbol{y}^*$ for the weak tangential gradient of $\boldsymbol{y}^*$. This notion  is revised in the final Appendix. Note that $D^{\rm t}\boldsymbol{y}^*$  is well defined in view of property (iv.1) in Definition \ref{def:deformation-perforated}.  
%Compared with \eqref{eq:Ve}--\eqref{eq:Pe}, the expressions in \eqref{eq:VV}--\eqref{eq:PP} seem  more approachable from the computational point of view as they do not involve  the topological degree.

\subsection{Main results}
\label{subsec:main}

 We now present  our main findings concerning the  approximation of singular minimizers. 
On the density $W$ we make the following standard assumptions:
\begin{enumerate}[(i)]
	\item \emph{Finiteness:} $W(\boldsymbol{F})<+\infty$ for all $\boldsymbol{F}\in\RNN_+$ and $W(\boldsymbol{F})=+\infty$ for all $\boldsymbol{F}\in\RNN\setminus \RNN_+$;
	\item \emph{Regularity:} $W$ is continuous in $\RNN$ and $W$ is of class $C^1$ in $\RNN_+$;
	\item \emph{Coercivity:} There exist a constant $c>0$ and a Borel function $g\colon (0,+\infty)\to [0,+\infty]$ with 
	\begin{equation}
		\label{eq:g}
		\lim_{\vartheta \to 0^+}g(\vartheta)=+\infty,\quad \lim_{\vartheta\to +\infty}\frac{g(\vartheta)}{\vartheta}=+\infty,
	\end{equation}
	such that
	\begin{equation}
		\label{eq:W-growth}
		W(\boldsymbol{F})\geq c|\boldsymbol{F}|^p+g(\det \boldsymbol{F})\quad \text{for all $\boldsymbol{F}\in\RNN_+$;}
	\end{equation}
	\item \emph{Polyconvexity:} $W$ is polyconvex, i.e., there exists a convex function $$\Phi\colon \prod_{i=1}^{N-1} \R^{\binom{N}{i}\times \binom{N}{i}}\times (0,+\infty)\to (0,+\infty)$$ such that 
	\begin{equation*}
		W(\boldsymbol{F})=\Phi(\mathbf{M}_{N-1}(\boldsymbol{F}),\det \boldsymbol{F}) \quad \text{for all $\boldsymbol{F}\in\RNN_+$,}
	\end{equation*}
	where $\mathbf{M}_{N-1}(\boldsymbol{F})$ collects all minors of $\boldsymbol{F}$ up to order $N-1$;
	\item \emph{Stress control:} There exist two constants $c_1 >0$ and $c_0\geq 0$ such that 
	\begin{equation*}
		|\boldsymbol{F}^\top DW(\boldsymbol{F})|\leq c_1 (W(\boldsymbol{F})+c_0) \quad \text{for all $\boldsymbol{F}\in\RNN_+$.}
	\end{equation*}
\end{enumerate}

In the following, given a sequence $(A_n)_n$ of subsets of  $\mathcal{A}_\varepsilon$    and $A  \in \mathcal{A}  $, we say that $A_n\to A$ up to ordering whenever there exists $m\in\N$ such that $\H^0(A_n)=m$ for all $n\in\N$, and we can write $A_n=\{ \boldsymbol{a}^n_1,\dots,\boldsymbol{a}^n_m \}$ for all $n\in \N$ and $A=\{ \boldsymbol{a}_1,\dots,\boldsymbol{a}_m \}$ in such a way that $\boldsymbol{a}^n_i\to \boldsymbol{a}_i$, as $n\to \infty$, for all $i=1,\dots,m$. Note that there might be repetitions among the points  $\boldsymbol{a}_1,\dots,\boldsymbol{a}_m$ when two sequences $(\boldsymbol{a}^n_i)_n$ and $(\boldsymbol{a}^n_j)_n$ have the same limit. The notion of  convergence up to ordering coincides with \cite[Definition 6.1]{moracorral}. Alternatively, one could consider the convergence induced by the Hausdorff distance as  in \cite{bresciani.friedrich}.

Our main results are formulated in the language of $\Gamma$-convergence. We begin with the compactness result.

\begin{theorem}[Compactness\EEE]
	\label{thm:equi-coercivity}
	Let $\varepsilon_n\to 0^+$ and let $((A_n,\boldsymbol{y}_n))_n$ be a sequence with $(A_n,\boldsymbol{y}_n)\in\mathcal{Q}_{\varepsilon_n}$  for all $n\in \N$ satisfying
	\begin{equation}
		\label{eq:equi-boundedness}
		\sup_{n\in\N} \mathcal{E}_{\varepsilon_n}(A_n,\boldsymbol{y}_n)<+\infty.
	\end{equation}
	Then, there exists $(A,\boldsymbol{y})\in \mathcal{Q}$ such that, up to subsequences, we have
	\begin{equation}
		\label{eq:convergence-A}
		\text{$A_n \to A$ \:\: up to ordering}
	\end{equation}
	and
	\begin{equation}
		\label{eq:convergence-y}
		\boldsymbol{y}_n \wk \boldsymbol{y} \:\:\text{in $W^{1,p}(\Omega_\delta^A;\RN)$,} \quad 
		\det D\boldsymbol{y}_n \wk \det D\boldsymbol{y}\:\: \text{in $L^1(\Omega_\delta^A)$ \quad for all $\delta>0$.}
	\end{equation} 
\end{theorem}
It is easy to see that, if \eqref{eq:convergence-A} holds, then $\Omega^A_\delta\subset \Omega^{A_n}_{\varepsilon_n}$ for $n\gg 1$ depending on $\delta$, so that \eqref{eq:convergence-y} implicitly  refers to the sequence $(\boldsymbol{y}_n\restr{\Omega^A_\delta})_{n\gg 1}$.

Next, we present the lower bound.

\begin{theorem}[Lower bound]
	\label{thm:lb}
Let  $\varepsilon_n\to 0^+$ and let $((A_n,\boldsymbol{y}_n))_n$ be a sequence with $(A_n,\boldsymbol{y}_n)\in\mathcal{Q}_{\varepsilon_n}$ for all $n\in \N$. Suppose that there exists  $(A,\boldsymbol{y})\in \mathcal{Q}$ for which \eqref{eq:convergence-A}--\eqref{eq:convergence-y} hold true. Then, we have
\begin{equation}\label{eq:lb}
	\liminf_{n\to \infty} \mathcal{E}_{\varepsilon_n}(A_n,\boldsymbol{y}_n)\geq \mathcal{E}(A,\boldsymbol{y}).
\end{equation}
\end{theorem}

Eventually, the third result establishes the optimality of the lower bound under an additional assumption of technical nature which, in a certain sense, constitutes a regularity requirement for admissible deformations.

\begin{theorem}[Optimality of the lower bound]\label{thm:ub}
Let  $\varepsilon_n\to 0^+$	and let $(A,\boldsymbol{y})\in \mathcal{Q}$. Additionally,  suppose that   
\begin{equation} \label{eq:conv-perimeter}
	\per \left(\imt(\boldsymbol{y},\boldsymbol{a})\right)=\liminf_{\substack{r\to 0^+ \\ B(\boldsymbol{a},r)\in\mathcal{U}_{\boldsymbol{y}}}} \per \left( \imt(\boldsymbol{y},B(\boldsymbol{a},r)) \right) \quad \text{for all $\boldsymbol{a}\in A$.}
\end{equation}  
 Then, there exists  a sequence $((A_n,\boldsymbol{y}_n))_n$ with $(A_n,\boldsymbol{y}_n)\in\mathcal{Q}_{\varepsilon_n}$ for all $n\in \N$ for which \eqref{eq:convergence-A}--\eqref{eq:convergence-y} hold true and we have
	\begin{equation}
	\label{eq:ub}
	\lim_{n\to \infty} \mathcal{E}_{\varepsilon_n}(A_n,\boldsymbol{y}_n)= \mathcal{E}(A,\boldsymbol{y}).
\end{equation}
\end{theorem}
In \eqref{eq:conv-perimeter}, the symbol $\mathcal{U}_{\boldsymbol{y}}$ stands for the class good domains introduced in Definition \ref{def:good} below. We mention that $B(\boldsymbol{a},r)\in\mathcal{U}_{\boldsymbol{y}}$ for  all $\boldsymbol{a}\in\Omega$ and \EEE almost all $0<r\ll 1$ by Lemma \ref{lem:good} below.  

We mention that when only the volume of the cavities is penalized, i.e., for $\lambda_{\rm p}=0$, the conclusion of Theorem \ref{thm:ub} still holds true without assuming \eqref{eq:conv-perimeter}. In that case, standard $\Gamma$-convergence arguments ensure  the  convergence of sequences of almost minimizers of $(\mathcal{E}_{\varepsilon_n})_n$ to minimizers of $\mathcal{E}$. In particular, for $\lambda_{\rm v}=\lambda_{\rm p}=0$, we recover \cite[Theorem 6.3]{henao.thesis}.

 Note that, because of \eqref{eq:conv-perimeter}, the existence of minimizers for the energy \eqref{eq:E} cannot be deduced from the abstract theory of $\Gamma$-convergence. To make up for this fact,   in Theorem \ref{thm:existence} below we give a direct proof for the existence of minimizers. \EEE  

Even though assumption \eqref{eq:conv-perimeter} seems well grounded,  so far, we are unfortunately   not  able to dispense with  it in proving Theorem \ref{thm:ub},  neither  to show that the limiting states identified in Theorem \ref{thm:equi-coercivity} satisfy this condition. In Section \ref{sec:example},  we illustrate the significance of assumption  \eqref{eq:conv-perimeter} by means of some examples.

\EEE

\section{Preliminaries}
\label{sec:prelim}

In this section, we recall some preliminary notions that will be employed within our analysis.
For convenience, we consider a bounded domain $G \subset \RN$  and  a  hypersurface $S\subset \RN$ of class $C^1$ for which a unit normal field $\boldsymbol{\nu}_S\colon S \to S(\boldsymbol{0},1)$ is assigned.  
 Later on, $G$ will be either $\Omega$ or $\Omega^A_\varepsilon$ for some $\varepsilon\in (0,\bar{\varepsilon})$ and $A\in\mathcal{A}_\varepsilon$,  while $S$ will be either $\partial U$ for some bounded domain $U\subset \RN$ of class $C^2$ with  $U\subset \subset \Omega$ or $\partial U \subset \Omega^A_\varepsilon$, or $S(\boldsymbol{a},\varepsilon)$ for some $\boldsymbol{a}\in A$.  

\subsection{Approximate differentiability} 
\label{subsec:approx-diff}
Recall that the  $N$-dimensional density of  a measurable  set $E\subset \RN$ at $\boldsymbol{x}_0\in\RN$ is defined as
\begin{equation*}
	\Theta^N(E,\boldsymbol{x}_0)\coloneqq \lim_{r\to 0^+} \frac{\leb(E\cap B(\boldsymbol{x}_0,r))}{\omega_N r^N}
\end{equation*}
whenever the limit exists. Here, $\omega_N\coloneqq \leb(B(\boldsymbol{0},1))$. 

 We refer to \cite[Subsection 3.1.2]{federer} for the notion of  approximate  differentiability and we define 
geometric domain and image as in \cite[Definition 2.4]{henao.moracorral.lusin}.

\begin{definition}[Geometric domain and image  for maps on domains]\label{def:geom}
Let $\boldsymbol{y}\in W^{1,p}(G;\RN)$ with $\det D \boldsymbol{y}(\boldsymbol{x})>0$ for almost all $\boldsymbol{x}\in G$. We define the geometric domain  of $\boldsymbol{y}$ as the set $\domg(\boldsymbol{y},G)$ of all points $\boldsymbol{x}\in G$ satisfying the following properties:
\begin{enumerate}[\rm (i)]
	\item $\boldsymbol{y}$ is approximately differentiable at $\boldsymbol{x}$ with approximate gradient $D\boldsymbol{y}(\boldsymbol{x})$ and $\det D \boldsymbol{y}(\boldsymbol{x})>0$;
	\item there exist a compact set $K\subset G$ with $\Theta^N(K,\boldsymbol{x})=1$  and a map $\boldsymbol{w}\in C^1(\RN;\RN)$ such that $\boldsymbol{y}\restr{K}=\boldsymbol{w}\restr{K}$ and $D\boldsymbol{y}\restr{K}=D\boldsymbol{w}\restr{K}$.
\end{enumerate}
Moreover, for every measurable set $E\subset G$, we set $\domg(\boldsymbol{y},E)\coloneqq E\cap \domg(\boldsymbol{y},G)$ and we define the geometric image of $E$ under $\boldsymbol{y}$ as $\img(\boldsymbol{y},E)\coloneqq \boldsymbol{y}\left( \domg(\boldsymbol{y},E)  \right)$.
\end{definition} 

 By standard properties of approximately differentiable maps \cite[Theorem 3, p.\ 217]{gms.cc1},  we see that $\domg(\boldsymbol{y},G)\cong G$. Also, if $\boldsymbol{y}$ is almost everywhere injective, then the restriction $\boldsymbol{y}\restr{\domg(\boldsymbol{y},G)}$ is actually injective \cite[Lemma 3]{henao.moracorral.fracture}.

The next result specializes Federer's area formula \cite[Theorem 1, p. 220]{gms.cc1} to our setting.

\begin{proposition}[Area formula for maps on domains]
	\label{prop:area-formula-domains}
Let $\boldsymbol{y}\in W^{1,p}(G;\RN)$ with $\det D\boldsymbol{y}(\boldsymbol{x})>0$ for almost all $\boldsymbol{x}\in G$. Suppose that $\boldsymbol{y}$ is almost everywhere injective. Then, for every measurable set $E\subset G$, we have
\begin{equation*}
	\label{eq:area-domain}
	\leb(\img(\boldsymbol{y},E))=\int_E \det D \boldsymbol{y}(\boldsymbol{x})\,\d\boldsymbol{x}.
\end{equation*}
Moreover, for every measurable function $\psi\colon \img(\boldsymbol{y},E)\to \R$, there holds
\begin{equation*}
	\label{eq:cov-domain}
	\int_{\img(\boldsymbol{y},E)} \psi(\boldsymbol{\xi})\,\d\boldsymbol{\xi}=\int_E \psi(\boldsymbol{y}(\boldsymbol{x}))\det D \boldsymbol{y}(\boldsymbol{x})\,\d\boldsymbol{x},
\end{equation*}
whenever one of the two integrals exists.
\end{proposition}

%From \eqref{eq:area-domain}, we see that  $\boldsymbol{y}\restr{\domg(\boldsymbol{y},G)}$ satisfies Lusin's property (N), that is, 
%\begin{equation*}
%\text{$\leb(\img(\boldsymbol{y},E))=0$\quad  whenever \quad $E\subset G$ satisfies $\leb(E)=0$.  	}
%\end{equation*}
%Also, by \cite[Remark 2.3(b)]{bresciani.friedrich.moracorral}, the map $\boldsymbol{y}$ satisfies Lusin's condition (N${}^{-1}$), that is,
%\begin{equation*}
%	\text{$\leb(\boldsymbol{y}^{-1}(F))=0$ \quad whenever \quad $F\subset \RN$ satisfies $\leb(F)=0$.}
%\end{equation*}
%In particular, this last property ensures the measurability of the composition $\psi\circ \boldsymbol{y}$ in \eqref{eq:cov-domain}. 

The concept of approximate tangential differentiability is recalled in the final Appendix. Given a map $\boldsymbol{y}\in W^{1,p}(S;\RN)$,  we denote  by $D^{\rm t}\boldsymbol{y}(\boldsymbol{x})$ and $J^{\rm t}\boldsymbol{y}(\boldsymbol{x})$ its weak tangential gradient and approximate tangential Jacobian, respectively.
 For convenience, similarly to Definition \ref{def:geom}, we introduce the geometric domain and image for maps on hypersurfaces.

\begin{definition}[Geometric domain and image for maps on hypersurfaces]\label{def:geom-surface}
Let $\boldsymbol{y}\in W^{1,p}(S;\RN)$ with $J^{\rm t}\boldsymbol{y}(\boldsymbol{x})>0$ for $\haus$-almost all $\boldsymbol{x}\in S$. We define the geometric domain of $\boldsymbol{y}$ as the set $\domg(\boldsymbol{y},S)$ of all points $\boldsymbol{x}\in S$ such that $\boldsymbol{y}$ is approximately tangentially differentiable at $\boldsymbol{x}$ with approximate tangential gradient $D^{\rm t}\boldsymbol{y}(\boldsymbol{x})$ and $J^{\rm t}\boldsymbol{y}(\boldsymbol{x})>0$. Moreover, for every $\haus$-measurable set $E\subset S$, we set $\domg(\boldsymbol{y},E)\coloneqq E\cap \domg(\boldsymbol{y},S)$ and we defined the geometric image of $E$ under $\boldsymbol{y}$ as $\img(\boldsymbol{y},E)\coloneqq \boldsymbol{y}(\domg(\boldsymbol{y},E))$.
\end{definition}

Similarly to before, see Proposition \ref{prop:watg} for details, we have that $\domg(\boldsymbol{y},S)\simeq S$.

The following version of the area formula for maps on hypersurfaces will be steadily employed in this work. The next result is a special case of Theorem \ref{thm:area-formula surface} whose proof is sketched in the Appendix. In particular, in this setting we observe that  $J^{\rm t}\boldsymbol{y}(\boldsymbol{x})=|(\cof D^{\rm t}\boldsymbol{y}(\boldsymbol{x}))\boldsymbol{\nu}_S(\boldsymbol{x})|$ for all $\boldsymbol{x}\in \domg(\boldsymbol{y},S)$.  We refer to   Lemma \ref{lem:j} for details.

\begin{proposition}[Area formula for  maps on hypersurfaces]
	\label{prop:cov-surface}
Let $\boldsymbol{y}\in W^{1,p}(S;\RN)$ with $J^{\rm t}\boldsymbol{y}(\boldsymbol{x})>0$ for $\haus$-almost all $\boldsymbol{x}\in S$. 	
Suppose that $\boldsymbol{y}$ is $\haus$-almost everywhere injective. 
Then, for every $\haus$-measurable set $E\subset S$, we have
\begin{equation}
	\label{eq:area-surface}
	\haus(\img(\boldsymbol{y},E))=\int_E |(\cof D^{\rm t}\boldsymbol{y}(\boldsymbol{x}))\boldsymbol{\nu}_S(\boldsymbol{x})|\,\d\haus(\boldsymbol{x}).
\end{equation}
Moreover, for every $\haus$-measurable function $\psi\colon \boldsymbol{y}(E)\to \R$, there holds
\begin{equation}
	\label{eq:cov-surface}
	\int_{\img(\boldsymbol{y},E)} {\psi}(\boldsymbol{\xi})\,\d\haus(\boldsymbol{\xi})=\int_E {\psi}(\boldsymbol{y}(\boldsymbol{x})) |(\cof D^{\rm t}\boldsymbol{y}(\boldsymbol{x}))\boldsymbol{\nu}_S(\boldsymbol{x})|\,\d\haus (\boldsymbol{x}),
\end{equation}
 whenever one of the two integrals is well defined. 
\end{proposition}

By passing, we mention that, as a consequence of the condition  $p>N-1$ and Morrey's embedding,
the maps $\boldsymbol{y}$ admits a representative $\overline{\boldsymbol{y}}\in C^0(S;\RN)$ which is $\haus$-almost everywhere tangentially
differentiable and satisfies Lusin's condition (N)  as given in  \eqref{eq:LusinN}. \EEE
We refer to Proposition \ref{prop:cont-repr} for a sketch of the proof of this result. 
 Thus, if \EEE  we apply Proposition \ref{prop:cov-surface} to the map $\overline{\boldsymbol{y}}$, then  the set $\img(\overline{\boldsymbol{y}},E)$ on the left-hand side of \eqref{eq:area-surface}--\eqref{eq:cov-surface} can be simply replaced by  $\overline{\boldsymbol{y}}(E)$. For example, when $E=S$, this yields 
\begin{equation}
	\label{eq:yS}
	\haus(\overline{\boldsymbol{y}}(S))\leq \|\cof D^{\rm t}\boldsymbol{y}\|_{L^1(S;\RNN)}\leq \|D^{\rm t}\boldsymbol{y}\|_{L^p(S;\RNN)}. 
\end{equation}

\EEE

\subsection{Topological degree}  
We recall the definition of  the  topological degree for Sobolev maps. We refer to \cite[Chapters 1--2]{fonseca.gangbo} for the case of continuous mappings.

\begin{definition}[Topological degree and topological image]
	\label{def:top-deg}
Let  $U\subset \RN$ be a bounded domain  of class  $C^1$   and let
  $\boldsymbol{y}\in W^{1,p}(\partial U;\RN)$. Denoting by $\overline{\boldsymbol{y}}\in C^0(\partial U;\RN)$ the continuous representative of $\boldsymbol{y}$, we define the topological degree $\deg(\boldsymbol{y},U,\boldsymbol{\xi})$ of $\boldsymbol{y}$ over $U$ at $\boldsymbol{\xi}\in\RN \setminus \overline{\boldsymbol{y}}(\partial U)$ as the topological degree of any continuous extension of $\overline{\boldsymbol{y}}$ to $\closure{U}$ over $U$ at $\boldsymbol{\xi}$. Moreover, we define the topological image of $U$ under $\boldsymbol{y}$ as 
\begin{equation*}
	\imt(\boldsymbol{y},U)\coloneqq \left\{ \boldsymbol{\xi}\in\RN\setminus \overline{\boldsymbol{y}}(\partial U):\:\deg(\boldsymbol{y},U,\boldsymbol{\xi})\neq 0  \right\}.
\end{equation*}
\end{definition}
In the previous definition, we take advantage of the fact that the topological degree is uniquely determined by boundary values \cite[Theorem 2.4]{fonseca.gangbo}.  
It is known that 
 $\imt(\boldsymbol{y},U)$ is an open and bounded set  with $\partial\,\imt(\boldsymbol{y},U)\subset\overline{\boldsymbol{y}}(\partial U)$ because of the continuity of the degree.  Note that, by \eqref{eq:yS},   $\haus(\overline{\boldsymbol{y}}(\partial U))<+\infty$  and, in turn,   $\leb(\partial\,\imt(\boldsymbol{y},U))\leq \leb(\overline{\boldsymbol{y}}(\partial U))=0$.  \EEE 
 
 The following formula for the distributional gradient of the degree will be instrumental for our analysis. We refer to \cite[Proposition 2.1]{mueller.spector.tang} for a proof.
 
 \begin{proposition}
 	\label{prop:degree-divergence}
 	Let $U\subset \RN$ be a bounded domain of class   $C^1$ \EEE  and $\boldsymbol{y}\in W^{1,p}(\partial U;\RN)$. Then
 	\begin{equation*}
 		\int_{\RN} \deg(\boldsymbol{y},U,\boldsymbol{\xi})\, \div \boldsymbol{\psi}(\boldsymbol{\xi})\,\d\boldsymbol{\xi}=\int_{\partial U} \boldsymbol{\psi}(\boldsymbol{y}(\boldsymbol{x}))\cdot (\cof  D^{\rm t}\boldsymbol{y}(\boldsymbol{x}))\boldsymbol{\nu}_U(\boldsymbol{x})\,\d\haus(\boldsymbol{x})
 	\end{equation*}	
 	for all $\boldsymbol{\psi}\in C^1(\RN;\RN)$.
 \end{proposition}

 In  what  follows, whenever $\boldsymbol{y}\in W^{1,p}(G;\RN)$ and $U\subset \RN$ is a bounded domain of class  $C^1$  with $\partial U\subset \closure{G}$ and  $\boldsymbol{y}^*\in W^{1,p}(\partial U;\RN)$,  we will simply write $\deg(\boldsymbol{y},U,\cdot)$ and $\imt(\boldsymbol{y},U)$ in place of $\deg(\boldsymbol{y}^*,U,\cdot)$ and $\imt(\boldsymbol{y}^*,U)$, respectively.

The next result has been established within the proof of \cite[Lemma 3.2, pp. 96--97]{sivaloganathan.spector}. 

\begin{lemma}[Positivity of the determinant]
	\label{lem:positivity}
	Let $(\boldsymbol{y}_n)_n$ be a sequence in $W^{1,p}(G;\RN)$ such that each $\boldsymbol{y}_n$ satisfies condition {\rm (INV)} and $\det D \boldsymbol{y}_n(\boldsymbol{x})>0$ for almost every $\boldsymbol{x}\in G$. Suppose that there exist $\boldsymbol{y}\in W^{1,p}(G;\RN)$ and $\vartheta\in L^1(G)$ with $\vartheta(\boldsymbol{x})>0$ for almost all $\boldsymbol{x}\in G$ such that
	\begin{equation*}
		\text{$\boldsymbol{y}_n\wk \boldsymbol{y}$ in $W^{1,p}(G;\RN)$,} \qquad \text{$\det D \boldsymbol{y}_n \wk  \vartheta  $ in $L^1(G)$.}
	\end{equation*}
	 Then, we have that $\det D \boldsymbol{y}(\boldsymbol{x})>0$ for almost all $\boldsymbol{x}\in G$.
\end{lemma}

\EEE

\subsection{Invertibility conditions}

In order to formulate the invertibility conditions, we need to  introduce some notation.   Given a bounded domain $U  \subset \RN$ of class $C^2$, the signed distance function $d_U\colon \R^N \to \R$ is defined by setting
\begin{equation*}
	d_U(\boldsymbol{x})\coloneqq \begin{cases}
		-\dist(\boldsymbol{x};\partial U) & \text{if $\boldsymbol{x}\in \closure{U},$} \\
		\dist(\boldsymbol{x};\partial U) & \text{if $\boldsymbol{x}\in \R^N \setminus \closure{U}.$}
	\end{cases}
\end{equation*}
Then,  there exists $ \bar{r}_U  >0$ such that the set $U_h\coloneqq \{\boldsymbol{x}\in \R^N:\:d_U(\boldsymbol{x})>h\}$ is a domain of class $C^2$ for every $h \in  (-\bar{r}_U,\bar{r}_U)$. We refer to \cite[Section~4]{ambrosio} for more details.

First, we state  condition (INV) introduced by M\"{u}ller and Spector  in \cite{mueller.spector}.

\begin{definition}[Condition INV]
	\label{def:INV}
	Let $\boldsymbol{y}\in W^{1,p}(G;\RN)$. The map $\boldsymbol{y}$ is termed to satisfy condition {\rm (INV)} if for every domain $U\subset \subset G$ of class $C^2$ and for almost all $h\in  (-\bar{r}_U,\bar{r}_U)$ such that $U_h\subset \subset G$ there hold
	\begin{enumerate}[\rm (i)]
		\item   $\boldsymbol{y}^*\in W^{1,p} (\partial U_h;\RN)$;  
		\item $\boldsymbol{y}(\boldsymbol{x})\in\imt(\boldsymbol{y},U_h)$ for almost all $\boldsymbol{x}\in U_h$;
		\item $\boldsymbol{y}(\boldsymbol{x})\notin \imt(\boldsymbol{y},U_h)$ for almost all $\boldsymbol{x}\in G\setminus U_h$.
	\end{enumerate}
\end{definition}

 It is worth to notice that, as a consequence of (ii) and the boundedness of the topological image, the map $\boldsymbol{y}$ in Definition \ref{def:INV} satisfies $\boldsymbol{y}\in L^\infty(U;\RN)$ for every domain $U\subset \subset G$ of class $C^2$. Equivalently, we find $\boldsymbol{y}\in L^\infty_{\rm loc}(G;\RN)$.

Next, we present the generalization of condition (INV) for  domains with holes due to Henao \cite{henao}. 

\begin{definition}[Condition INV${}^{\boldsymbol{\prime}}$]
	\label{def:INV-prime}
	Let $\boldsymbol{y}\in W^{1,p}(G;\RN)$. The map $\boldsymbol{y}$ is termed to satisfy condition {\rm (INV${}^{\prime}$)} if for every domain $U\subset \RN$ of class $C^2$ with $\partial U\subset G$ and for almost all $h\in  (-\bar{r}_U,\bar{r}_U) $ such that $\partial U_h\subset G$ there hold
	\begin{enumerate}[\rm (i)]
		\item  $ \boldsymbol{y}^*\in W^{1,p} (\partial U_h;\RN)$;
		\item $\boldsymbol{y}(\boldsymbol{x})\in\imt(\boldsymbol{y},U_h)$ for almost all $\boldsymbol{x}\in U_h\cap G$;
		\item $\boldsymbol{y}(\boldsymbol{x})\notin \imt(\boldsymbol{y},U_h)$ for almost all $\boldsymbol{x}\in G\setminus  U_h $.
	\end{enumerate}
\end{definition}

Roughly speaking, condition (INV${}^\prime$) requires the three properties defining condition (INV) to be satisfied for almost all domains $U\subset \RN$ of class $C^2$ such that $\partial U\subset G$ but not necessarily $U\subset G$. In particular, we immediately realize that  (INV${}^\prime$) implies  (INV). However, the converse does not hold and we refer to \cite[Section 6]{sivaloganathan.spector.tilakraj} for an example of a deformation satisfying  (INV)  but not (INV${}^\prime$).

 Similarly to before, we observe that the map $\boldsymbol{y}$ in Definition \ref{def:INV-prime} satisfies $\boldsymbol{y}\in L^\infty(U \cap G;\RN)$ for every bounded domain $U\subset \RN$ of class $C^2$ with $\partial U \subset G$.  In general, this property   is stronger than $\boldsymbol{y}\in L^\infty_{\rm loc}(G;\RN)$.  

We introduce the class of good domains.

\begin{definition}[Good domains]
	\label{def:good}
	Let $\boldsymbol{y}\in W^{1,p}(G;\RN)$ with $\det D \boldsymbol{y}(\boldsymbol{x})>0$ for almost all $\boldsymbol{x}\in G$. We define the class $\mathcal{U}_{\boldsymbol{y}}^\prime$ of good domains for $\boldsymbol{y}$ as the collection of all domains $U\subset \RN$ of class $C^2$ with $\partial U\subset G$ satisfying the following properties:
	\begin{enumerate}[\rm (i)]
		\item   $\boldsymbol{y}^*\in W^{1,p} (\partial U;\RN)$;
		\item $\partial U \simeq \domg(\boldsymbol{y},\partial U)$,  $\cof D \boldsymbol{y}\in L^1(\partial U;\RNN)$,  and $ D^{\rm t}\boldsymbol{y}^*\EEE\simeq D\boldsymbol{y}(\boldsymbol{I}-\boldsymbol{\nu}_U\otimes \boldsymbol{\nu}_U)$ on $\partial U$;
		\item $\boldsymbol{y}(\boldsymbol{x})\in\imt(\boldsymbol{y},U)$ for almost all $\boldsymbol{x}\in U \cap G$ and $\boldsymbol{y}(\boldsymbol{x})\notin \imt(\boldsymbol{y},U)$ for almost all $\boldsymbol{x}\in G \setminus U$; 
		\item There holds
		\begin{equation*}
			\lim_{ r\EEE\to 0^+} \frac{1}{ r }  \int_{-r}^{r} \left | \int_{\partial U_h} |\cof D\boldsymbol{y}|\,\d\haus - \int_{\partial U} |\cof D\boldsymbol{y}|\,\d\haus \right |\,\d h=0.
		\end{equation*}
		\item For all $\boldsymbol{\psi}\in \mathcal{D}(\RN;\RN)$, there holds
		\begin{equation*}
			\lim_{ r  \to 0^+} \frac{1}{ r }  \int_{-r}^r  \left | \int_{\partial U_h} \boldsymbol{\psi}\circ \boldsymbol{y}\cdot \left(\cof D \boldsymbol{y} \right)\boldsymbol{\nu}_{U_h}\,\d\haus - \int_{\partial U} \boldsymbol{\psi}\circ \boldsymbol{y}\cdot \left(\cof D \boldsymbol{y} \right)\boldsymbol{\nu}_{U}\,\d\haus \right | \,\d h=0
		\end{equation*}
	\end{enumerate}
	Moreover, we denote by $\mathcal{U}_{\boldsymbol{y}}$ the class of all $U\in\mathcal{U}_{\boldsymbol{y}}^\prime$ with $U\subset G$.
\end{definition}

Some comments are in order.  By (ii) and Definition \ref{def:geom}, the map $\boldsymbol{y}$ its approximately differentiable with approximate gradient $D\boldsymbol{y}(\boldsymbol{x})$ at  $\haus$-almost every $\boldsymbol{x}\in \partial U$, so that the two conditions in (ii) involving pointwise values of $D\boldsymbol{y}$ are meaningful. The same holds for the values of $D\boldsymbol{y}$ on $\partial U_h$ in (iv)--(v)   since, by the coarea formula, $\domg(\boldsymbol{y},G)\cong G$ entails $\domg(\boldsymbol{y},\partial U_h)\simeq \partial U_h$ for almost every   $h\in (-\bar{r}_U,\bar{r}_U)$. Additionally, the last condition in (ii) yields $\d^{\rm t}_{\boldsymbol{x}}\boldsymbol{y}^*=\d_{\boldsymbol{x}}\boldsymbol{y}\restr{\Tan_{\boldsymbol{x}}\partial U}$, where  $\d^{\rm t}_{\boldsymbol{x}}\boldsymbol{y}^*$, $\d_{\boldsymbol{x}}\boldsymbol{y}$, and $\Tan_{\boldsymbol{x}}\partial U$ are the approximate tangential differential of $\boldsymbol{y}^*$, the approximate differential of $\boldsymbol{y}$, and the tangent space to $\partial U$  at  $\boldsymbol{x}$, respectively, for $\haus$-almost every $\boldsymbol{x}\in \partial U$  (see the Appendix for more details).  Thus, recalling that $\det D \boldsymbol{y}(\boldsymbol{x})>0$ for $\haus$-almost every $\boldsymbol{x}\in \partial U$ because of (ii) and Definition \ref{def:geom}, we see that $\d_{\boldsymbol{x}}^{\rm t}\boldsymbol{y}^*$ is injective at any such point and, in turn, Lemma \ref{lem:j} and Binet's formula for cofactors yield the identity
\begin{equation}
	\label{eq:cof}
	\text{$J^{\rm t}\boldsymbol{y}^*\simeq|(\cof D^{\rm t}\boldsymbol{y}^*)\boldsymbol{\nu}_U|\simeq|(\cof D\boldsymbol{y})(\boldsymbol{\nu}_U \otimes \boldsymbol{\nu}_U)\boldsymbol{\nu}_U|=|(\cof D\boldsymbol{y})\boldsymbol{\nu}_U|$ on $\partial U$.}
\end{equation}

The next result can be proved as \cite[Lemma 2]{henao.moracorral.fracture} and \cite[Lemma 2.16]{henao.moracorral.lusin}. 

\begin{lemma}[Abundance of good domains]\label{lem:good}
Let $\boldsymbol{y}\in W^{1,p}(G;\RN)$ with $\det D \boldsymbol{y}(\boldsymbol{x})>0$ for almost all $\boldsymbol{x}\in G$ satisfy condition {\rm (INV${}^\prime$)}. Then, for each bounded domain $U\subset \RN$ of class $C^2$ with $\partial U\subset G$, we have that 
\begin{equation*}
	\text{$U_h \in\mathcal{U}_{\boldsymbol{y}}^\prime$ \quad for almost all $h\in  (-\bar{r}_U,\bar{r}_U)$.}
\end{equation*}
Additionally, if $U \subset G$, then we actually have $U_h \in\mathcal{U}_{\boldsymbol{y}}$ for almost all $h\in  (-\bar{r}_U,\bar{r}_U)$.	
\end{lemma}

The next proposition adapts \cite[Proposition 2.29]{bresciani.friedrich.moracorral} to the present setting. Its proof works 
 as the one in the reference and, hence, we omit it.

\begin{proposition}[Properties of topological image]
	\label{prop:topim}
Let $\boldsymbol{y}\in W^{1,p}(G;\RN)$ with $\det D \boldsymbol{y}(\boldsymbol{x})>0$ for almost all $\boldsymbol{x}\in G$ satisfy condition {\rm (INV${}^\prime$)}. Then:
\begin{enumerate}[\rm (i)]
		\item For all $U\in\mathcal{U}_{\boldsymbol{y}}^\prime$, we have
	\begin{equation*}
		\left (\closure{\imt(\boldsymbol{y},U)} \right )^\circ=\imt(\boldsymbol{y},U).
	\end{equation*}
	\item For all $U_1,U_2\in\mathcal{U}_{\boldsymbol{y}}^\prime$ with $\haus(\partial U_1\cap \partial U_2)=0$, we have
	\begin{align*}
		%\label{eq:imt-inclusion}
		\imt(\boldsymbol{y},U_1)\subset \imt(\boldsymbol{y},U_2) \quad &\text{whenever $U_1\subset U_2$,}\\
		%\label{eq:imt-intersection}
		\imt(\boldsymbol{y},U_1)\cap \imt(\boldsymbol{y},U_2)=\emptyset \quad &\text{whenever $U_1\cap U_2=\emptyset$.}
	\end{align*}
\end{enumerate}
\end{proposition}

 The following definition allows to rigorously describe the cavities.    

\begin{definition}[Topological image of a point]
	\label{def:imt-point}
	Let $\boldsymbol{y}\in W^{1,p}({G};\RN)$ with $\det D \boldsymbol{y}(\boldsymbol{x})>0$ for almost all $\boldsymbol{x}\in G$ satisfy condition {\rm (INV)}. The topological image of a point $\boldsymbol{a}\in G$ is defined as
	\begin{equation*}
		\imt(\boldsymbol{y},\boldsymbol{a})\coloneqq \bigcap \left\{\closure{\imt(\boldsymbol{y},B(\boldsymbol{a},r))}: r\in \left(0,\dist(\boldsymbol{a};\partial  G)\right) ,\:B(\boldsymbol{a},r)\in \mathcal{U}_{\boldsymbol{y}}   \right\}.
	\end{equation*}
	%The set of cavitation points of $\boldsymbol{y}$ is given by $C_{\boldsymbol{y}}\coloneqq \left\{ \boldsymbol{a}\in G:\:\leb(\imt(\boldsymbol{y},\boldsymbol{a}))>0  \right\}$.
\end{definition}
 
By arguing as in \cite[Remark 2.33(c)]{bresciani.friedrich.moracorral} with the aid of Proposition \ref{prop:topim}(ii), we see that   
\begin{equation}
	\label{eq:imtop-point-int}
	\imt(\boldsymbol{y},\boldsymbol{a})=\bigcap_{j\in\N} \closure{\imt(\boldsymbol{y},U_j)}
\end{equation}
for each sequence $(U_j)_j$ in $\mathcal{U}_{\boldsymbol{y}}$ with $\boldsymbol{a}\in U_j$ and $U_{j+1}\subset \subset U_{j}$.  This shows that  $\imt(\boldsymbol{y},U)$ is a nonempty compact set.  Moreover, as $\leb(\partial \,\imt(\boldsymbol{y},U_j))=0$ for all $j\in\N$, we find
\begin{equation}
	\label{eq:imtop-point-int2}
	\leb(\imt(\boldsymbol{y},\boldsymbol{a}))=\lim_{j\to \infty} \leb(\imt(\boldsymbol{y},U_j)).
\end{equation}

\subsection{Distributional determinant}
We recall the definition of the distributional determinant.

\begin{definition}[Distributional determinant]
	\label{def:distributional-determinant}
	Let $\boldsymbol{y}\in W^{1,p}(G;\RN)$ with $\det D \boldsymbol{y}(\boldsymbol{x})>0$ for almost all $\boldsymbol{x}\in G$ satisfy condition {\rm (INV)}. The distributional determinant  of $\boldsymbol{y}$ is the distribution $\Det D \boldsymbol{y}\in\mathcal{D}'(G)$ defined by setting
	\begin{equation*}
		(\Det D \boldsymbol{y})(\varphi)\coloneqq -\frac{1}{N}\int_\Omega \big((\adj D\boldsymbol{y})\boldsymbol{y}\big)\cdot D\varphi\,\d\boldsymbol{x}\quad \text{for all $\varphi\in \mathcal{D}(G)$.}
	\end{equation*}	
\end{definition}

In our analysis, we will deal with maps whose distributional determinant is represented by a given measure. For the ease of reference, we register a simple observation.

\begin{lemma}[Distributional determinant]
	\label{lem:Det}
	Let $\boldsymbol{y}\in W^{1,p}({G};\RN)$ with $\det D \boldsymbol{y}(\boldsymbol{x})>0$ for almost all $\boldsymbol{x}\in G$ satisfy condition {\rm (INV)}. Suppose that
	\begin{equation*}
		\Det D \boldsymbol{y}=(\det D \boldsymbol{y})\leb + \sum_{\boldsymbol{a}\in A} \kappa_{\boldsymbol{a}}^{\boldsymbol{y}} \delta_{\boldsymbol{a}} \quad \text{in $\mathcal{M}_{\rm b}(G)$}
	\end{equation*}
	for some  finite set $A\subset G$ and $\{ \kappa_{\boldsymbol{a}}^{\boldsymbol{y}}:\:\boldsymbol{a}\in A \}\subset (0,+\infty)$. Then:
	\begin{enumerate}[\rm (i)]
		\item $(\Det D \boldsymbol{y})(U)=\leb(\imt(\boldsymbol{y},U))$ for all $U\in\mathcal{U}_{\boldsymbol{y}}$;
		\item $A=\{\boldsymbol{a}\in G:\:\leb(\imt(\boldsymbol{y},\boldsymbol{a}))>0 \}$ and $\kappa_{\boldsymbol{a}}^{\boldsymbol{y}}=\leb(\imt(\boldsymbol{y},\boldsymbol{a}))$ for all $\boldsymbol{a}\in A$;
		\item 
		 $	\imt(\boldsymbol{y},U)\cong \img(\boldsymbol{y},U) \cup  \bigcup_{\boldsymbol{a}\in A \cap U} \imt(\boldsymbol{y},\boldsymbol{a})$ for all $U\in\mathcal{U}_{\boldsymbol{y}}$,
		where the two sets on the right-hand side are $\leb$-essentially disjoint. \EEE 
	\end{enumerate}
\end{lemma}
\begin{proof}
Claim (i) can be proved as in \cite[Corollary 9.2]{mueller.spector}.  It  is a consequence of \cite[Theorem 9.1]{mueller.spector} for which a corrected proof is given in \cite[Theorem A.1]{sivaloganathan.spector.tilakraj}. 

For the proof of (ii), let $\boldsymbol{a}\in G$ and consider a sequence $(r_j)_j$ with $r_j\to 0^+$, as $j\to \infty$, such that $B(\boldsymbol{a},r_j)\in\mathcal{U}_{\boldsymbol{y}}$ for all $j\in\N$. The existence of such a sequence in ensured by Lemma \ref{lem:good}.  For $\boldsymbol{a} \notin A$, we set  $\kappa_{\boldsymbol{a}}^{\boldsymbol{y}} = 0$.  
Then,  by \eqref{eq:imtop-point-int},  
\begin{equation*}
	\begin{split}
		\kappa_{\boldsymbol{a}}^{\boldsymbol{y}}&=(\Det D \boldsymbol{y})(\{ \boldsymbol{a} \})=(\Det D \boldsymbol{y}) \left( \textstyle  \bigcap_{j\in\N} B(\boldsymbol{a},r_j) \right)\\
		&=\lim_{j\to \infty} (\Det D \boldsymbol{y})(B(\boldsymbol{a},r_j))
		=\lim_{j\to \infty} \leb(\imt(\boldsymbol{y},B(\boldsymbol{a},r_j)))\\&= \leb \left( \textstyle  \bigcap_{j\in\N} \closure{\imt(\boldsymbol{y},B(\boldsymbol{a},r_j))}   \right)=\leb(\imt(\boldsymbol{y},\boldsymbol{a})).
	\end{split}
\end{equation*}
This shows the claim. 
In the previous equation, we used that  $(\Det D \boldsymbol{y})(B(\boldsymbol{a},r_j))=\leb(\imt(\boldsymbol{y},B(\boldsymbol{a},r_j)))$ in view of  claim (i),    $\leb(\partial\, \imt(\boldsymbol{y},B(\boldsymbol{a},r_j)))=0$, and, eventually,  $\imt(\boldsymbol{y},B(\boldsymbol{a},r_{j+1}))\subset \imt(\boldsymbol{y},B(\boldsymbol{a},r_j))$  by  Proposition \ref{prop:topim}(ii).

To prove (iii), we argue as in the proof of \cite[Theorem 3.2(iii)]{henao.moracorral.lusin}. First, recall that $\boldsymbol{y}$ is almost everywhere injective by \cite[Lemma 3.4]{mueller.spector}.
From claims (i)--(ii) and Proposition \ref{prop:area-formula-domains}, we see that
\begin{equation*}
	\leb(\imt(\boldsymbol{y},U))=\leb(\img(\boldsymbol{y},U))+\sum_{\boldsymbol{a}\in A \cap U} \leb(\imt(\boldsymbol{y},\boldsymbol{a})).
\end{equation*} 
Note that $\img(\boldsymbol{y},U)\subset \closure{\imt(\boldsymbol{y},U)}\cong\imt(\boldsymbol{y},U)$ by \cite[Lemma 2.27]{bresciani.friedrich.moracorral}. Also, the same result entails
\begin{equation*}
	\imt(\boldsymbol{y},\boldsymbol{a})\subset \closure{\imt(\boldsymbol{y},B(\boldsymbol{a},r))}\cong \imt(\boldsymbol{y},B(\boldsymbol{a},r))\subset \RN \setminus \img(\boldsymbol{y},\Omega\setminus B(\boldsymbol{a},r))
\end{equation*}   
for all $\boldsymbol{a}\in A\cap U$ and $r>0$ with $B(\boldsymbol{a},r)\in\mathcal{U}_{\boldsymbol{y}}$. Hence, 
\begin{equation*}
	\leb \left( \imt(\boldsymbol{y},\boldsymbol{a}) \cap \img(\boldsymbol{y},\Omega) \right) \leq \leb \left( \img(\boldsymbol{y},\Omega) \setminus \img(\boldsymbol{y},\Omega\setminus B(\boldsymbol{a},r)) \right)=\leb(\img(\boldsymbol{y},B(\boldsymbol{a},r))).
\end{equation*}
Letting $r\to 0^+$ along a suitable sequence, with the aid of Proposition \ref{prop:area-formula-domains}, we  find $\leb(\imt(\boldsymbol{y},\boldsymbol{a}) \cap \img(\boldsymbol{y},\Omega))=0$. Therefore, the conclusion follows.
\end{proof}

The next lemma can be proved by arguing as in the proof of \cite[Lemma 3.2, pp. 94--95]{sivaloganathan.spector}. 

\begin{lemma}[Weak*-continuity of distributional determinants]
	\label{lem:weakstar}
Let $(\boldsymbol{y}_n)_n$ be a sequence in $W^{1,p}(G;\RN)$ such that each $\boldsymbol{y}_n$ satisfies condition {\rm (INV)},  $\det D \boldsymbol{y}_n(\boldsymbol{x})>0$ for almost all $\boldsymbol{x}\in G$, and  $\Det D \boldsymbol{y}_n\in\mathcal{M}_{\rm b}(G)$. Suppose that there exists $\boldsymbol{y}\in W^{1,p}(G;\RN)$ such that 
\begin{equation*}
	\text{$\boldsymbol{y}_n \wk \boldsymbol{y}$ in $W^{1,p}(G;\RN)$.}
\end{equation*}
Then, $\Det D \boldsymbol{y}\in\mathcal{M}_{\rm b}(G)$ and  we have
\begin{equation*}
	\text{$\Det D \boldsymbol{y}_n \wks \Det D \boldsymbol{y}$ in $\mathcal{M}_{\rm b}(G)$.}
\end{equation*}
\end{lemma}

\section{Existence of minimizers}
\label{sec:existence}

In this section, we prove the existence of minimizers of the functional $\mathcal{E}$ in \eqref{eq:E}.  The result extends   \cite[Theorem 6.3]{henao.thesis}   by possibly accounting for the volume of the cavities (as done in \cite[Theorem 4.2]{sivaloganathan.spector} for fixed flaw points) as well as their perimeter. As we mentioned in  Subsection \ref{subsec:overview}, this result can be easily inferred from the theory in \cite{henao.moracorral.lusin}. Here, we provide an alternative and more self-contained proof. 

\begin{theorem}[Existence of minimizers]\label{thm:existence}
Suppose that $\mathcal{Q}$ is nonempty. Then, the functional $\mathcal{E}$ in \eqref{eq:E} admits minimizers within the class $\mathcal{Q}$.
\end{theorem} 

The proof of the theorem relies on the following convergence result.

\begin{proposition}[Convergence]
	\label{prop:convergence}
Let $( (A_n,\boldsymbol{y}_n))_n$ be a sequence in $\mathcal{Q}$. Suppose that there exist $B\in\mathcal{A}$, $\boldsymbol{y}\in W^{1,p}(\Omega;\RN)$, and $\vartheta\in L^1(\Omega)$ with $\vartheta(\boldsymbol{x})>0$ for almost all $\boldsymbol{x}\in \Omega$ such that
\begin{equation*}
	\text{$A_n\to B$ up to ordering,} \quad \text{$\boldsymbol{y}_n\wk \boldsymbol{y}$ in $W^{1,p}(\Omega;\RN)$,} \quad \text{$\det D \boldsymbol{y}_n\to \vartheta$ in $L^1(\Omega)$.}
\end{equation*}
Then, $\vartheta\cong\det D \boldsymbol{y}$ and there exists $A\subset B$ such that $\boldsymbol{y}\in \mathcal{Y}^{A}$, i.e., $(A,\boldsymbol{y})\in\mathcal{Q}$. Moreover, 
\begin{equation}
	\label{eq:cimt}
	\text{$\sum_{\boldsymbol{a}\in A_n} \chi_{\imt(\boldsymbol{y}_n,\boldsymbol{a})} \to \sum_{\boldsymbol{a}\in A}  \chi_{\imt(\boldsymbol{y},\boldsymbol{a})} $ \quad  in $L^1(\RN)$.}
\end{equation}
\end{proposition}
\begin{proof}
We divide the proof into two steps.

\emph{Step 1.} We begin by identifying the set $A$ and showing that $(A,\boldsymbol{y})\in\mathcal{Q}$. 	
By assumption, there exists $m\in\N$ with $\H^0(A_n)=m$ for all $n\in\N$ and  $\H^0(B)\leq m$ for which we can write $A_n=\{ \boldsymbol{a}^n_1,\dots,\boldsymbol{a}^n_m \}$ for all $n\in\N$ and $B=\{ \boldsymbol{a}_1,\dots,\boldsymbol{a}_m \}$ in such a way that $\boldsymbol{a}^n_i\to \boldsymbol{a}_i$, as $n\to \infty$, for all $i=1,\dots,m$. By  Lemma \ref{lem:positivity}, we know that $\det D \boldsymbol{y}(\boldsymbol{x})>0$ for almost all $\boldsymbol{x}\in\Omega$. Moreover,  \cite[Lemma 3.3]{mueller.spector} ensures that $\boldsymbol{y}$ satisfies (INV) and, in turn,  that  $\boldsymbol{y}$ is  almost everywhere injective,  see  \cite[Lemma 3.4]{mueller.spector}.

By applying  Lemma \ref{lem:weakstar}, we find that  $\Det D \boldsymbol{y}\in\mathcal{M}_{\rm b}(\Omega)$  and we also have 
\begin{equation}
	\label{eq:dy}
	\text{$\Det D \boldsymbol{y}_n \wks \Det D \boldsymbol{y}$ in $\mathcal{M}_{\rm b}(\Omega)$.}
\end{equation}
In particular, the sequence $(\Det D \boldsymbol{y}_n)_n$ is bounded in $\mathcal{M}_{\rm b}(\Omega)$. By assumption  (iv)  in Definition \ref{def:deformation},  each measure $\Det D \boldsymbol{y}_n$ admits a representation of the form 
\begin{equation}
	\label{eq:dn}
	\Det D \boldsymbol{y}_n=(\det D \boldsymbol{y}_n)\leb + \sum_{i=1}^{m} \kappa_i^n \delta_{\boldsymbol{a}^n_i} \quad \text{in $\mathcal{M}_{\rm b}(\Omega)$}
\end{equation}
for some $\kappa_1^n,\dots,\kappa_m^n>0$. 
As the boundedness of $(\Det D \boldsymbol{y}_n)_n$ entails that of each sequence $(\kappa^n_i)_n$, we find $\kappa_1,\dots,\kappa_m\geq 0$ such that, up to subsequences, $\kappa^n_i\to \kappa_i$, as $n\to \infty$, for all $i=1,\dots,m$.
Passing to the limit, as $n\to \infty$, in  \eqref{eq:dn}  and combining the result with \eqref{eq:dy}, we deduce
\begin{equation*}
	\label{eq:ddy}
	\Det D \boldsymbol{y}=\vartheta \leb + \sum_{i=1}^m \kappa_i\delta_{\boldsymbol{a}_i} \quad \text{in $\mathcal{M}_{\rm b}(\Omega)$.}
\end{equation*}
 Thus, \cite[Lemma 8.1]{mueller.spector} yields $\vartheta  \cong  \det D \boldsymbol{y}$. 

 For convenience, we set $I\coloneqq \{1,\dots,m\}$. Then, we define  $A\coloneqq \{ \boldsymbol{a}_i:\:i\in I, \:\kappa_i>0  \}$.
 Also, for each $\boldsymbol{a}\in A$, we set $I_{\boldsymbol{a}}\coloneqq \{i\in I:\:\boldsymbol{a}_i=\boldsymbol{a}\}$ and $\kappa_{\boldsymbol{a}}\coloneqq \sum_{i\in I_{\boldsymbol{a}}} \kappa_i$. Note that we might have $\kappa_i=0$ for some $i\in I_{\boldsymbol{a}}$. 
  With this notation, the previous equation becomes
\begin{equation*}
	\Det D \boldsymbol{y}=(\det D \boldsymbol{y})\leb + \sum_{\boldsymbol{a}\in {A}} \kappa_{\boldsymbol{a}} \delta_{\boldsymbol{a}}.
\end{equation*}
 As clearly also $\boldsymbol{y}^*\simeq \boldsymbol{d}^*$ on $\Gamma$ holds, we conclude   $(A,\boldsymbol{y})\in\mathcal{Q}$.

\emph{Step 2.} We are left to prove \eqref{eq:cimt}. Setting $J\coloneqq \{ i\in I:\:\boldsymbol{a}_i\notin A \}$, we have $\kappa_i=0$ for all $i\in J$. Moreover, we can write $I=J \cup \bigcup_{\boldsymbol{a}\in A} I_{\boldsymbol{a}}$ and, in turn, 
\begin{equation}
	\label{eq:bhbh}
	\sum_{\boldsymbol{a}\in A_n} \chi_{\imt(\boldsymbol{y}_n,\boldsymbol{a})}=\sum_{i\in J} \chi_{\imt(\boldsymbol{y}_n,\boldsymbol{a}^n_i)}+ \sum_{\boldsymbol{a}\in A} \sum_{i\in I_{\boldsymbol{a}}} \chi_{\imt(\boldsymbol{y}_n,\boldsymbol{a}^n_i)}. 
\end{equation} 
For $i\in J$, considering \eqref{eq:dn} and using Lemma \ref{lem:Det}(ii), we have  $\leb(\imt(\boldsymbol{y}_n,\boldsymbol{a}^n_i))=\kappa^n_i \to \kappa_i=0$, or, equivalently, $\chi_{\imt(\boldsymbol{y}_n,\boldsymbol{a}^n_i)}\to 0$ in $L^1(\RN)$, as $n\to \infty$. If we prove that
\begin{equation}
	\label{eq:s}
	\text{$\sum_{i\in I_{\boldsymbol{a}}}\chi_{\imt(\boldsymbol{y}_n,\boldsymbol{a}^n_i)}\to \chi_{\imt(\boldsymbol{y},\boldsymbol{a})}$ in $L^1(\RN)$, \quad as $n\to \infty$, \quad for all $\boldsymbol{a}\in A$,}
\end{equation}
then, letting $n\to \infty$ in \eqref{eq:bhbh}, we establish  \eqref{eq:cimt}. 

To prove \eqref{eq:s}, fix $\boldsymbol{a}\in A$. By applying \cite[Lemma 2.9]{mueller.spector} in combination with a diagonal argument, we select a sequence  $(r_l)_l$  of positive radii with $r_l\to 0^+$, as $l\to \infty$, and a not relabeled subsequence of $(\boldsymbol{y}_{n})_n$ such that, for all $l\in\N$,
\begin{equation}
	\label{eq:sb}
	B(\boldsymbol{a},r_l)\in \mathcal{U}_{\boldsymbol{y}} \cap \bigcap_{n \in \N} \mathcal{U}_{\boldsymbol{y}_{n}}, \qquad  \text{$\boldsymbol{y}^*_{n}\wk \boldsymbol{y}^*$ in $W^{1,p}(S(\boldsymbol{a},r_l);\RN)$, \quad as $n\to \infty$.}
\end{equation}
Given $l\in\N$, consider $n\gg 1$ depending on $l$  such that $\boldsymbol{a}^{n}_i \in B(\boldsymbol{a},r_l)$ for all $i\in I_{\boldsymbol{a}}$.
We write
\begin{equation}
	\label{eq:ea}
	\begin{split}
		\left 	\|\sum_{i\in I_{\boldsymbol{a}}}\chi_{\imt(\boldsymbol{y}_n,\boldsymbol{a}^n_i)}-\chi_{\imt(\boldsymbol{y},\boldsymbol{a})}\right \|_{L^1(\RN)}&\leq \left  \|\sum_{i\in I_{\boldsymbol{a}}}\chi_{\imt(\boldsymbol{y}_n,\boldsymbol{a}^n_i)}-\chi_{\imt(\boldsymbol{y}_n,B(\boldsymbol{a},r_l))} \right \|_{L^1(\RN)}\\
		&+\| \chi_{\imt(\boldsymbol{y}_n,B(\boldsymbol{a},r_l))} - \chi_{\imt(\boldsymbol{y},B(\boldsymbol{a},r_l))}\|_{L^1(\RN)}\\
		&+\|\chi_{\imt(\boldsymbol{y},B(\boldsymbol{a},r_l))} - \chi_{\imt(\boldsymbol{y},\boldsymbol{a}))}\|_{L^1(\RN)}.
	\end{split}
\end{equation}
For the first term on the right-hand side of \eqref{eq:ea}, note that
the family of sets $\{ \imt(\boldsymbol{y}_{n},\boldsymbol{a}^{n}_i):\:i\in I_{\boldsymbol{a}} \}$ is $\haus$-essentially disjoint in view of \cite[Lemma 2.34]{bresciani.friedrich.moracorral} (cf.~\cite[Lemma 7.6]{mueller.spector}). Hence,
\begin{equation*}
	\sum_{i\in I_{\boldsymbol{a}}}\chi_{\imt(\boldsymbol{y}_{n},\boldsymbol{a}^{n}_i)}\simeq  \chi_{ \bigcup_{ i \in I_{\boldsymbol{a}}}\imt(\boldsymbol{y}_{n},\boldsymbol{a}_i^{n})}. 
\end{equation*}
Thus, by Lemma \ref{lem:Det}(iii),
\begin{equation*}
	\begin{split}
		\left  \|\sum_{i\in I_{\boldsymbol{a}}}\chi_{\imt(\boldsymbol{y}_n,\boldsymbol{a}^n_i)}-\chi_{\imt(\boldsymbol{y}_n,B(\boldsymbol{a},r_l))} \right \|_{L^1(\RN)}&=\leb \left( \imt(\boldsymbol{y}_n,B(\boldsymbol{a},r_l)) \setminus \bigcup_{ i \in I_{\boldsymbol{a}}}  \imt(\boldsymbol{y}_n,\boldsymbol{a}^n_i)  \right)\\
		&=  \leb(\img(\boldsymbol{y},B(\boldsymbol{a},r_l)))  =\int_{B(\boldsymbol{a},r_l)} \det D \boldsymbol{y}_n(\boldsymbol{x})\,\d \boldsymbol{x},
	\end{split}
\end{equation*} 
where in the last line we applied Proposition \ref{prop:area-formula-domains}. The second summand on the right-hand side of \eqref{eq:ea} goes to zero, as $n\to \infty$, because of \eqref{eq:sb} and \cite[Lemma 2.19]{bresciani.friedrich.moracorral}. For the third one, we have
\begin{equation*}
	\|\chi_{\imt(\boldsymbol{y},B(\boldsymbol{a},r_l))} - \chi_{\imt(\boldsymbol{y},\boldsymbol{a}))}\|_{L^1(\RN)}=\leb(\imt(\boldsymbol{y},B(\boldsymbol{a},r_l)))-\leb(\imt(\boldsymbol{y},\boldsymbol{a})).
\end{equation*}
Hence, letting $n\to \infty$ in \eqref{eq:ea}, we find
\begin{equation*}
	 \begin{split}
	 	\limsup_{n\to \infty} \left 	\|\sum_{i\in I_{\boldsymbol{a}}}\chi_{\imt(\boldsymbol{y}_n,\boldsymbol{a}^n_i)}-\chi_{\imt(\boldsymbol{y},\boldsymbol{a})}\right \|_{L^1(\RN)}&\leq \int_{B(\boldsymbol{a},r_l)} \det D \boldsymbol{y}(\boldsymbol{x})\,\d\boldsymbol{x}\\
	 	&+ \leb(\imt(\boldsymbol{y},B(\boldsymbol{a},r_l)))-\leb(\imt(\boldsymbol{y},\boldsymbol{a})).
	 \end{split}
\end{equation*}
As $l\to \infty$, recalling \eqref{eq:imtop-point-int2}, the right-hand side goes to zero. This concludes the proof of \eqref{eq:s}.
\end{proof}

With Proposition \ref{prop:convergence} at hand, the proof of Theorem \ref{thm:existence} is achieved by a standard application of the direct method.

\begin{proof}[Proof of Theorem \ref{thm:existence}]
Without loss of generality,  we   assume $\alpha\coloneqq \inf_{\mathcal{Q}} \mathcal{F}\in [0,+\infty)$. 
Let $((A_n,\boldsymbol{y}_n))_n$ be a sequence in $\mathcal{Q}$ which is minimizing  for $\mathcal{E}$, i.e., such that $\mathcal{E}(A_n,\boldsymbol{y}_n)\to \alpha$, as $n\to \infty$. 
In view of the confinement given by $H$, there exists $B\in\mathcal{A}$ such that $A_n\to B$ up to ordering for a not relabeled subsequence. 
From the boundedness of $(\mathcal{E}(A_n,\boldsymbol{y}_n))_n$, using assumption \eqref{eq:W-growth}, item (iii) of Definition \ref{def:deformation}, and the Poincar\'{e} inequality with trace term, we obtain 
\begin{equation*}
	\sup_{n\in\N} \left\{ \| \boldsymbol{y}_n \EEE \|_{W^{1,p}(\Omega;\RN)}+\|g(\det D \boldsymbol{y}_n)\|_{L^1(\Omega)}  \right\}<+\infty.
\end{equation*}	
Thus,   in view of  \eqref{eq:g},  we find $\boldsymbol{y}\in W^{1,p}(\Omega;\RN)$ and  $\vartheta\in L^1(\Omega)$, such that, up to subsequences, $\boldsymbol{y}_n \wk \boldsymbol{y}$ in $W^{1,p}(\Omega;\RN)$ and $\det D \boldsymbol{y}_n \wk \vartheta$ in $L^1(\Omega)$. In particular, a standard contradiction argument based on \eqref{eq:g}--\eqref{eq:W-growth} shows that $\vartheta(\boldsymbol{x})>0$ for almost all $\boldsymbol{x}\in \Omega$. Therefore, by applying Proposition \ref{prop:convergence}, we find $A\subset B$ such that $(A,\boldsymbol{y})\in\mathcal{Q}$ and we deduce $\vartheta\cong \det D \boldsymbol{y}$ as well as \eqref{eq:cimt}. 

We claim that $(A,\boldsymbol{y})$ is a minimizer of $\mathcal{E}$. First, owing to the weak continuity of minors and the polyconvexity of $W$, we find
\begin{equation}
	\label{eq:m1}
	\liminf_{n\to \infty} \mathcal{W}(\boldsymbol{y}_n)\geq \mathcal{W}(\boldsymbol{y}). 
\end{equation}
Second, observe that each family of sets $\{ \imt(\boldsymbol{y}_n,\boldsymbol{a}):\:\boldsymbol{a}\in A_n  \}$ for  $n\in\N$ as well as $\{\imt(\boldsymbol{y},\boldsymbol{a}):\:\boldsymbol{a}\in A  \}$ is $\haus$-essentially disjoint because of \cite[Lemma 2.34]{bresciani.friedrich.moracorral},  see also \cite[Lemma 7.6]{mueller.spector}.   Thus, \eqref{eq:cimt} gives
\begin{equation}
	\label{eq:m2}
	\begin{split}
		\lim_{n\to \infty}\sum_{\boldsymbol{a}\in A_n} \leb(\imt(\boldsymbol{y}_n,\boldsymbol{a}))&= \lim_{n\to \infty} \leb \left( \bigcup_{\boldsymbol{a}\in A_n} \imt(\boldsymbol{y}_n,\boldsymbol{a})  \right)\\
		&= \leb \left(\bigcup_{\boldsymbol{a}\in A} \imt(\boldsymbol{y},\boldsymbol{a}) \right)=\sum_{\boldsymbol{a}\in {A}} \leb \left( \imt(\boldsymbol{y},\boldsymbol{a})\right),
	\end{split}
\end{equation}
and, thanks to the lower semicontinuity of the perimeter, 
\begin{equation}
	\label{eq:m3}
	\begin{split}
		\liminf_{n\to \infty} \sum_{\boldsymbol{a}\in {A}_n} \per \left( \imt(\boldsymbol{y}_n,\boldsymbol{a}) \right)&=\liminf_{n\to \infty} \per \left( \bigcup_{\boldsymbol{a}\in {A}_n} {\imt(\boldsymbol{y}_n,\boldsymbol{a})}  \right)\\
		&\geq \per \left( \bigcup_{\boldsymbol{a}\in {A}}  {\imt(\boldsymbol{y},\boldsymbol{a})} \right)
		=\sum_{\boldsymbol{a}\in {A}}  \per \left( \imt(\boldsymbol{y},\boldsymbol{a})\right).
	\end{split}
\end{equation}
Eventually, the combination of \eqref{eq:m1}--\eqref{eq:m3} leads to
\begin{equation*}
	\alpha=\liminf_{n\to \infty} \mathcal{E}(A_n,\boldsymbol{y}_n)\geq \mathcal{E}(A,\boldsymbol{y})
\end{equation*}
which yields the claim.
\end{proof}

\EEE 

\section{Deformations on perforated domains}
\label{sec:perforated}

Throughout  this section,  let $\varepsilon \in (0,\bar{\varepsilon})$ and $A\in \mathcal{A}_\varepsilon$ be fixed.

\subsection{Properties of geometric and topological images}
We investigate the relationship between geometric and topological images in the case of deformations on perforated domains.
The first result generalizes  \cite[Lemma 2.27]{bresciani.friedrich.moracorral} to this situation.

\begin{lemma}[Geometric image and topological image]
	\label{lem:img-vs-imt}
Let  $\boldsymbol{y}\in W^{1,p}(\Omega^A_\varepsilon;\RN)$ satisfy property {\rm (iv)} in the Definition \ref{def:deformation-perforated}. Then, the following holds:
\begin{enumerate}[\rm (i)]
	\item For every $U\in\mathcal{U}_{\boldsymbol{y}}^\prime$, denoting by $\overline{\boldsymbol{y}}$ the continuous representative of the trace of  $\boldsymbol{y}$ on $\partial U$, we have
	\begin{equation*}
		\img(\boldsymbol{y},U\cap \Omega^A_\varepsilon)\subset \closure{\imt(\boldsymbol{y},U)}, \qquad
		\img(\boldsymbol{y},\Omega^A_\varepsilon \setminus {U})\subset \RN \setminus \imt(\boldsymbol{y},U),
	\end{equation*}
	and
	\begin{equation*}
		\partial\, \imt(\boldsymbol{y},U)=\overline{\boldsymbol{y}}(\partial U).
	\end{equation*}
	\item For all $\boldsymbol{a}\in A$,  we have
	\begin{equation*}
		\img(\boldsymbol{y},\Omega^A_\varepsilon)\subset \RN \setminus \imt(\boldsymbol{y},B(\boldsymbol{a},\varepsilon)).
	\end{equation*}
\end{enumerate}
\end{lemma}
\begin{proof}
The proof of (i)  is identical to the one in \cite[Lemma 2.27]{bresciani.friedrich.moracorral}. For (ii), we adopt the same argument  for proving the first inclusion in (i). We repeat it for the convenience of the reader. 
Let $\boldsymbol{a}\in A$. By assumption, $\boldsymbol{y}(\boldsymbol{x})\notin \imt(\boldsymbol{y},B(\boldsymbol{a},\varepsilon))$ for almost all $\boldsymbol{x}\in \Omega^A_\varepsilon$. Thus, setting $E\coloneqq \{ \boldsymbol{x}\in \domg(\boldsymbol{y},\Omega^A_\varepsilon):\: \boldsymbol{y}(\boldsymbol{x})\notin \imt(\boldsymbol{y},B(\boldsymbol{a},\varepsilon))  \}$, we find $E\cong \Omega^A_\varepsilon$. Fix $\boldsymbol{x}\in \domg(\boldsymbol{y},\Omega^A_\varepsilon)$. Then, $\Theta^N(E, \boldsymbol{x}  )=\Theta^N(\Omega^A_\varepsilon,\boldsymbol{x})=1$ and, by \cite[Lemma 1]{henao.moracorral.fracture}, we find $\Theta^N(\boldsymbol{y}(E),\boldsymbol{y}(\boldsymbol{x}))=1$. This yields $\Theta^N(\RN \setminus \imt(\boldsymbol{y},B(\boldsymbol{a},\varepsilon)),\boldsymbol{y}(\boldsymbol{x}))=1$ and $\Theta^N( \imt(\boldsymbol{y},B(\boldsymbol{a},\varepsilon)),\boldsymbol{y}(\boldsymbol{x}))=0$. Being $\imt(\boldsymbol{y},B(\boldsymbol{a},\varepsilon))$ open, we  deduce $\boldsymbol{y}(\boldsymbol{x})\notin \imt(\boldsymbol{y},B(\boldsymbol{a},\varepsilon))$. 
\end{proof}

The next result provides the analogue of \cite[Lemma 3.5]{mueller.spector}.

\begin{lemma}[Degree and perimeter of topological image]
	\label{lem:degree-perimeter-imt}
Let $\boldsymbol{y}\in W^{1,p}(\Omega^A_\varepsilon;\RN)$ satisfy properties {\rm (i)} and {\rm (iv)} in Definition \ref{def:deformation-perforated}. Then, the following holds: 
\begin{enumerate}[\rm (i)]
	\item For all $U\in\mathcal{U}_{\boldsymbol{y}}^\prime$, denoting by $\overline{\boldsymbol{y}}$  the continuous representative of the trace of $\boldsymbol{y}$ on $\partial U$, we have
	\begin{equation*}
		\text{$\deg(\boldsymbol{y},U,\boldsymbol{\xi})\in \{0,1\}$ for all $\boldsymbol{\xi}\in \RN\setminus \overline{\boldsymbol{y}}(\partial U)$.}
	\end{equation*}
		Moreover,
	\begin{equation*}
		\partial^* \imt(\boldsymbol{y},U)\simeq \img(\boldsymbol{y},\partial U)\simeq \overline{\boldsymbol{y}}(\partial U).
	\end{equation*}
	\item For all $\boldsymbol{a}\in A$, denoting by $\overline{\boldsymbol{y}}$  the continuous representative of the trace   $\boldsymbol{y}^*$   of $\boldsymbol{y}$ on $S(\boldsymbol{a},\varepsilon)$, we have 
	\begin{equation}
		\label{eq:d}
		\text{$\deg(\boldsymbol{y},B(\boldsymbol{a},\varepsilon),\boldsymbol{\xi})\in \{0,1\}$ for all $\boldsymbol{\xi}\in \RN\setminus \overline{\boldsymbol{y}}(S(\boldsymbol{a},\varepsilon))$.}
	\end{equation}
	Moreover, 
	\begin{equation}
		\label{eq:ded}
		\partial^* \imt(\boldsymbol{y},B(\boldsymbol{a},\varepsilon))\simeq 
		 \img(\boldsymbol{y}^*,S(\boldsymbol{a},\varepsilon)) \simeq \overline{\boldsymbol{y}}(S(\boldsymbol{a},\varepsilon)).
	\end{equation}
	Eventually,  it holds that  
	\begin{equation}\label{eq:volume-imt-B}
		\leb(\imt(\boldsymbol{y},B(\boldsymbol{a},\varepsilon)))=\frac{1}{N}\int_{S(\boldsymbol{a},\varepsilon)}  \boldsymbol{y}^*  \cdot (\cof D^{\rm t} \boldsymbol{y}^*\EEE )\boldsymbol{\nu}_{B(\boldsymbol{a},\varepsilon)}\,\d\haus
	\end{equation}
	and
	\begin{equation}
		\label{eq:perimeter-imt-B}
		\per \left( \imt(\boldsymbol{y},B(\boldsymbol{a},\varepsilon))  \right)=\int_{S(\boldsymbol{a},\varepsilon)} |(\cof D^{\rm t} \boldsymbol{y}^* )\boldsymbol{\nu}_{B(\boldsymbol{a},\varepsilon)}|\,\d\haus. 
	\end{equation}
\end{enumerate}	
\end{lemma}
\begin{proof}
The proof of (i) works as in \cite[Lemma 3.5, Steps 1--5]{mueller.spector}. For (ii),  following the same arguments in  \cite[Lemma 3.5, Steps 1--4]{mueller.spector} with the aid of property (iv.2) in Definition \ref{def:deformation-perforated},  we show that $\deg(\boldsymbol{y},B(\boldsymbol{a},\varepsilon),\boldsymbol{\xi})\in \{0,1,-1\}$ for all $\boldsymbol{\xi}\in \RN\setminus \overline{\boldsymbol{y}}(S(\boldsymbol{a},\varepsilon))$. In view of property   (iv.3)  in Definition~\ref{def:deformation-perforated},   claim  \eqref{eq:d} follows. Thus, $\deg(\boldsymbol{y},B(\boldsymbol{a},\varepsilon),\cdot)=\chi_{\imt(\boldsymbol{y},B(\boldsymbol{a},\varepsilon))}$. By Proposition~\ref{prop:degree-divergence}, we find
\begin{equation}
	\label{eq:a}
	\int_{\imt(\boldsymbol{y},B(\boldsymbol{a},\varepsilon))} \div \boldsymbol{\psi}\,\d\boldsymbol{\xi}=\int_{S(\boldsymbol{a},\varepsilon)} \boldsymbol{\psi}\circ \boldsymbol{y}\cdot (\cof D^{\rm t} \boldsymbol{y}^*)\boldsymbol{\nu}_{B(\boldsymbol{a},\varepsilon)}\,\d \haus
\end{equation} 
for all $\boldsymbol{\psi}\in C^1(\RN;\RN)$. Setting 
\begin{equation*}
	\widetilde{\boldsymbol{\nu}}( \boldsymbol{y}^* (\boldsymbol{x}))\coloneqq \frac{(\cof D^{\rm t} \boldsymbol{y}^*(\boldsymbol{x}))\boldsymbol{\nu}_{B(\boldsymbol{a},\varepsilon)}(\boldsymbol{x})  }{|(\cof D^{\rm t} \boldsymbol{y}^*(\boldsymbol{x}))\boldsymbol{\nu}_{B(\boldsymbol{a},\varepsilon)}(\boldsymbol{x})|} \quad \text{for  all $\boldsymbol{x}\in  \domg(\boldsymbol{y}^*,S(\boldsymbol{a},\varepsilon))$,}
\end{equation*}
 with  the aid of  Proposition \ref{prop:cov-surface}  and \eqref{eq:LusinN},  the right--hand side of \eqref{eq:a} can be rewritten  as
\begin{equation*}
	\int_{ \img(\boldsymbol{y}^*,\, S(\boldsymbol{a},\varepsilon)) } \boldsymbol{\psi}\cdot \widetilde{\boldsymbol{\nu}}\,\d\haus = \int_{\overline{\boldsymbol{y}}(S(\boldsymbol{a},\varepsilon ))}\boldsymbol{\psi}\cdot \widetilde{\boldsymbol{\nu}}\,\d\haus.
\end{equation*} 
This gives \eqref{eq:ded}. Eventually, \eqref{eq:volume-imt-B} follows by applying Proposition \ref{prop:degree-divergence} with $\boldsymbol{\psi}(\boldsymbol{\xi})\coloneqq \frac{1}{N} \boldsymbol{\xi}$, while \eqref{eq:perimeter-imt-B} is a consequence of Proposition \ref{prop:cov-surface} and \eqref{eq:ded}. 
\end{proof}

The next lemma is the  analogue  of \cite[Proposition 2.29]{bresciani.friedrich.moracorral} for perforated domains. 

\begin{lemma}[Properties of topological image]
	\label{lem:imt}
Let $\boldsymbol{y}\in W^{1,p}(\Omega^A_\varepsilon;\RN)$ satisfy properties {\rm (i)} and {\rm (iv)} in Definition \ref{def:deformation-perforated}. Then, the following holds: 
\begin{enumerate}[\rm (i)]
	\item For all $U\in\mathcal{U}_{\boldsymbol{y}}^\prime$, we have 
	\begin{equation*}
		\bigcup_{\boldsymbol{a}\in A \cap U} \imt(\boldsymbol{y},B(\boldsymbol{a},\varepsilon))\subset \imt(\boldsymbol{y},U).
	\end{equation*}
	\item For all $\boldsymbol{a},\boldsymbol{b}\in A$ with $\boldsymbol{a}\neq \boldsymbol{b}$, we have
	\begin{equation*}
		\imt(\boldsymbol{y},B(\boldsymbol{a},\varepsilon)) \cap \imt(\boldsymbol{y},B(\boldsymbol{b},\varepsilon))=\emptyset.
	\end{equation*}
\end{enumerate}
\end{lemma}
\begin{proof}
We prove (i).  
 Let $U\in\mathcal{U}_{\boldsymbol{y}}^\prime$ and $\boldsymbol{a}\in A\cap U$. Let $\overline{\boldsymbol{y}}\in C^0(\closure{U};\RN)$ be such that $\overline{\boldsymbol{y}}\simeq  \boldsymbol{y}^* \EEE$ on  $\partial U \cup S(\boldsymbol{a},\varepsilon)$, where we denote both traces of $\boldsymbol{y}$ on $\partial U$ and $S(\boldsymbol{a},\varepsilon)$ by $\boldsymbol{y}^*$. In view   of  property (iv.5) in Definition \ref{def:deformation-perforated} and Lemma \ref{lem:degree-perimeter-imt}(ii), we have 
 \begin{equation*}
 	\haus(\img(\boldsymbol{y},\Omega^A_\varepsilon) \cap \img(\boldsymbol{y}^*,S(\boldsymbol{a},\varepsilon)))=0, \quad \img(\boldsymbol{y}^*,S(\boldsymbol{a},\varepsilon))\simeq \overline{\boldsymbol{y}}(S(\boldsymbol{a},\varepsilon)),
 \end{equation*}
 so that  
\begin{equation*}
	\haus(\overline{\boldsymbol{y}}(S(\boldsymbol{a},\varepsilon))\cap \img(\boldsymbol{y},\Omega^A_\varepsilon))=0.
\end{equation*} 
In particular, $\haus(\overline{\boldsymbol{y}}(S(\boldsymbol{a},\varepsilon))\cap \img(\boldsymbol{y},\partial U))=0$. Since $\img(\boldsymbol{y},\partial U)\simeq \overline{\boldsymbol{y}}(\partial U)$ by Lemma \ref{lem:degree-perimeter-imt}(i),  we get
\begin{equation}
	\label{eq:inc1}
	\haus(\overline{\boldsymbol{y}}(S(\boldsymbol{a},\varepsilon))\cap \overline{\boldsymbol{y}}(\partial U))=0.
\end{equation}
By Lemma \ref{lem:img-vs-imt}(ii), we have
\begin{equation*}
	 \img(\boldsymbol{y},\partial U)\subset \img(\boldsymbol{y},\Omega^A_\varepsilon)  \subset \RN \setminus \imt(\boldsymbol{y},B(\boldsymbol{a},\varepsilon)).	
\end{equation*}
From this, using  Lemma \ref{lem:img-vs-imt}(i), the  continuity of $\overline{\boldsymbol{y}}$, and the  density of $\img(\boldsymbol{y},\partial U)$ in $\overline{\boldsymbol{y}}(\partial U)$ given by Lemma \ref{lem:degree-perimeter-imt}(i), we find
\begin{equation*} 
\partial\, \imt(\boldsymbol{y},U) = 	\overline{\boldsymbol{y}}(\partial U)\subset \RN \setminus \imt(\boldsymbol{y},B(\boldsymbol{a},\varepsilon)).  
\end{equation*}
Let $T$ be the set of $\boldsymbol{x}\in S(\boldsymbol{a},\varepsilon)$ such that $\overline{\boldsymbol{y}}(\boldsymbol{x})=\boldsymbol{y}^*(\boldsymbol{x})$.  Given that $\boldsymbol{y}\in L^\infty(A(\boldsymbol{a},\varepsilon,\bar{\varepsilon});\RN)$  because of condition (INV${}^\prime$), in view of \cite[Proposition 3.65]{ambrosio.fusco.pallara} and \cite[Remark 4.4.5]{ziemer.wdf},  we can assume
\begin{equation*}
		\overline{\boldsymbol{y}}(\boldsymbol{x})=\boldsymbol{y}^*(\boldsymbol{x})=  \lim_{\substack{\boldsymbol{z}\to \boldsymbol{x}\\ 
		\boldsymbol{z}\in E}} \boldsymbol{y}(\boldsymbol{x}),
\end{equation*}
 for some set $E\subset U \cap  \Omega^A_\varepsilon$ with $\Theta^N(E,\boldsymbol{x})=1/2$. We have  
 $T\simeq S(\boldsymbol{a},\varepsilon)$. Fix $\boldsymbol{x}\in T$. As $\domg(\boldsymbol{y},U \cap \Omega^A_\varepsilon)\cong U\cap \Omega^A_\varepsilon$, there exists a sequence $(\boldsymbol{x}_n)_n$ in $\domg(\boldsymbol{y},U\cap \Omega^A_\varepsilon)$ such that $\boldsymbol{y}(\boldsymbol{x}_n)\to \overline{\boldsymbol{y}}(\boldsymbol{x})$, as $n\to \infty$.  Moreover,   $\boldsymbol{y}(\boldsymbol{x}_n)\in \img(\boldsymbol{y},U \cap \Omega^A_\varepsilon)\subset  \closure{\imt(\boldsymbol{y},U)}$ for all $n\in\N$ by Lemma \ref{lem:img-vs-imt} (i),  so that  we deduce  $\overline{\boldsymbol{y}}(\boldsymbol{x})\in \closure{\imt(\boldsymbol{y},U)}$. This shows that $\overline{\boldsymbol{y}}(T)\subset \closure{\imt(\boldsymbol{y},U)}$. Using the density of $T$  in $S(\boldsymbol{a},\varepsilon)$ and the continuity of $\overline{\boldsymbol{y}}$, we obtain
\begin{equation}
	\label{eq:inc3}
	\partial\,\imt(\boldsymbol{y},B(\boldsymbol{a},\varepsilon)) =   \overline{\boldsymbol{y}}(S(\boldsymbol{a},\varepsilon))\subset \closure{ \imt(\boldsymbol{y},U)},
\end{equation}
where we resorted to the last inclusion in Lemma \ref{lem:img-vs-imt}(i).  Now,  
owing to \eqref{eq:inc1}--\eqref{eq:inc3}, by applying \cite[Lemma A.1]{mueller.spector} with $A=\imt(\boldsymbol{y},B(\boldsymbol{a},\varepsilon))$ and $D=\RN \setminus \closure{\imt(\boldsymbol{y},U)}$, we find $\imt(\boldsymbol{y},B(\boldsymbol{a},\varepsilon))\subset \closure{\imt(\boldsymbol{y},U)}$ and,  by Lemma \ref{prop:topim}(i),    $\imt(\boldsymbol{y},B(\boldsymbol{a},\varepsilon))\subset \imt(\boldsymbol{y},U)$.  Therefore, claim (i) follows.  

For proving (ii), let $U_{\boldsymbol{a}},U_{\boldsymbol{b}}\in\mathcal{U}_{\boldsymbol{y}}^\prime $ with $\closure{U}_{\boldsymbol{a}} \cap \closure{U}_{\boldsymbol{b}}=\emptyset$ be such that $B(\boldsymbol{a},\varepsilon)\subset U_{\boldsymbol{a}}$ and $B(\boldsymbol{b},\varepsilon)\subset U_{\boldsymbol{b}}$. Then,   by combining claim (i) and Proposition \ref{prop:topim}(ii),  we find
\begin{equation*}
	\imt(\boldsymbol{y},B(\boldsymbol{a},\varepsilon)) \cap \imt(\boldsymbol{y},B(\boldsymbol{b},\varepsilon))\subset \imt(\boldsymbol{y},U_{\boldsymbol{a}}) \cap \imt(\boldsymbol{y},U_{\boldsymbol{b}})=\emptyset.
\end{equation*}
\end{proof}

The following result registers a simple consequence of the properties of excision   and additivity of the topological degree.

\begin{lemma}
	\label{lem:deg-splitting}
Let $\boldsymbol{y}\in W^{1,p}({\Omega}^A_\varepsilon;\RN)$ satisfy properties  {\rm (i)} and {\rm (iv)} in Definition \ref{def:deformation-perforated}. Also, let $U\in\mathcal{U}_{\boldsymbol{y}}^\prime $ and let $\{U_i:\:i\in I\}\subset \mathcal{U}_{\boldsymbol{y}}^\prime$ be a finite family with $U_i\subset \subset U$ for all $i\in I$ satisfying $\closure{U}_i\cap \closure{U}_j=\emptyset$ for all $i\neq j$. 
Then, there holds
\begin{equation*}
	\imt(\boldsymbol{y},U)\cong \imt\left (\boldsymbol{y},U\setminus \bigcup_{i\in I} \closure{U}_i\right ) \cup \bigcup_{i\in I} \imt(\boldsymbol{y},U_i).
\end{equation*}
\end{lemma}
\begin{proof}
Let $\overline{\boldsymbol{y}}\in C^0(\closure{U};\RN)$ be such that $\overline{\boldsymbol{y}}\simeq  \boldsymbol{y}^* $ on $ \partial U\cup (\bigcup_{i\in I} \partial U_i)\EEE $, where we denote by  the same symbol  $\boldsymbol{y}^*$  the traces of  $\boldsymbol{y}$ on  both  $\partial U$ and $\partial U_i$, respectively. First, let $\boldsymbol{\xi}\in \imt(\boldsymbol{y},U)$ and suppose that $\boldsymbol{\xi}\notin \bigcup_{ i\in I } \overline{\boldsymbol{y}}(\partial U_i)$. By the excision property \cite[Theorem 2.7(2)]{fonseca.gangbo} and the additivity of the degree \cite[Theorem 2.7(1)]{fonseca.gangbo}, we have  	
\begin{equation}
	\label{eq:dd}
		\deg(\boldsymbol{y},U,\boldsymbol{\xi})=\deg\left (\boldsymbol{y},U\setminus \left (\textstyle \bigcup_{i\in I}\overline{\boldsymbol{y}}( \partial U_i)\right),\boldsymbol{\xi}\right )=\deg(\boldsymbol{y},U\setminus (\textstyle \bigcup_{i\in I} \closure{U}_i),\boldsymbol{\xi})+\displaystyle \sum_{i\in I} \deg(\boldsymbol{y},U_i,\boldsymbol{\xi}).
\end{equation}
As the left-hand side equals one  and all the summands on the right-hand side are either one or zero in view of Lemma \ref{lem:degree-perimeter-imt}(i), exactly one of the summands on the right-hand side has to be one and all the others zero. This shows that
\begin{equation}
	\label{eq:i1}
	\imt(\boldsymbol{y},U)\subset  \imt\left (\boldsymbol{y},U\setminus \bigcup_{ i \in I} \closure{U}_i\right ) \cup \bigcup_{i\in I} \imt(\boldsymbol{y},U_i) \cup \bigcup_{i\in I} \overline{\boldsymbol{y}}(\partial U_i).
\end{equation} 
Conversely, let $\boldsymbol{\xi}\in \imt\left (\boldsymbol{y},U\setminus \textstyle \bigcup_{ i \in I} \closure{U_i}\right ) \cup \bigcup_{i\in I} \imt(\boldsymbol{y},U_i) $ and suppose that $\boldsymbol{\xi}\notin \overline{\boldsymbol{y}}(\partial U)$. By excision and additivity, equation \eqref{eq:dd} still holds true. As the left-hand side is either zero or one while one of  the  summands on the right-hand side are positive, again by Lemma \ref{lem:degree-perimeter-imt}(i), we have that $\deg(\boldsymbol{y},U,\boldsymbol{\xi})=1$, i.e., $\boldsymbol{\xi}\in\imt(\boldsymbol{y},U)$. This shows that
\begin{equation}
	\label{eq:i2}
	\imt\left (\boldsymbol{y},U\setminus \textstyle \bigcup_{ i\in I} \closure{U_i}\right ) \cup \bigcup_{i\in I} \imt(\boldsymbol{y},U_i)  \subset \imt(\boldsymbol{y},U)\cup \overline{\boldsymbol{y}}(\partial U).
\end{equation} 
 The two inclusions in     \eqref{eq:i1}--\eqref{eq:i2} yield the claim  given that  $\leb(\overline{\boldsymbol{y}}(\partial U))=\leb(\overline{\boldsymbol{y}}(\partial U_i))=0$ for all $i\in I$.
\end{proof}

\subsection{Extended distributional determinant}
We introduce a generalization of the distributional determinant for deformations defined on perforated domains.  Let $\varepsilon\in (0,\bar{\varepsilon})$ and $A\in\mathcal{A}_\varepsilon$. Recalling that $\widetilde{\Omega}^A_\varepsilon\coloneqq \Omega\setminus \bigcup_{\boldsymbol{a}\in A} B(\boldsymbol{a},\varepsilon)$,  we define    
\begin{equation*}
	\mathcal{F}(\widetilde{\Omega}^A_\varepsilon)\coloneqq \left\{ \varphi\restr{\widetilde{\Omega}^A_\varepsilon}:\:\varphi \in \mathcal{D}(\Omega)   \right\}.
\end{equation*} 
This is a vector space which is naturally endowed with a locally convex topology.

\begin{definition}[Extended distributional determinant]
	\label{def:extended-distributional-determinant}
 Let $\boldsymbol{y}\in W^{1,p}(\Omega^A_\varepsilon;\RN)$ satisfy properties {\rm (i)}, {\rm (ii)}, and {\rm (iv)} in Definition \ref{def:deformation-perforated}. The extended distributional determinant of $\boldsymbol{y}$ is the functional $\Det^A_\varepsilon D \boldsymbol{y}\colon \mathcal{F}(\widetilde{\Omega}^A_\varepsilon)\to \R$ defined by setting
\begin{equation*}
	\begin{split}
		(\Det^A_\varepsilon D \boldsymbol{y})(\varphi)&\coloneqq -\frac{1}{N} \int_{\Omega^A_\varepsilon} ((\adj D \boldsymbol{y})\boldsymbol{y})\cdot D \varphi\,\d\boldsymbol{x}\\
			& \ \ \ -\frac{1}{N} \sum_{\boldsymbol{a}\in A} \int_{S(\boldsymbol{a},\varepsilon)} \boldsymbol{y}^* \cdot (\cof  D^{\rm t}\boldsymbol{y}^*)\boldsymbol{\nu}_{B(\boldsymbol{a},\varepsilon)}\,\varphi \,\d\haus
	\end{split}
\end{equation*}
for all $\varphi\in\mathcal{F}(\widetilde{\Omega}^A_\varepsilon)$. 
\end{definition}

 Concering the previous definition, note that  $\boldsymbol{y}\in L^\infty( K\cap \Omega^A_\varepsilon ;\RN)$ for every compact set $K\subset \Omega$ as a consequence of condition (INV${}^\prime$). Thus,  $(\adj D \boldsymbol{y})\boldsymbol{y}\in L^1(K\cap  \Omega^A_\varepsilon;  \RN)$ for every compact set $K\subset \Omega$ and, in turn, the first integral on the right-hand side is well defined. The same holds for the second integral as $\boldsymbol{y}^*\in W^{1,p}(S(\boldsymbol{a},\varepsilon))$ entails $\boldsymbol{y}^*\in L^\infty(S(\boldsymbol{a},\varepsilon);\RN)$ and $\cof D^{\rm t}\boldsymbol{y}^*\in L^1(S(\boldsymbol{a},\varepsilon);\RNN)$.

Formally, the definition of $\Det^A_\varepsilon D \boldsymbol{y}$ is deduced from the identity $\det D \boldsymbol{y}=\frac{1}{N}\div \left( (\adj D \boldsymbol{y})\boldsymbol{y} \right)$ after integration by parts over $\Omega^A_\varepsilon$ against test functions $\varphi \in\mathcal{F}(\widetilde{\Omega}^A_\varepsilon)$. In particular, for $\boldsymbol{y}\in C^2(\Omega;\RN)$, we find   $\Det^A_\varepsilon D\boldsymbol{y}=(\det D \boldsymbol{y}) \leb$ in $\mathcal{M}_{\rm b}(\widetilde{\Omega}^A_\varepsilon)$.  
We mention that $\Det^A_\varepsilon D \boldsymbol{y}$ as given in Definition \ref{def:extended-distributional-determinant} belongs to the dual space of   $\mathcal{F}(\widetilde{\Omega}^A_\varepsilon)$ although this observation will not be relevant for our discussion.

 Note that $(\Det^A_\varepsilon D\boldsymbol{y})(\varphi)=(\Det D \boldsymbol{y})(\varphi)$ for all $\varphi \in \mathcal{D}(\Omega^A_\varepsilon)$, where $\varphi$ is trivially identified with its extension to $\widetilde{\Omega}^A_\varepsilon$ by zero. Thus,  if we know that  $\Det^A_\varepsilon D \boldsymbol{y}\in\mathcal{M}_{\rm b}(\widetilde{\Omega}^A_\varepsilon)$, then we easily see $\Det D\boldsymbol{y}\in\mathcal{M}_{\rm b}(\Omega^A_\varepsilon)$ and $(\Det^A_\varepsilon D\boldsymbol{y})(U)=(\Det D\boldsymbol{y})(U)$ for all open sets $U\subset \subset \Omega^A_\varepsilon$.

The next result is the analogue of \cite[Corollary 9.2]{mueller.spector} for the extended distributional determinant.

\begin{proposition}[Extended distributional determinant]
	\label{prop:Det-ext}
Let $\boldsymbol{y}\in W^{1,p}(\Omega^A_\varepsilon;\RN)$ satisfy properties {\rm (i)}, {\rm(ii)}, and {\rm (iv)} in Definition \ref{def:deformation-perforated}. Then, $\Det^A_\varepsilon D\boldsymbol{y}\in \mathcal{M}_{\rm b}(\widetilde{\Omega}^A_\varepsilon)$ and we have 	\begin{equation}
		\label{eq:Det-ext}
		(\Det^A_\varepsilon D\boldsymbol{y})(U  \cap \widetilde{\Omega}^A_\varepsilon)  =\leb(\imt(\boldsymbol{y},U \cap \Omega^A_\varepsilon)) \quad \text{for all $U\in\mathcal{U}_{\boldsymbol{y}}^\prime$.}
	\end{equation} 
In particular,  we have
	\begin{equation*}
		\imt(\boldsymbol{y},U\cap \Omega^A_\varepsilon)\cong \img(\boldsymbol{y},U\cap \Omega^A_\varepsilon) \quad \text{for all $U\in\mathcal{U}_{\boldsymbol{y}}^\prime$},
	\end{equation*}
whenever 	$\Det^A_\varepsilon D\boldsymbol{y}=(\det D \boldsymbol{y})\leb$ in $\mathcal{M}_{\rm b}(\widetilde{\Omega}^A_\varepsilon)$.
\end{proposition}

\begin{proof}	
We divide the proof into five steps.	

\emph{Step 1.} Consider $\mu\in \mathcal{D}'(\Omega\setminus A)$ defined as 
\begin{equation*}
	\mu(\varphi)\coloneqq -\frac{1}{N} \int_{\Omega^A_\varepsilon} ((\adj D \boldsymbol{y})\boldsymbol{y})\cdot D \varphi\,\d\boldsymbol{x} \quad \text{for all $\varphi \in \mathcal{D}(\Omega \setminus A)$.}
\end{equation*}	
We prove that $\mu\in \mathcal{M}_{\rm b}(\Omega \setminus A)$ by arguing as in \cite[Lemma 8.1]{mueller.spector}. Let $(\rho_\delta)_{\delta>0}$ be a sequence of standard mollifiers, so that $\rho_\delta (\boldsymbol{x})=h_\delta(|\boldsymbol{x}|)$ for all $\boldsymbol{x}\in\RN$, where $h\in\mathcal{D}(\R)$ and $h_\delta(t)\coloneqq \delta^{-N}h(t/\delta)$. Then, $(\mu\ast \rho_\delta)(\boldsymbol{x})=\mu(\rho_\delta^{\boldsymbol{x}})$, where $\rho_\delta^{\boldsymbol{x}}(\boldsymbol{z})\coloneqq \rho_\delta(\boldsymbol{z}-\boldsymbol{x})$ for all $\boldsymbol{x}\in  \R^N  $. If $\boldsymbol{x}\in B(\boldsymbol{a},\varepsilon)\setminus \{\boldsymbol{a}\}$ for some $\boldsymbol{a}\in A$, then $\closure{B}(\boldsymbol{x},\delta)\subset B(\boldsymbol{a},\varepsilon)$ for $\delta \ll 1$ so that $(\mu\ast \rho_\delta)(\boldsymbol{x})=0$. Instead, if $\boldsymbol{x}\in \Omega^A_\varepsilon$, then $\closure{B}(\boldsymbol{x},\delta)\subset \Omega^A_\varepsilon$ for $\delta\ll 1$ and the same computation  as  in \cite[Lemma 8.1]{mueller.spector} yields $(\mu\ast \rho_\delta)(\boldsymbol{x})\geq 0$. Therefore,  the claim follows  by \cite[Theorem 3.9.4]{ziemer.wdf}  as,  for $\varphi \ge 0$,
\begin{equation*}
	\mu(\varphi)=\lim_{\delta\to 0^+} \int_{\Omega \setminus A} (\mu\ast \rho_\delta) \varphi\,\d\boldsymbol{x}\geq \int_{\Omega \setminus A} \left ( \liminf_{\delta \to 0^+} \mu\ast \rho_\delta \right ) \varphi\,\d\boldsymbol{x}\geq 0,
\end{equation*}
 thanks to   Fatou's lemma. 

\emph{Step 2.} Next, for   $\boldsymbol{a}\in A$ fixed, we show that 
\begin{align}
	\label{eq:T-ball}
	\mu({B}(\boldsymbol{a},r)\setminus \{ \boldsymbol{a} \})&=\leb \big(\imt(\boldsymbol{y},B(\boldsymbol{a},r))\big)
\end{align}
for all $r\in (\varepsilon,2\varepsilon)$ such that $B(\boldsymbol{a},r)\in\mathcal{U}_{\boldsymbol{y}}^\prime$. 
Note that, for any such $r$, we have $B(\boldsymbol{a},r)\cap B(\boldsymbol{b},\varepsilon)=\emptyset$ for all $\boldsymbol{b}\in A \setminus \{\boldsymbol{a} \}$  as $A\in\mathcal{A}_\varepsilon$.   
Also here,  we adapt the arguments in \cite[Lemma 8.1]{mueller.spector}. We continuously extend $\mu$  to the space 
 of functions in $ W^{1,\infty}(\Omega\setminus A)$   whose support is contained in $\Omega \setminus A$  without renaming it. 
For $\delta \in (0,(\varepsilon/2)\wedge (2\varepsilon-r))$, let $f_\delta\colon (0,+\infty)\to [0,1]$ be given by 
\begin{equation*}
	f_\delta(t)\coloneqq \begin{cases}
		0 & \text{if $t\in (0,\delta]\cup [r,+\infty)$, }\\
		\frac{1}{\delta}(t-\delta) & \text{if $t\in (\delta,2\delta)$,}\\
		1 & \text{if $t\in [2\delta,r-\delta]$,}\\		
		-\frac{1}{\delta}(t-r+\delta)+1 & \text{if $t\in (r-\delta,r)$.}
	\end{cases}
\end{equation*} 
Define $\varphi_\delta\in W^{1,\infty}(\Omega \setminus A)$ as $\varphi_\delta(\boldsymbol{x})\coloneqq f_\delta(|\boldsymbol{x}-\boldsymbol{a}|)$. Then,  
$\varphi_\delta(\boldsymbol{x}) \to \chi_{{B}(\boldsymbol{a},r)\setminus \{ \boldsymbol{a}  \}}(\boldsymbol{x})$, as $\delta \to 0^+$, for all $\boldsymbol{x}\in \Omega \setminus A$. Thus,   by the dominated convergence theorem,
\begin{equation*}
	\lim_{\delta \to 0^+} \mu(\varphi_\delta)=\lim_{\delta \to 0^+} \int_{\Omega \setminus A} \varphi_\delta\,\d\mu=\int_{\Omega \setminus A} \chi_{{B}(\boldsymbol{a},r)\setminus \{ \boldsymbol{a}  \}}\,\d\mu =\mu({B}(\boldsymbol{a},r)\setminus \{ \boldsymbol{a}  \}).
\end{equation*} 
We compute
\begin{equation*}
	\begin{split}
		\mu(\varphi_\delta)&=\frac{1}{\delta N} \int_{  A(\boldsymbol{a},r-\delta,r)  } (\adj D\boldsymbol{y}(\boldsymbol{x}))\boldsymbol{y}(\boldsymbol{x})\cdot  \frac{\boldsymbol{x}-\boldsymbol{a}}{|\boldsymbol{x}-\boldsymbol{a}|}\,\d \boldsymbol{x}\\
		&=\frac{1}{\delta} \int_{r-\delta}^r \int_{S(\boldsymbol{a},s)}\frac{ \boldsymbol{y}^* (\boldsymbol{x})}{N}\cdot (\cof D^{\rm t} \boldsymbol{y}^* (\boldsymbol{x}))\boldsymbol{\nu}_{B(\boldsymbol{a},s)}(\boldsymbol{x})\,\d \haus(\boldsymbol{x})\,\d s,
	\end{split}
\end{equation*}
where we applied the coarea formula. Passing to the limit, as $\delta\to 0^+$, in the previous equation with the aid of Definition \ref{def:good}(v) by choosing $\boldsymbol{\psi}\in\mathcal{D}(\RN;\RN)$ such that $\boldsymbol{\psi}=\frac{1}{N}\boldsymbol{ id}$ in a neighborhood of $S(\boldsymbol{a},r)$, we find
\begin{equation*}
	\lim_{\delta\to 0^+} \mu(\varphi_\delta)=\int_{S(\boldsymbol{a},r)} \frac{ \boldsymbol{y}^*(\boldsymbol{x})}{N}\cdot \left( \cof D^{\rm t} \boldsymbol{y}^* (\boldsymbol{x}) \right)\boldsymbol{\nu}_{B(\boldsymbol{a},r)}(\boldsymbol{x})\,\d \haus(\boldsymbol{x})=\leb(\imt(\boldsymbol{y},B(\boldsymbol{a},r))),
\end{equation*}
where we  used    {eq:volume-imt-B}. This proves  \eqref{eq:T-ball}. 

\emph{Step 3.} We show that  
\begin{equation}
	\label{eq:TU}
	\mu(U\setminus A)=\leb(\imt(\boldsymbol{y},U)) \quad \text{for all $U\in\mathcal{U}_{\boldsymbol{y}}^\prime$.}
\end{equation}
For each $\boldsymbol{a}\in A\cap U$, we
choose $r_{\boldsymbol{a}}\in (\varepsilon,2\varepsilon)$  such that
\begin{equation*}
	B(\boldsymbol{a},r_{\boldsymbol{a}})\subset U, \quad B(\boldsymbol{a},r_{\boldsymbol{a}})\in \mathcal{U}_{\boldsymbol{y}}^\prime, \quad \mu(S(\boldsymbol{a},r_{\boldsymbol{a}}))=0.
\end{equation*} 
Note that, for each $\boldsymbol{a}\in A\cap U$, the third property is valid for all but a countable number of radii (see, e.g., \cite[Proposition 2.16]{maggi}). Also, $ U \setminus \textstyle  \bigcup_{\boldsymbol{a}\in A} \closure{B}(\boldsymbol{a},r_{\boldsymbol{a}})\in\mathcal{U}_{\boldsymbol{y}}$.    Given the definition of $\mu$ and Lemma \ref{lem:Det}(ii), we have $\mu(V)=(\Det D \boldsymbol{y})(V)=\leb(\imt(\boldsymbol{y},V))$ for all $V\in\mathcal{U}_{\boldsymbol{y}}$. In particular, 
 \begin{equation}
 	\label{eq:ETET}
 	\mu(U \setminus \textstyle  \bigcup_{\boldsymbol{a}\in A} \closure{B}(\boldsymbol{a},r_{\boldsymbol{a}}))=\leb\left (\imt(\boldsymbol{y}, U \setminus \textstyle  \bigcup_{\boldsymbol{a}\in A} \closure{B}(\boldsymbol{a},r_{\boldsymbol{a}}) )\right ).
 \end{equation}
 Altogether, \eqref{eq:T-ball} and \eqref{eq:ETET} give  
\begin{equation*}
	\begin{split}
		\mu(U\setminus A)
		&=\mu(U \setminus \textstyle  \bigcup_{\boldsymbol{a}\in A} \closure{B}(\boldsymbol{a},r_{\boldsymbol{a}}))+\displaystyle \sum_{\boldsymbol{a}\in A\cap U} \mu({B}(\boldsymbol{a},r_{\boldsymbol{a}})\setminus \{ \boldsymbol{a} \})\\
		&=\leb\left (\imt(\boldsymbol{y}, U \setminus \textstyle  \bigcup_{\boldsymbol{a}\in A} \closure{B}(\boldsymbol{a},r_{\boldsymbol{a}}) )\right ) +\displaystyle \sum_{\boldsymbol{a}\in A\cap U} \leb(\imt(\boldsymbol{y},B(\boldsymbol{a},r_{\boldsymbol{a}})))\\
		&=\leb(\imt(\boldsymbol{y},U)),
	\end{split}
\end{equation*}
 where in the last line  we applied  Lemma \ref{lem:deg-splitting}.  

\emph{Step 4.} To conclude the proof of \eqref{eq:Det-ext}, observe that
\begin{equation*}
	\Det^A_\varepsilon D\boldsymbol{y}= \mu-\frac{1}{N} \sum_{\boldsymbol{a}\in A}  \boldsymbol{y}^* \cdot (\cof D^{\rm t} \boldsymbol{y}^*)\boldsymbol{\nu}_{B(\boldsymbol{a},\varepsilon)} \haus\mres{S(\boldsymbol{a},\varepsilon)} \quad \text{in $\mathcal{M}_{\rm b}(\widetilde{\Omega}^A_\varepsilon)$}.
\end{equation*}
Let $U\in\mathcal{U}_{\boldsymbol{y}}^\prime$. By combining \eqref{eq:volume-imt-B} and \eqref{eq:TU}  with the aid of Lemma \ref{lem:deg-splitting}, we obtain
\begin{equation*}
	\begin{split}
		(\Det^A_\varepsilon D\boldsymbol{y})(U\cap \widetilde{\Omega}^A_\varepsilon)&=\leb(\imt(\boldsymbol{y},U))- \sum_{\boldsymbol{a}\in A\cap U} \leb \big(\imt(\boldsymbol{y},B(\boldsymbol{a},\varepsilon))\big)\\
		&=\leb(\imt(\boldsymbol{y}, U\cap \Omega^A_\varepsilon)).
	\end{split}
\end{equation*}

\emph{Step 5.} Suppose that $\Det^A_\varepsilon D \boldsymbol{y}=(\det D \boldsymbol{y})\leb$ in $\mathcal{M}_{\rm b}(\widetilde{\Omega}^A_\varepsilon)$. From Lemma \ref{lem:img-vs-imt} and Lemma \ref{lem:deg-splitting}, we deduce that 
\begin{equation*}
	\img(\boldsymbol{y},U\cap \Omega^A_\varepsilon) \subset \closure{\imt(\boldsymbol{y},U)}\setminus \left( \bigcup_{\boldsymbol{a}\in A \cap U} \imt(\boldsymbol{y},B(\boldsymbol{a},\varepsilon)) \right) \cong \imt(\boldsymbol{y},U\cap \Omega^A_\varepsilon).
\end{equation*}
Then, \eqref{eq:Det-ext} yields
\begin{equation*}
	\leb(\imt(\boldsymbol{y},U\cap \Omega^A_\varepsilon))=(\Det^A_\varepsilon D \boldsymbol{y})(U  \cap \widetilde{\Omega}^A_\varepsilon) =\int_{U\cap \Omega^A_\varepsilon} \det D \boldsymbol{y} \,\d\boldsymbol{x}=\leb(\img(\boldsymbol{y},U\cap \Omega^A_\varepsilon)),
\end{equation*}
where on the right-hand side we applied Proposition \ref{prop:area-formula-domains}.  This concludes the proof.  
\end{proof}

\section{Proof of the main results}
\label{sec:proof}

\subsection{Equi-coercivity}

We begin by proving the equi-coercivity result. The proof follows the lines of  \cite[Theorem 6.3]{henao.thesis}.

\begin{proof}[Proof of Theorem \ref{thm:equi-coercivity}] 
 Up to passing to a not relabeled subsequence, we can assume that $\H^0(A_n) = m$ for all $n \in \N$ for some $m \le M$.  We write
$A_n=\{ \boldsymbol{a}_1^n,\dots,\boldsymbol{a}_m^n\}$. By compactness,  we find $\boldsymbol{a}_1,\dots,\boldsymbol{a}_m\in H$ such that $\boldsymbol{a}^n_i\to \boldsymbol{a}_i$, as $n\to \infty$, for all $i=1,\dots,m$ along a not relabeled subsequence. We define $A\coloneqq \{\boldsymbol{a}_1,\dots,\boldsymbol{a}_m \}$ where some of the elements might be repeated.

Let $\delta\in \big( 0, \bar{\varepsilon}  \big)$ be sufficiently small in order to have $B(\boldsymbol{a},\delta)\cap B(\boldsymbol{b},\delta)=\emptyset$ for all distinct $\boldsymbol{a},\boldsymbol{b}\in A$. By choosing $\bar{n}(\delta)\in \N$ so that $\varepsilon_n\leq \delta/4$ and $|\boldsymbol{a}^n_i-\boldsymbol{a}_i|\leq \delta/2$ for all $i=1,\dots,m$ and $n\geq \bar{n}(\delta)$, we find
\begin{equation*}
	 \closure{B}(\boldsymbol{a}^n_i,\varepsilon_n)\subset B(\boldsymbol{a}_i,\delta)\quad \text{for all $i=1,\dots,m$ and  $n\geq \bar{n}(\delta)$.}
\end{equation*} 
Thus,
\begin{equation*}
	\Omega_\delta^A\subset  \Omega_{\varepsilon_n}^{A_n} \quad \text{for all $n\geq \bar{n}(\delta)$.}
\end{equation*}
From \eqref{eq:equi-boundedness}, using \eqref{eq:W-growth} and the Poincar\'{e} inequality with trace term, taking into account  item (iii) in Definition \ref{def:deformation}, \EEE  we deduce 
\begin{equation*}
	\sup_{n\geq \bar{n}(\delta)} \left\{ \|\boldsymbol{y}_n\|_{W^{1,p}(\Omega_\delta^A;\RN)}+\|g(\det D \boldsymbol{y}_n)\|_{L^1(\Omega_\delta^A)}   \right\}<+\infty.
\end{equation*}
Owing to this bound, we find a subsequence $(\boldsymbol{y}_{n_k})_k$, possibly depending on $\delta$, and two functions $\boldsymbol{y}_\delta\in W^{1,p}(\Omega_\delta^A;\RN)$ and $\vartheta_\delta\in L^1(\Omega_\delta^A)$ such that $\boldsymbol{y}_{n_k}\wk \boldsymbol{y}_\delta$ in $W^{1,p}(\Omega_\delta^A;\RN)$ and $\det D \boldsymbol{y}_{n_k}\wk \vartheta_\delta$ in $L^1(\Omega_\delta^A)$. 
By  Definition \ref{def:deformation-perforated}(i) and  \cite[Lemma 4.2]{henao}, the map  $\boldsymbol{y}_\delta$ satisfies (INV${}^\prime$). A standard contradiction argument which makes use of \eqref{eq:g} shows that $\vartheta_\delta(\boldsymbol{x})>0$ for almost all $\boldsymbol{x}\in \Omega^A_\delta$.  Then, Lemma \ref{lem:positivity} yields $\det D \boldsymbol{y}_\delta(\boldsymbol{x})>0$ for almost all $\boldsymbol{x}\in \Omega^A_\delta$. From Lemma \ref{lem:weakstar},  we find  \EEE 
\begin{equation*}
	\text{$\Det D \boldsymbol{y}_{n_k}\wks \Det D\boldsymbol{y}_\delta$ \quad in $\mathcal{M}_{\rm b}(\Omega_\delta^A).$}
\end{equation*}
As 
\begin{equation*}
	\text{$\Det D \boldsymbol{y}_{n_k}=\Det^{A_{n_k}}_{\varepsilon_{n_k}}D\boldsymbol{y}_{n_k}=(\det D \boldsymbol{y}_{n_k})\leb \quad \text{in $\mathcal{M}_{\rm b}({\Omega}^{A_{n_k}}_{\varepsilon_{n_k}})$}$ \quad for all $k\in\N$,}
\end{equation*}
we deduce that $\Det D\boldsymbol{y}_\delta=\vartheta_\delta \leb $ in $\mathcal{M}_{\rm b}(\Omega^A_\delta)$.  Moreover, by applying \cite[Lemma 8.1]{mueller.spector}, we conclude $\vartheta_\delta \cong \det D \boldsymbol{y}_\delta$ in $\Omega^A_\delta$. \EEE 

Now, let $(\delta_l)_l$ be a decreasing sequence with $\delta_l\to 0^+$.
Combining the reasoning above with a diagonal argument, we select a subsequence $(\boldsymbol{y}_{n_k})_k$ for which
\begin{equation*}
	\text{$\boldsymbol{y}_{n_k}\wk \boldsymbol{y}_{\delta_l}$ in $W^{1,p}(\Omega_{\delta_l}^A;\RN)$, \quad  $\det D \boldsymbol{y}_{n_k}\wk \det D \boldsymbol{y}_{\delta_l}$ in $L^1(\Omega^A_{\delta_l})$, \quad \text{as $k\to \infty$}, \quad for all $\l\in\N$,}
\end{equation*} 
where $\boldsymbol{y}_{\delta_l}\in W^{1,p}(\Omega^A_{\delta_l};\RN)$ for every $l\in \N$. We define $\boldsymbol{y}\colon \Omega \setminus A\to \RN$ by setting $\boldsymbol{y}(\boldsymbol{x})\coloneqq \boldsymbol{y}_{\delta_l}(\boldsymbol{x})$ for $\boldsymbol{x}\in\Omega_{\delta_l}^A$. The definition is well posed,  $\boldsymbol{y}\in W^{1,p}(\Omega\setminus A;\RN)$, and \eqref{eq:convergence-y} clearly holds. 

We are left to prove that $\boldsymbol{y}\in\mathcal{Y}^A$. By the characterization of Sobolev maps in terms of absolutely continuity on almost every line \cite[Theorem 11.45]{leoni}, we immediately realize that $\boldsymbol{y}\in W^{1,p}(\Omega;\RN)$. As $\boldsymbol{y}_{\delta_l}$ satisfies condition (INV${}^\prime$) for every $l\in \N$, we deduce that $\boldsymbol{y}$ satisfies (INV).  Also the boundary conditions are sastisfied  by the weak continuity of the trace operator.  Trivially, $\det D \boldsymbol{y}(\boldsymbol{x})>0$ almost all $\boldsymbol{x}\in \Omega$.  We show that $\Det D \boldsymbol{y}$ has the  structure  required in Definition \ref{def:deformation}(iv).   By \cite[Lemma 8.1]{mueller.spector}, we have $\Det D \boldsymbol{y}=(\det D \boldsymbol{y})\leb+\mu$ for some $\mu\in\mathcal{M}_{\rm b}(\Omega)$. For every $l\in\N$, we have
\begin{equation*}
	\begin{split}
		\Det D\boldsymbol{y}=\Det D \boldsymbol{y}_{\delta_l}=(\det D \boldsymbol{y}_{\delta_l})\leb 
		=(\det D \boldsymbol{y})\leb \quad \text{in $\mathcal{M}_{\rm b}(\Omega^A_{\delta_l})$},
	\end{split}
\end{equation*}
so that $\mu(\Omega_{\delta_l}^A)=0$. Therefore,
\begin{equation*}
	\begin{split}
		\mu(\Omega)=\mu(A)+\mu\left( \bigcup_{l=1}^\infty \Omega_{\delta_l}^A \right)
		\leq \mu(A)+\sum_{l=1}^\infty \mu\left(  \Omega_{\delta_l}^A \right)=\mu(A),
	\end{split}
\end{equation*}
which entails $\mu=\sum_{\boldsymbol{a}\in A} \kappa_{\boldsymbol{a}} \delta_{\boldsymbol{a}}$ with $\kappa_{\boldsymbol{a}}\coloneqq \mu(\{\boldsymbol{a}\}) $ for all $\boldsymbol{a}\in A$.
\end{proof}

\subsection{Lower bound}

Next, we move to the proof of the lower bound.

\begin{proof}[Proof of Theorem \ref{thm:lb}]
%As in the proof of Theorem \ref{thm:equi-coercivity}, for 	$\delta\in \big(0,\dist(K;\partial \Omega)\big)$, we find $\bar{n}(\delta)\in \N$ such that 
%\begin{equation*}
%	\bigcup_{\boldsymbol{b}\in A_n} B(\boldsymbol{\boldsymbol{b}},\varepsilon_n)\subset \subset \bigcup_{\boldsymbol{a}\in A} {B}(\boldsymbol{a},\delta) \quad \text{for all $n\geq \bar{n}(\delta)$.}
%\end{equation*}
%In particular,  \eqref{eq:localization} holds.
%
%By standard polyconvexity arguments, we find
%\begin{equation*}
%	\liminf_{n\to \infty} \int_{\Omega_{\varepsilon_n}^{A_n}} W(D\boldsymbol{y}_n)\,\d\boldsymbol{x}\geq \liminf_{n\to \infty} \int_{\Omega_{\delta}^A} W(D\boldsymbol{y}_n)\,\d\boldsymbol{x}\geq \int_{\Omega_\delta^A} W(D\boldsymbol{y})\,\d\boldsymbol{x}.
%\end{equation*}
%Thus, letting $\delta\to 0^+$ with the aid of the monotone convergence theorem, we find
%\begin{equation*}
%	\liminf_{n\to \infty} \mathcal{W}_{\varepsilon_n}(A_n,\boldsymbol{y}_n)\geq \mathcal{W}(\boldsymbol{y}).
%\end{equation*}

Without loss of generality, we assume that \eqref{eq:equi-boundedness} holds true. In view of \eqref{eq:W-growth}, we have
\begin{equation*}
	\sup_{n\in\N} \int_{\Omega_{\varepsilon_n}^{A_n}} g(\det D \boldsymbol{y}_n)\,\d\boldsymbol{x}\leq \sup_{n\in\N} \mathcal{E}_{\varepsilon_n}(A_n,\boldsymbol{y}_n) <+\infty.
\end{equation*}
Recalling \eqref{eq:g}, this entails the equi-integrability of the sequence formed by
the extension of the functions $\det D \boldsymbol{y}_n$ to $\Omega$ by zero.

Let $\delta\in \left(0,\bar{\varepsilon}\right)$. As in the proof of Theorem \ref{thm:equi-coercivity}, there exists $\bar{n}(\delta)\in\N$ such that $\Omega^A_{\delta}\subset \Omega^{A_n}_{\varepsilon_n}$ for all $n\geq \bar{n}(\delta)$. Given \eqref{eq:convergence-y}, by standard polyconvexity arguments, we find
\begin{equation*}
	\liminf_{n\to \infty} \int_{\Omega^{A_n}_{\varepsilon_n}} W(D\boldsymbol{y}_n)\,\d\boldsymbol{x}\geq \liminf_{n\to \infty} \int_{\Omega^{A}_{\delta}} W(D\boldsymbol{y}_n)\,\d\boldsymbol{x} \geq \int_{\Omega^{A}_{\delta}} W(D\boldsymbol{y})\,\d\boldsymbol{x}.
\end{equation*}
Thus, letting $\delta \to 0^+$ with the aid of the monotone convergence theorem, we obtain
\begin{equation}
	\label{eq:W-lsc}
	\liminf_{n\to \infty} \mathcal{W}_{\varepsilon_n} (A_n,\boldsymbol{y}_n)\geq  \mathcal{W}(\boldsymbol{y}).  
\end{equation}
Now, let $m\coloneqq \mathscr{H}^0(A)$ and write $A=\{ \boldsymbol{a}_1,\dots,\boldsymbol{a}_m \}$. For every $n\in\N$, we select $\boldsymbol{a}^n_1,\dots,\boldsymbol{a}^n_m\in A_n$ such that $\boldsymbol{a}^n_i\to \boldsymbol{a}_i$, as $n\to \infty$, for all $i=1,\dots,m$. Let $(\delta_l)_l$ be a decreasing sequence  with $\delta_l\to 0^+$.   By applying \cite[Lemma 4.1]{henao} together with a diagonal argument,  we select a subsequence $(\boldsymbol{y}_{n_k})_k$ and, for every $l\in \N$,  we fix $r_l\in (\delta_l,2\delta_l)$ satisfying
\begin{equation*}
	B(\boldsymbol{a}_i,r_l) \in \mathcal{U}_{\boldsymbol{y}} \cap \bigcap_{n\geq \bar{n}(\delta_l)} \mathcal{U}_{\boldsymbol{y}_n}^\prime  \quad \text{for all $i=1,\dots,m$.}
\end{equation*}
and
\begin{equation}
	\label{eq:wk-bdry}
	\text{$ \boldsymbol{y}_{n_k}^*\wk \boldsymbol{y}^*$ in $W^{1,p}(S(\boldsymbol{a}_i,r_{l});\RN)$, \:as $k\to \infty$,\: for all $i=1,\dots,m$. }
\end{equation}
 Then, we estimate
\begin{equation}
	\label{eq:triangle}
	\begin{split}
		\|\chi_{\imt(\boldsymbol{y}_{n_k},B(\boldsymbol{a}_i^{n_k},\, \varepsilon_{n_k}))}&-\chi_{\imt(\boldsymbol{y},\boldsymbol{a}_i)}\|_{L^1(\RN)}\\
		&\leq \|\chi_{\imt(\boldsymbol{y}_{n_k},B(\boldsymbol{a}_i^{n_k},\varepsilon_{n_k}))}-\chi_{\imt(\boldsymbol{y}_{n_k},B(\boldsymbol{a}_i,r_l))}\|_{L^1(\RN)}\\
		&+\|\chi_{\imt(\boldsymbol{y}_{n_k},B(\boldsymbol{a}_i,r_l))}-\chi_{\imt(\boldsymbol{y},B(\boldsymbol{a}_i,r_l))}\|_{L^1(\RN)}\\
		&+\|\chi_{\imt(\boldsymbol{y},B(\boldsymbol{a}_i,r_l))}-\chi_{\imt(\boldsymbol{y},\boldsymbol{a}_i)}\|_{L^1(\RN)}.
	\end{split}
\end{equation}
We look at the first summand of the right-hand side of \eqref{eq:triangle}. As $B(\boldsymbol{a}^{n_k}_i,\varepsilon_{n_k})\subset \subset B(\boldsymbol{a}_i,r_l)$  for $k\gg 1$,  it  holds that $\imt(\boldsymbol{y}_{n_k},B(\boldsymbol{a}_i^{n_k},\varepsilon_{n_k}))\subset \imt(\boldsymbol{y}_{n_k},B(\boldsymbol{a}_i,r_l))$ by  Lemma \ref{lem:imt}(i).   Thus,  by applying Lemma \ref{lem:deg-splitting},  Proposition \ref{prop:Det-ext} with Definition \ref{def:deformation-perforated}(v), and Proposition \ref{prop:area-formula-domains}, we estimate   
\begin{equation*}
	\begin{split}
		\| \chi_{\imt(\boldsymbol{y}_{n_k},B(\boldsymbol{a}_i^{n_k},\varepsilon_{n_k}))}-\chi_{\imt(\boldsymbol{y}_{n_k},B(\boldsymbol{a}_i,r_l))} \|_{L^1(\RN)}&=\leb \big(\imt(\boldsymbol{y}_{n_k},B(\boldsymbol{a}_i,r_l)\setminus \closure{B}(\boldsymbol{a}_i^{n_k},\varepsilon_{n_k}))\big)\\
		&=\leb \big(\img(\boldsymbol{y}_{n_k},B(\boldsymbol{a}_i,r_l)\setminus \closure{B}(\boldsymbol{a}_i^{n_k},\varepsilon_{n_k}))\big)\\
		&=\int_{B(\boldsymbol{a}_i,r_l)\setminus \closure{B}(\boldsymbol{a}_i^{n_k},\varepsilon_{n_k})} \det D \boldsymbol{y}_{n_k}\,\d\boldsymbol{x} \\
		&\leq \int_{B(\boldsymbol{a}_i,r_l)} \det D \boldsymbol{y}_{n_k}\,\d\boldsymbol{x}.
	\end{split}
\end{equation*}
Hence,  given the equi-integrability of the extensions of $\det D \boldsymbol{y}_{n_k}$ to $\Omega$ by zero,  the first summand on the right-hand side of \eqref{eq:triangle} goes to zero, as $l\to \infty$. The same holds for 
the second summand  because of \eqref{eq:wk-bdry}  and \cite[Lemma 2.19]{bresciani.friedrich.moracorral} as well as for the third one in view of \eqref{eq:imtop-point-int2}.   Therefore, we proved that
\begin{equation}
	\label{eq:main-convergence}
	\text{$\chi_{\imt(\boldsymbol{y}_{n_k},B(\boldsymbol{a}_i^{n_k},\, \varepsilon_{n_k}))}\to \chi_{\imt(\boldsymbol{y},\boldsymbol{a}_i)}$  in $L^1(\RN)$, \quad  as $k\to \infty$, \quad for all $i=1,\dots,m$.}
\end{equation}
It is not restrictive to assume that the inferior limits of the two  sequences $\left( \mathcal{V}_{\varepsilon_n}(A_n,\boldsymbol{y}_n) \right)_n$ and $\left( \mathcal{P}_{\varepsilon_n}(A_n,\boldsymbol{y}_n) \right)_n$ are both attained as limits along the subsequences indexed by $(n_k)_k$. In this case, 
\eqref{eq:main-convergence} immediately gives 
\begin{equation*}
	\begin{split}
		 \lim_{k\to \infty}\mathcal{V}_{\varepsilon_{n_k}}(A_{n_k},\boldsymbol{y}_{n_k})&=\lim_{k\to \infty} \sum_{\boldsymbol{a}\in A_{n_k}} \leb \big(\imt(\boldsymbol{y}_{n_k},B(\boldsymbol{a},\varepsilon_{n_k}))\big)\\
		 &\geq  \sum_{i=1}^m \lim_{k\to \infty} \leb \big(\imt(\boldsymbol{y}_{n_k},B(\boldsymbol{a}_i^{n_k},\varepsilon_{n_k}))\big) = \sum_{i=1}^m \leb(\imt(\boldsymbol{y},\boldsymbol{a}_i))=\mathcal{V}(A,\boldsymbol{y}).
	\end{split}
\end{equation*}
Similarly, given the lower semicontinuity of the perimeter, we obtain
\begin{equation*}
	\begin{split}
		\lim_{k\to \infty} \mathcal{P}_{\varepsilon_{n_k}}(A_{n_k},\boldsymbol{y}_{n_k})&=\lim_{k\to \infty} \sum_{\boldsymbol{a}\in A_{n_k}} \per \left ( \imt(\boldsymbol{y}_{n_k},B(\boldsymbol{a},\varepsilon_{n_k})) \right) \\
		&\geq \sum_{i=1}^m \liminf_{k\to \infty} \per \left ( \imt(\boldsymbol{y}_{n_k},B(\boldsymbol{a}_i^{n_k},\varepsilon_{n_k})) \right) \geq \sum_{i=1}^m \per \left( \imt(\boldsymbol{y},\boldsymbol{a}_i) \right)=\mathcal{P}(A,\boldsymbol{y}).
	\end{split}
\end{equation*}
The combination of the last two inequalities with \eqref{eq:W-lsc} yields \eqref{eq:lb}. 
\end{proof}

\subsection{Optimality of the lower  bound}

\label{subsec:ub}

For convenience, we register a preliminary result before moving to the proof of the optimality of the lower bound.

\begin{lemma}\label{lem:recovery}
Let $(A,\boldsymbol{y})\in\mathcal{Q}$. Then, for every $\varepsilon \in (0,\bar{\varepsilon})$  such that  $B(\boldsymbol{a},\varepsilon)\in\mathcal{U}_{\boldsymbol{y}}$ for all $\boldsymbol{a}\in {A}$, there holds $\boldsymbol{y}\restr{\Omega^A_\varepsilon} \in \mathcal{Y}^{A}_\varepsilon$.
\end{lemma}
\begin{proof}
We check that $\boldsymbol{y}^A_\varepsilon\coloneqq \boldsymbol{y}\restr{\Omega^A_\varepsilon}$  satisfies all properties in Definition \ref{def:deformation-perforated}. Properties (i) and (iii) are trivial, while property (ii) holds as $\boldsymbol{y}$ satisfies condition (INV).    We look at (iv). From  items (i)--(iii) in Definition \ref{def:good}, we see that (iv.1)--(iv.2) are satisfied. Also, we find $D^{\rm t}\boldsymbol{y}^*\simeq D\boldsymbol{y}(\boldsymbol{I}-\boldsymbol{\nu}_{B(\boldsymbol{a},\varepsilon)}\otimes \boldsymbol{\nu}_{B(\boldsymbol{a},\varepsilon)})$ on $S(\boldsymbol{a},\varepsilon)$ which, together with \eqref{eq:cof} and the fact that $\det (\cof D \boldsymbol{y}(\boldsymbol{x}))=(\det D \boldsymbol{y}(\boldsymbol{x}))^{N-1}>0$ for all $\boldsymbol{x}\in \domg(\boldsymbol{y},S(\boldsymbol{a},\varepsilon))$, shows that $\domg(\boldsymbol{y}^*,S(\boldsymbol{a},\varepsilon))\simeq S(\boldsymbol{a},\varepsilon)$ as required in item (iv.5).
 Next, observe that $\boldsymbol{y}^*(\boldsymbol{x})=\boldsymbol{y}(\boldsymbol{x})$ for all  $\boldsymbol{x}\in \domg(\boldsymbol{y},S(\boldsymbol{a},\varepsilon))$ and $\boldsymbol{a}\in A $  in view of \cite[Remark 4.4.5]{ziemer.wdf}. Recalling that $\boldsymbol{y}\restr{\domg(\boldsymbol{y},\Omega)}$ is injective by \cite[Lemma 3]{henao.moracorral.fracture}, we deduce that (iv.4)--(iv.5) hold true. 
By Lemma \ref{lem:degree-perimeter-imt}(i),  $\boldsymbol{y}^A_\varepsilon$ satisfies  (iv.3) in Definition \ref{def:deformation-perforated} and we are left to check (v). From Proposition~\ref{prop:Det-ext}, we know that $\Det^A_\varepsilon D\boldsymbol{y}^A_\varepsilon\in\mathcal{M}_{\rm b}(\widetilde{\Omega}^A_\varepsilon)$ and, for every $U\in\mathcal{U}^\prime_{\boldsymbol{y}^A_\varepsilon}$, we have,  by Lemma \ref{lem:Det}(i), 
\begin{equation*}
	\begin{split}
		(\Det^A_\varepsilon D \boldsymbol{y}^A_\varepsilon)(U\cap \widetilde{\Omega}^A_\varepsilon)&=\leb(\imt(\boldsymbol{y},U\cap \Omega^A_\varepsilon))=(\Det D \boldsymbol{y})(U\cap \Omega^A_\varepsilon)\\
		&=\int_{U \cap \Omega^A_\varepsilon} \det D \boldsymbol{y}\,\d\boldsymbol{x}=\int_{U \cap \Omega^A_\varepsilon} \det D \boldsymbol{y}^A_\varepsilon\,\d\boldsymbol{x}.
	\end{split}
\end{equation*}	
This proves  $\Det^A_\varepsilon D \boldsymbol{y}^A_\varepsilon=(\det D \boldsymbol{y})\leb$ in $\mathcal{M}_{\rm b}(\widetilde{\Omega}^A_\varepsilon)$.
\end{proof}

We are ready to prove the upper bound.

\begin{proof}[Proof of Theorem \ref{thm:ub}]
Set $m\coloneqq \mathscr{H}^0(A)$ and define $A_n\coloneqq A$ for all $n\in\N$. Thus, \eqref{eq:convergence-A} trivially holds and $A_n\in\mathcal{A}_{\varepsilon_n}$ for $n\gg 1$. 
The construction of the deformations works as follows. Fix $n\in \N$ and choose $r_n\in (\varepsilon_n- \delta_n,\varepsilon_n +  \delta_n)$ for some $0< \delta_n\ll  \varepsilon_n $ such that $B(\boldsymbol{a},r_n)\in \mathcal{U}_{\boldsymbol{y}}\EEE$ for all $\boldsymbol{a}\in A$. Let $\phi_n\in C^\infty([0,+\infty))$ be a  strictly increasing function  satisfying 
\begin{equation}
	\label{eq:phi-value}
	\phi_n(0)=0  \qquad \phi_n(\varepsilon_n)=r_n, \qquad \text{$\phi_n(t)=t$ \quad  for all $t\geq 2\varepsilon_n$.}
\end{equation}
For $\delta_n \ll  \varepsilon_n  $, a function $\phi_n$ as above can be constructed to additionally have
\begin{equation}
	\label{eq:phi-uniform}
\left |  \frac{\phi_n(t)}{t}-1 \right | + |\phi_n'(t)-1|   \leq \frac{1}{n} \quad \text{for all $t\geq 0$.}	
\end{equation}
This can be done by mollifying a piecewise affine function which is affine with slope exactly one in  $(\varepsilon_n- \delta_n,\varepsilon_n+\delta_n)$. In this way,  the value of the function at $\varepsilon_n$ remains unchanged after mollification. 

%\begin{figure}
%		\begin{tikzpicture}[scale=.7]
%		\begin{axis}[
%			xmin=0,xmax=5,
%			ymin=0,ymax=5,
%			axis x line=bottom,
%			axis y line=middle,
%			axis line style=->,
%			%xlabel={$x$},
%			%ylabel={$y$},
%			xlabel style={align=center},
%			xticklabels={$0$, $\varepsilon_n-\eta_n$,$\varepsilon_n$,$\varepsilon_n+\eta_n$,$2\varepsilon_n$},
%			xtick={0,1.6,2,2.4,4},
%			ytick={2.1,2.5,2.9,4},
%			yticklabels={$r_n-\eta_n$,$r_n$,$r_n+\eta_n$,$2\varepsilon_n$}
%			]
%			\addplot+ [
%			sharp plot,
%			] coordinates {
%				(0,0) (1.6,2.1)  (2.4,2.9) (4,4)
%			};
%		\end{axis}
%	\end{tikzpicture}
%	\begin{tikzpicture}[scale=.7]
%		\begin{axis}[
%			xmin=0,xmax=5,
%			ymin=0,ymax=5,
%			axis x line=bottom,
%			axis y line=middle,
%			axis line style=->,
%			%xlabel={$x$},
%			%ylabel={$y$},
%			xlabel style={align=center},
%			xticklabels={$0$, $\varepsilon_n-\eta_n$,$\varepsilon_n$,$\varepsilon_n+\eta_n$,$2\varepsilon_n$},
%			xtick={0,1.6,2,2.4,4},
%			ytick={1.9,2.3,2.7,4},
%			yticklabels={$r_n-\eta_n$,$r_n$,$r_n+\eta_n$,$2\varepsilon_n$}
%			]
%%			\addplot+ [
%%			smooth,
%%			] coordinates {
%%				(0,0) (1.4,1.9) 
%%			};
%		\end{axis}
%	
%	\end{tikzpicture}
%\end{figure}

We define a smooth diffeomorphism $\boldsymbol{f}_n\colon \Omega \to  \Omega$ by setting
\begin{equation*}
	\boldsymbol{f}_n(\boldsymbol{x})\coloneqq \begin{cases}
		\boldsymbol{a}+\phi_n(|\boldsymbol{x}-\boldsymbol{a}|)\frac{\boldsymbol{x}-\boldsymbol{a}}{|\boldsymbol{x}-\boldsymbol{a}|} & \text{if $\boldsymbol{x}\in \closure{B}(\boldsymbol{a},2\varepsilon_n)$ for some $\boldsymbol{a}\in A$,}\\
		\boldsymbol{x} & \text{if $\boldsymbol{x}\in \Omega\setminus \textstyle \bigcup_{\boldsymbol{a}\in A} \closure{B}(\boldsymbol{a},2\varepsilon_n)$.}
	\end{cases}
\end{equation*}  
Note that
\begin{equation}
	\label{eq:fn-fix}
	\text{$\boldsymbol{f}_n(\boldsymbol{a})=\boldsymbol{a}$, \qquad $\boldsymbol{f}_n(B(\boldsymbol{a},\varepsilon_n))=B(\boldsymbol{a},r_n)$ \qquad for all $\boldsymbol{a}\in A$, }
\end{equation}
 where the first claim follows  by \eqref{eq:phi-uniform}.  
With the aid of \cite[Proposition A.1 and Corollary A.3]{bresciani.friedrich.moracorral}, we compute
\begin{equation*}
	D\boldsymbol{f}_n(\boldsymbol{x})=\frac{\phi_n(|\boldsymbol{x}-\boldsymbol{a}|)}{|\boldsymbol{x}-\boldsymbol{a}|} \left( \boldsymbol{I}- \frac{\boldsymbol{x}-\boldsymbol{a}}{|\boldsymbol{x}-\boldsymbol{a}|}\otimes \frac{\boldsymbol{x}-\boldsymbol{a}}{|\boldsymbol{x}-\boldsymbol{a}|}  \right) + \phi'_n(|\boldsymbol{x}-\boldsymbol{a}|) \frac{\boldsymbol{x}-\boldsymbol{a}}{|\boldsymbol{x}-\boldsymbol{a}|}\otimes \frac{\boldsymbol{x}-\boldsymbol{a}}{|\boldsymbol{x}-\boldsymbol{a}|}
\end{equation*}
and
\begin{equation*}
	\det D \boldsymbol{f}_n(\boldsymbol{x})= \left( \frac{\phi_n(|\boldsymbol{x}-\boldsymbol{a}|)}{|\boldsymbol{x}-\boldsymbol{a}|}  \right)^{N-1} \phi_n'(|\boldsymbol{x}-\boldsymbol{a}|)
\end{equation*} 
for all $\boldsymbol{a}\in A$ and $\boldsymbol{x}\in B(\boldsymbol{a},2\varepsilon_n)$. Thus, $\boldsymbol{f}_n$ is orientation-preserving, i.e.,  $\det D \boldsymbol{f}_n(\boldsymbol{x})>0$ for all $\boldsymbol{x}\in\Omega$. Moreover, as a consequence of \eqref{eq:phi-uniform}, we obtain
\begin{equation}
	\label{eq:fn-lip}
	\|\boldsymbol{f}_n-\boldsymbol{id}\|_{W^{1,\infty}(\Omega;\RN)}\leq \frac{ 2(\varepsilon_n+1)+\sqrt{N}}{n} \quad \text{for all $n\in \N$.}
\end{equation}
Define the deformation $\widetilde{\boldsymbol{y}}_n\coloneqq \boldsymbol{y}\circ \boldsymbol{f}_n\in W^{1,p}(\Omega;\RN)$. By the chain rule, $D \widetilde{\boldsymbol{y}}_n\cong (D\boldsymbol{y})\circ \boldsymbol{f}_n (D\boldsymbol{f}_n)$ and, hence,
\begin{equation*}
	\det D\widetilde{\boldsymbol{y}}_n(\boldsymbol{x})=\det D\boldsymbol{y}(\boldsymbol{f}_n(\boldsymbol{x}))\det D \boldsymbol{f}_n(\boldsymbol{x})>0 \quad \text{for almost all $\boldsymbol{x}\in \Omega$.}
\end{equation*}
Trivially, $\widetilde{\boldsymbol{y}}_n^*\simeq \boldsymbol{d}^*$  on $\Gamma$. Also,  $\widetilde{\boldsymbol{y}}_n$ satisfies condition (INV) by \cite[Theorem~9.1]{mueller.spector} and,  owing to \cite[Theorem 6.5]{sivaloganathan.spector} and \eqref{eq:fn-fix},  it holds
\begin{equation*}
	\Det D\widetilde{\boldsymbol{y}}_n=(\det D\widetilde{\boldsymbol{y}}_n )\leb +  \sum_{\boldsymbol{a}\in A} \kappa_{\boldsymbol{a}} \delta_{\boldsymbol{a}} \quad \text{in $\mathcal{M}_{\rm b}(\Omega)$.}
\end{equation*}
 
Therefore,  $\widetilde{\boldsymbol{y}}_n\in\mathcal{Y}^{A}$ and $(A,\widetilde{\boldsymbol{y}}_n)\in\mathcal{Q}$. Moreover, as in \cite[Proposition 5.2]{moracorral} and \cite[Theorem 9.1]{mueller.spector},  we find
\begin{equation}
	\label{eq:imt-recovery}
	\imt(\widetilde{\boldsymbol{y}}_n,\boldsymbol{a})=\imt(\boldsymbol{y},\boldsymbol{a}), \qquad \imt(\widetilde{\boldsymbol{y}}_n,B(\boldsymbol{a},\varepsilon_n))=\imt(\boldsymbol{y},B(\boldsymbol{a},r_n)) \quad \text{for all $\boldsymbol{a}\in A$,}
\end{equation}
where we resorted again to \eqref{eq:fn-fix}.

Now, let  $\boldsymbol{a}\in A$ and denote by  $\widetilde{\boldsymbol{y}}_n^*$   the trace of $\widetilde{\boldsymbol{y}}_n$ on $S(\boldsymbol{a},\varepsilon_n)$ and by   $\boldsymbol{y}^*$  the trace of $\boldsymbol{y}$ on $S(\boldsymbol{a},r_n)$, respectively. In view of \eqref{eq:phi-value}, we have $ \widetilde{\boldsymbol{y}}_n^*\simeq  \boldsymbol{y}^* \circ \boldsymbol{f}_n$ on $S(\boldsymbol{a},\varepsilon_n)$. As $B(\boldsymbol{a},r_n)\in\mathcal{U}_{\boldsymbol{y}}$, we have $B(\boldsymbol{a},\varepsilon_n)\in\mathcal{U}_{\widetilde{\boldsymbol{y}}_n}$. Therefore, by Lemma \ref{lem:recovery}, we find  $\boldsymbol{y}_n\coloneqq \widetilde{\boldsymbol{y}}_n\restr{\Omega^{A_n}_{\varepsilon_n}}\in \mathcal{Y}_{\varepsilon_n}^{A_n}$ and, in turn, $(A_n,{\boldsymbol{y}}_n)\in\mathcal{Q}_{\varepsilon_n}$. By observing that  $\boldsymbol{y}_n=\boldsymbol{y}$ in $\Omega^{A}_{2\varepsilon_n}$, we immediately realize that \eqref{eq:convergence-y} holds.

We move to the proof of \eqref{eq:ub}. First, thanks to assumption (v) on $W$ and   \cite[Proposition 1.5]{dalmaso.lazzaroni}, we know that there exist $\gamma\in (0,1)$ such that 
\begin{equation}\label{eq:stress-control}
	W(\boldsymbol{F}\boldsymbol{G})\leq \frac{N}{N-1} \left( W(\boldsymbol{F}) + c_0 \right) \quad \text{for all $\boldsymbol{F},\boldsymbol{G}\in\RNN_+$ with $|\boldsymbol{G}-\boldsymbol{I}|<\gamma$.}
\end{equation} 
Note that \eqref{eq:fn-lip} yields
\begin{equation}
	\label{eq:fn}
	\|D\boldsymbol{f}_n-\boldsymbol{I}\|_{L^\infty(\Omega;\RNN)}<\gamma, \qquad \textstyle  \frac{1}{2}<\|\det D \boldsymbol{f}_n\|_{L^\infty(\Omega)}<\frac{3}{2} \qquad \text{for $n\gg 1$.}
\end{equation}
Thus, using \eqref{eq:stress-control}--\eqref{eq:fn} and the change-of-variable formula, for $n\gg 1$ we estimate 
\begin{equation*}
	\begin{split}
		\sum_{\boldsymbol{a}\in A} \int_{ A(\boldsymbol{a},\varepsilon_n,2\varepsilon_n)\EEE} W(D\boldsymbol{y}_n)\,\d\boldsymbol{x}&= \sum_{\boldsymbol{a}\in A} \int_{B(\boldsymbol{a},2\varepsilon_n)} W((D\boldsymbol{y})\circ \boldsymbol{f}_n(D\boldsymbol{f}_n))\,\d\boldsymbol{x}\\
		&\leq \frac{N}{N-1} \sum_{\boldsymbol{a}\in A} \int_{ A(\boldsymbol{a},\varepsilon_n,2\varepsilon_n)\EEE}  W((D\boldsymbol{y})\circ \boldsymbol{f}_n)\,\d\boldsymbol{x}  + \frac{c_0 \omega_N 2^N N M}{N-1} \varepsilon_n^N\\
		&\leq   \frac{2N}{N-1} \sum_{\boldsymbol{a}\in A} \int_{ A(\boldsymbol{a},\varepsilon_n,2\varepsilon_n)\EEE}  W((D\boldsymbol{y})\circ \boldsymbol{f}_n) \det D \boldsymbol{f}_n\,\d\boldsymbol{x} + \frac{c_0 \omega_N 2^N N M}{N-1} \varepsilon_n^N\\
		&\leq \frac{2N}{N-1} \sum_{\boldsymbol{a}\in A} \int_{ A(\boldsymbol{a},\varepsilon_n,2\varepsilon_n)\EEE}  W(D\boldsymbol{y})\,\d\boldsymbol{x} + \frac{c_0 \omega_N 2^N N M}{N-1} \varepsilon_n^N,
	\end{split}
\end{equation*}
where we recall   $\omega_N = \leb(B(\boldsymbol{0},1))$. The right-hand side of the previous equation tends to zero, as $n\to \infty$. Therefore, recalling that $\boldsymbol{y}_n=\boldsymbol{y}$ in $\Omega^{A_n}_{2\varepsilon_n}$, we obtain 
\begin{equation*}
	\lim_{n\to \infty} \mathcal{W}_{\varepsilon_n}(A_n,\boldsymbol{y}_n) = \lim_{n\to \infty} \int_{\Omega^{A_n}_{2\varepsilon_n}} W(D\boldsymbol{y})\,\d\boldsymbol{x} +  \lim_{n\to \infty}  \sum_{\boldsymbol{a}\in A_n} \int_{A(\boldsymbol{a},\varepsilon_n,2\varepsilon_n)} W(D \boldsymbol{y}_n \EEE)\,\d\boldsymbol{x}=\mathcal{W}(\boldsymbol{y})
\end{equation*} 
thanks to the monotone convergence theorem.

Next, note that \eqref{eq:imt-recovery}  yields 
\begin{equation*}
	\imt(\boldsymbol{y}_n,B(\boldsymbol{a},\varepsilon_n))=\imt(\boldsymbol{y},B(\boldsymbol{a},r_n)) \quad \text{for all $\boldsymbol{a}\in A$ and all $n\in\N$.}
\end{equation*}
Thus,   in view of \eqref{eq:imtop-point-int2},  
\begin{equation*}
	\lim_{n\to \infty} \leb(\imt(\boldsymbol{y}_n,B(\boldsymbol{a},\varepsilon_n)))=	\lim_{n\to \infty}  \leb \left(  \imt(\boldsymbol{y},B(\boldsymbol{a},r_n)) \right)=\leb(\imt(\boldsymbol{y},\boldsymbol{a}))
\end{equation*} 
and, thanks to assumption \eqref{eq:conv-perimeter}, also
\begin{equation*}
	\lim_{n\to \infty} \per(\imt(\boldsymbol{y}_n,B(\boldsymbol{a},\varepsilon_n)))=	\lim_{n\to \infty}  \per \left(  \imt(\boldsymbol{y},B(\boldsymbol{a},r_n)) \right)=\per(\imt(\boldsymbol{y},\boldsymbol{a})).
\end{equation*} 
The last two equations immediately give
\begin{equation*}
	\lim_{n\to \infty} \mathcal{V}_{\varepsilon_n}(A_n,\boldsymbol{y}_n)=\mathcal{V}(A,\boldsymbol{y}), \qquad \lim_{n\to \infty} \mathcal{P}_{\varepsilon_n}(A_n,\boldsymbol{y}_n)=\mathcal{P}(A,\boldsymbol{y}),
\end{equation*}
which concludes the proof.
\end{proof}

%\begin{comment}

\section{Examples}
\label{sec:example}

We conclude by presenting a few examples of deformations  related to   assumption \eqref{eq:conv-perimeter}. In what follows, we will employ the $q$-norm on $\RN$  defined for $q\in [1,\infty]$ as
\begin{equation*}
	|\boldsymbol{x}|_q\coloneqq \begin{cases}
		\left ( \sum_{i=1}^N |x_i|^q \right )^{\frac{1}{q}} & \text{if $q\in [1,\infty)$,} \\
		\displaystyle \max_{i=1,\dots,N} |x_i| & \text{if $q=\infty$.}
	\end{cases}
\end{equation*}
The corresponding ball and sphere are denoted by $B_q(\boldsymbol{a},r)$ and $S_q(\boldsymbol{a},r)$,  respectively,  with $\boldsymbol{a}\in\RN$ and $r>0$. We use $\closure{B}_q(\boldsymbol{a},r)$ for the closure of $B_q(\boldsymbol{a},r)$ and we set $A_q(\boldsymbol{a},r,R)\coloneqq B_q(\boldsymbol{a},R)\setminus \closure{B}_q(\boldsymbol{a},r)$ for $R>r$. Similarly, $\closure{A}_q(\boldsymbol{a},r,R)$ stands for the closure of $A_q(\boldsymbol{a},r,R)$.  In agreement with the previously employed notation, we omit the subscript when $q=2$. \EEE

Radial deformations with respect to the Euclidean norm are trivially shown to fulfill \eqref{eq:conv-perimeter}. The first example shows an instance where this holds true also for another norm corresponding to a cavity of different shape.

\begin{center}
	\begin{figure}[t]
		\begin{tikzpicture}[scale=2.5]
			\draw[xshift=-40pt,fill=black!20] (0,1) -- (1,0) -- (0,-1) -- (-1,0) -- (0,1);
			\draw[xshift=-40pt,fill=black] (0,0) circle (.5pt);
			\node[xshift=-40pt] at (-.85,-1.2) {$\Omega$};
			\draw[xshift=40pt,fill=black!20] (0,1) -- (1,0) -- (0,-1) -- (-1,0) -- (0,1);
			\def\b{.33};
			\node[color=blue] at (-1.2,-.35) {\footnotesize $S(\boldsymbol{0},r)$};
			\node[color=blue] at (1.43,-.62) {\footnotesize $\boldsymbol{y}(S(\boldsymbol{0},r))$};
			\draw[xshift=40pt,fill=white] (0,\b) -- (\b,0) -- (0,-\b) -- (-\b,0) -- (0,\b);
			\draw[thick,->] (-.5,.5) to [out=30,in=150] node[above,midway] {$\boldsymbol{y}$} (.5,.5);
			\node[xshift=40pt] at (0.8,-1.2) {$\boldsymbol{y}(\Omega \setminus \{ \boldsymbol{0} \})$};
			\usetikzlibrary{calc}
			\usetikzlibrary{math}
			\def\r{.25};
			\draw[xshift=-40pt,thick,color=blue] (0,0) circle (\r);
			\draw[xshift=40pt,thick,color=blue] plot [smooth, samples=100, domain=0:90, variable=\t] ({(1-\b)*\r*cos(\t)+((\b*cos(\t))/(cos(\t)+sin(\t)))},{(1-\b)*\r*sin(\t)+((\b*sin(\t))/(cos(\t)+sin(\t)))});
			\draw[xshift=40pt,thick,color=blue] plot [smooth, samples=100, domain=0:90, variable=\t] ({(1-\b)*\r*cos(\t)+((\b*cos(\t))/(cos(\t)+sin(\t)))},{-(1-\b)*\r*sin(\t)-((\b*sin(\t))/(cos(\t)+sin(\t)))});
			\draw[xshift=40pt,thick,color=blue] plot [smooth, samples=100, domain=0:90, variable=\t] ({-(1-\b)*\r*cos(\t)-((\b*cos(\t))/(cos(\t)+sin(\t)))},{(1-\b)*\r*sin(\t)+((\b*sin(\t))/(cos(\t)+sin(\t)))});
			\draw[xshift=40pt,thick,color=blue] plot [smooth, samples=100, domain=0:90, variable=\t] ({-(1-\b)*\r*cos(\t)-((\b*cos(\t))/(cos(\t)+sin(\t)))},{-(1-\b)*\r*sin(\t)-((\b*sin(\t))/(cos(\t)+sin(\t)))});
		\end{tikzpicture}
		\caption{The deformation in Example \ref{ex:radial} for $b=1/2$.}
		\label{fig:square-cavity}
	\end{figure}
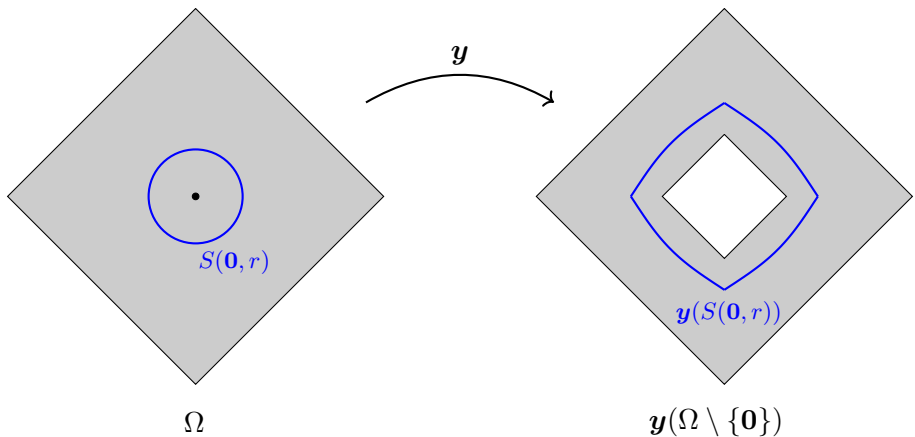
\end{center}
 
\begin{example}[Radial deformation]\label{ex:radial}
Let $N=2$, $p>1$, and $b\in (0,1)$. Given $\Omega\coloneqq B_1(\boldsymbol{0},1)$, we define $\boldsymbol{y}\colon \Omega \to \R^2$ by setting
\begin{equation*}
	\label{eq:radial}
	\boldsymbol{y}(\boldsymbol{x})\coloneqq ((1-b)|\boldsymbol{x}|_1+b) \frac{\boldsymbol{x}}{|\boldsymbol{x}|_1}\quad \text{for all $\boldsymbol{x}\neq \boldsymbol{0}$.}
\end{equation*}
In view of the results in \cite[Appendix]{bresciani.friedrich.moracorral}, we have $\boldsymbol{y}\in\mathcal{Y}^A$ for $A\coloneqq \{ \boldsymbol{0} \}$ with $\imt(\boldsymbol{y},\boldsymbol{0})=\closure{B}_1(\boldsymbol{0},b)$.
Let $r\in (0,1)$.  Then,  $\per \left( \imt(\boldsymbol{y},B(\boldsymbol{0},r)) \right)=\H^1(\boldsymbol{y}(S(\boldsymbol{0},r)))$ thanks to \cite[Lemma 2.28]{bresciani.friedrich.moracorral}. \EEE  By symmetry, we have 
\begin{equation*}
	\H^1(\boldsymbol{y}(S(\boldsymbol{0},r)))=4\H^1(\boldsymbol{y}(S(\boldsymbol{0},r))\cap Q),
\end{equation*}
where $Q$ denotes the first quadrant of the Cartesian plane. To compute the length on the right-hand side, we employ polar coordinates. By writing $\boldsymbol{v}_r(t)\coloneqq \boldsymbol{y}\left ( r(\cos t, \sin t)^\top \right )$ for  $t\in (0,\pi/2)$, we find
\begin{equation*}
	\boldsymbol{v}_r(t)= \left( (1-b)r+ \frac{b}{\cos t+  \sin t} \right) (\cos t, \sin t)^\top.
\end{equation*}   
Then,
\begin{equation*}
	 \H^1(\boldsymbol{y}(S(\boldsymbol{0},r)\cap Q)) =\int_0^{\pi/2} |\boldsymbol{v}_r'(t)|\,\d t,
\end{equation*}
where  
\begin{equation*}
\begin{split}
		|\boldsymbol{v}'_r(t)|=\sqrt{  \frac{2b^2  }{(\cos t + \sin t)^4}  + \O(r)  }=\frac{  \sqrt{2}  b}{(\cos t+\sin t)^2} + \O(r), \quad \text{as $r\to 0^+$,} \quad \text{for all $t\in (0,\pi/2)$.}
	\end{split}
\end{equation*}
Thus, by letting $r\to 0^+$, with the aid of the dominated convergence theorem, we obtain
\begin{equation*}
	\lim_{r\to 0^+} \H^1(\boldsymbol{y}(S(\boldsymbol{0},r))=4 \int_0^{\pi/2} \frac{\sqrt{2}b}{(\cos t + \sin t)^2} \d t=4\sqrt{2}b=\H^1(S_1(\boldsymbol{0},b))=\per \left( \imt(\boldsymbol{y},\boldsymbol{0}) \right),
\end{equation*}
 which shows that $\boldsymbol{y}$ fulfills \eqref{eq:conv-perimeter}. \EEE 
%Here, we noted that
%\begin{equation*}
%		\int_0^{\pi/2} \frac{1}{(\cos t + \sin t)^2} \d t=\int_0^{\pi/2} \frac{1}{1+\sin(2t)}\,\d t=\int_0^{\pi/2}\frac{1-\sin(2t)}{\cos^2(2t)}\,\d t=1.
%\end{equation*}
\end{example}

In the next example, the reference configuration is stretched before opening a radial cavity. This touches the question of the possible dependence of assumption \eqref{eq:conv-perimeter} on the choice of the reference configuration. 

\begin{figure}
	\begin{tikzpicture}[scale=1.6]
		\usetikzlibrary{calc,math};
		\draw[xshift=-90pt,fill=black!20] (-1,-1) rectangle (1,1);
		\def\r{.35};
		\draw[xshift=-90pt,fill=black] (0,0) circle (.5pt);
		\draw[xshift=-90pt,color=blue,thick] (0,0) circle (\r);
		\node[color=blue] at (-2.7,-.5) {\footnotesize $S(\boldsymbol{0},r)$};
		\node at (-3.2,-1.2) {$\Omega$};
		\node at (.6,-1.2) {$\boldsymbol{f}(\Omega)$};
		\node at (4.7,-1.2) {$\boldsymbol{y}(\Omega \setminus {\boldsymbol{0}})$};
		\draw[thick,->] (-2,0) to [out=30,in=150] node[above,midway] {$\boldsymbol{f}$} (-1.2,0);
		\draw[thick,->] (-2,.9) to [out=30,in=150] node[above,midway] {$\boldsymbol{y}$} (3,.9);
		\draw[thick,->] (2.2,0) to [out=30,in=150] node[above,midway] {$\boldsymbol{u}$} (3,0);
		\draw[fill=black!20] (-1,-1) rectangle (2,1);
		\draw[fill=black] (0,0) circle (.5pt);
		\node[color=blue] at (.8,-.5) {\footnotesize $\boldsymbol{f}(S(\boldsymbol{0},r))$};
		\draw[color=blue,thick] plot [smooth,samples=100,domain=-90:90,variable=\t] ({2*\r*cos(\t)},{\r*sin(\t)});
		\draw[color=blue,thick] plot [smooth,samples=100,domain=-90:90,variable=\t] ({-\r*cos(\t)},{\r*sin(\t)});
		\draw[xshift=120pt,fill=black!20] (-1,-1) rectangle (2,1);
		\def\b{.55};
		\draw[xshift=120pt,fill=white] (0,0) circle (\b);
		\draw[xshift=120pt,color=blue,thick] plot [smooth,samples=100,domain=-90:90,variable=\t] ({-((1-\b)*\r+\b)*cos(\t)},{((1-\b)*\r+\b)*sin(\t)});
		\draw[xshift=120pt,color=blue,thick] plot [smooth,samples=100,domain=-90:90,variable=\t] ({2*(1-\b)*\r*cos(\t)+((2*\b*cos(\t))/(sqrt(3*(cos(\t))^2+1)))},{(1-\b)*\r*sin(\t)+((\b*sin(\t))/(sqrt(3*(cos(\t))^2+1)))});
		\node[color=blue] at (5.4,-.5) {\footnotesize $\boldsymbol{y}(S(\boldsymbol{0},r))$};
	\end{tikzpicture}
	\caption{The deformation in Example \ref{ex:cor} for $b=1/2$.}
	\label{fig:chofref}
\end{figure}
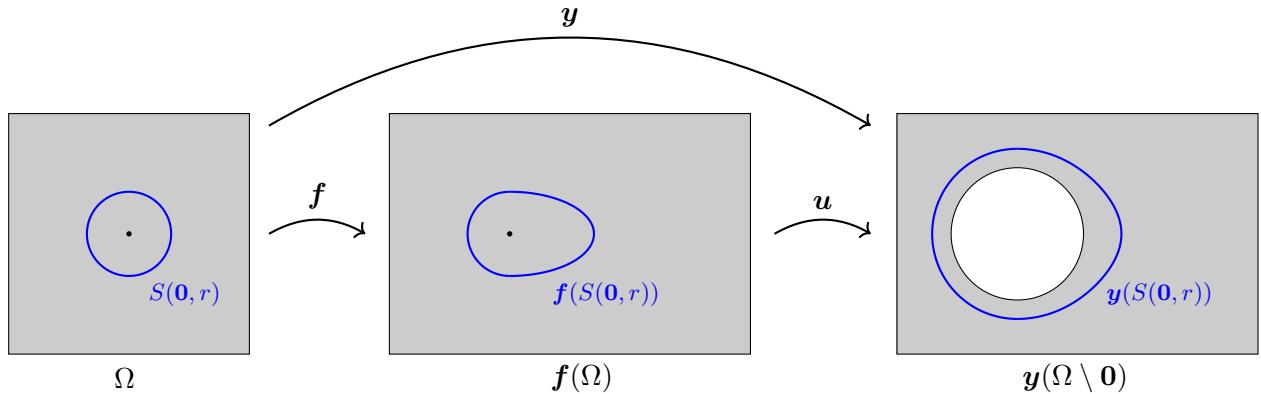

\begin{example}[Change of reference configuration]\label{ex:cor}
Let $N=2$, $p>1$, and $b\in (0,1)$. Given $\Omega\coloneqq B_\infty(\boldsymbol{0},1)=(-1,1)^2$, consider the piecewise affine map $\boldsymbol{f}\colon \Omega \to \R^2$ defined as
\begin{equation*}
	\boldsymbol{f}(\boldsymbol{x})\coloneqq \begin{cases}
		(2x_1,x_2)^\top & \text{if $x_1\geq 0$,}\\
		\boldsymbol{x} & \text{if $x_1 <  0$.}
	\end{cases}
\end{equation*}
Thus, $\boldsymbol{f}(\Omega)=(-1,2)^2$. Define $\boldsymbol{u}\colon (-1,2)^2\to \R^2$ as
\begin{equation*}
	\boldsymbol{u}(\boldsymbol{z})\coloneqq \begin{cases}
		((1-b)|\boldsymbol{z}|+b)\frac{\boldsymbol{z}}{|\boldsymbol{z}|} & \text{if $0<|\boldsymbol{z}|<1$,}\\
		\boldsymbol{z} & \text{if $|\boldsymbol{z}| \ge  1$.}
	\end{cases}
\end{equation*}
Then, set $\boldsymbol{y}\coloneqq \boldsymbol{u}\circ \boldsymbol{f}$. In view of the results in \cite[Appendix]{bresciani.friedrich.moracorral}, we deduce $\boldsymbol{y}\in\mathcal{Y}^A$ with $A\coloneqq \{ \boldsymbol{0}\}$ and $\imt(\boldsymbol{y},\boldsymbol{0})=\closure{B}(\boldsymbol{0},b)$.
Let $r\in (0,1)$   and note that $\per \left( \imt(\boldsymbol{y},B(\boldsymbol{0},r)) \right)=\H^1(\boldsymbol{y}(S(\boldsymbol{0},r)))$ thanks to \cite[Lemma 2.28]{bresciani.friedrich.moracorral}. \EEE By symmetry
\begin{equation*}
	\H^1(\boldsymbol{y}(S(\boldsymbol{0},r)))=\pi ((1-b)r+b) + 2 \H^1(\boldsymbol{y}(S(\boldsymbol{0},r)\cap Q)),
\end{equation*}
where  $Q$ denotes the first quadrant of the Cartesian plane.
To compute $\H^1(\boldsymbol{y}(S(\boldsymbol{0},r)\cap Q))$ we use polar coordinates. By writing $\boldsymbol{v}_r(t)\coloneqq \boldsymbol{y}(r(\cos t, \sin t)^\top)$ for $t\in (0,\pi/2)$, we have
\begin{equation*}
	\boldsymbol{v}_r(t)=\left( (1-b)r + \frac{b}{\sqrt{3\cos^2 t +1}}  \right) (2\cos t, \sin t)^\top
\end{equation*}
and
\begin{equation*}
	\H^1(\boldsymbol{y}(S(\boldsymbol{0},r)\cap Q))=\int_0^{\pi/2} |\boldsymbol{v}_r'(t)|\,\d t.
\end{equation*}
A direct computation gives
\begin{equation*}
	|\boldsymbol{v}_r'(t)|=\sqrt{\frac{4b^2}{(3\cos^2 t +1)^2}+\O(r)}=\frac{2b}{3\cos^2t+1}+\O(r), \quad \text{as $r\to 0^+$,}\quad \text{for all $t\in (0,\pi/2)$.}
\end{equation*}
Thus, letting $r\to 0^+$, with the aid of the dominated convergence theorem, we find
\begin{equation*}
	\lim_{r\to 0^+}\H^1(\boldsymbol{y}(S(\boldsymbol{0},r)\cap Q))=2b\int_0^{\pi/2} \frac{1}{3 \cos^2 t + 1} dt=\frac{\pi b}{2}
\end{equation*}
and
\begin{equation*}
		\lim_{r\to 0^+}\H^1(\boldsymbol{y}(S(\boldsymbol{0},r))=	\lim_{r\to 0^+} \left (\pi ((1-b)r+b) + \pi b \right)= 2\pi b =\H^1(S(\boldsymbol{0},b))=\per \left(\imt(\boldsymbol{y},\boldsymbol{0})\right).
\end{equation*}
 Thus, $\boldsymbol{y}$ satisfies \eqref{eq:conv-perimeter}. 
\end{example}

In the  third   example, we consider the superposition of a radial deformation with a piecewise smooth map. The resulting deformation constitutes the main building block of the famous example in \cite[Section 11]{mueller.spector}.

\begin{figure}
	\begin{tikzpicture}[scale=1.8]
		\draw[xshift=-90pt,fill=black!20] (-1,-1) rectangle (1,1);
		\draw[xshift=-90pt,fill=black] (0,0) circle (.5pt);
		\def\r{.3};
		\draw[xshift=-90pt,color=blue,thick] (0,0) circle (\r);
		\draw[xshift=0pt,fill=black!20] (-1,-1) rectangle (1,1);
		\draw[->,thick] (-2,0) to[out=30,in=150] node[above,midway] {$\boldsymbol{u}$} (-1.2,0);
		\draw[->,thick] (1.2,0) to[out=30,in=150] node[above,midway] {$\boldsymbol{g}$} (2,0);
		\draw[->,thick] (-2.1,.9) to[out=30,in=150] node[above,midway] {$\boldsymbol{y}$} (2.1,.9);
		\draw[color=blue,thick] plot [smooth,samples=30,domain=-45:45,variable=\t] ({.5*\r*cos(\t)+.5},{.5*\r*sin(\t)+.5*tan(\t)});
		\draw[color=blue,thick] plot [smooth,samples=30,domain=-45:45,variable=\t] ({.5*\r*sin(\t)+.5*tan(\t)},{.5*\r*cos(\t)+.5});
		\draw[color=blue,thick] plot [smooth,samples=30,domain=-45:45,variable=\t] ({-.5*\r*cos(\t)-.5},{.5*\r*sin(\t)+.5*tan(\t)});
		\draw[color=blue,thick] plot [smooth,samples=30,domain=-45:45,variable=\t] ({.5*\r*sin(\t)+.5*tan(\t)},{-.5*\r*cos(\t)-.5});
		\draw[xshift=0pt,fill=white] (-.5,-.5) rectangle (.5,.5);
		\draw[xshift=90pt,fill=black!20] (-1,-1) rectangle (1,1);
		\draw[xshift=90pt,color=blue,thick] plot[smooth,samples=30,domain=0:42.7,variable=\t] ({1-.5*tan(\t)+.5*\r^2*sin(\t)*cos(\t)},{.5*\r*sin(\t)+.5*tan(\t)});
		\draw[xshift=90pt,color=blue,thick] plot[smooth,samples=30,domain=0:42.7,variable=\t] ({-1+.5*tan(\t)-.5*\r^2*sin(\t)*cos(\t)},{.5*\r*sin(\t)+.5*tan(\t)});
		\draw[xshift=90pt,color=blue,thick] plot[smooth,samples=30,domain=0:42.7,variable=\t] ({1-.5*tan(\t)+.5*\r^2*sin(\t)*cos(\t)},{-.5*\r*sin(\t)-.5*tan(\t)});
		\draw[xshift=90pt,color=blue,thick] plot[smooth,samples=30,domain=0:42.7,variable=\t] ({-1+.5*tan(\t)-.5*\r^2*sin(\t)*cos(\t)},{-.5*\r*sin(\t)-.5*tan(\t)});
		\draw[xshift=90pt,color=blue,thick] plot[smooth,samples=30,domain=0:42.7,variable=\t] ({.5*\r*sin(\t)+.5*tan(\t)},{1-.5*tan(\t)+.5*\r^2*sin(\t)*cos(\t)});
		\draw[xshift=90pt,color=blue,thick] plot[smooth,samples=30,domain=0:42.7,variable=\t] ({-.5*\r*sin(\t)-.5*tan(\t)},{1-.5*tan(\t)+.5*\r^2*sin(\t)*cos(\t)});
		\draw[xshift=90pt,color=blue,thick] plot[smooth,samples=30,domain=0:42.7,variable=\t] ({.5*\r*sin(\t)+.5*tan(\t)},{-1+.5*tan(\t)-.5*\r^2*sin(\t)*cos(\t)});
		\draw[xshift=90pt,color=blue,thick] plot[smooth,samples=30,domain=0:42.7,variable=\t] ({-.5*\r*sin(\t)-.5*tan(\t)},{-1+.5*tan(\t)-.5*\r^2*sin(\t)*cos(\t)});
		\draw[xshift=90pt,fill=black!20,color=black!20] (.5,.5) rectangle (.9,.9);
		\draw[xshift=90pt,fill=black!20,color=black!20] (.51,-.51) rectangle (.9,-.9);
		\draw[xshift=90pt,fill=black!20,color=black!20] (-.51,.51) rectangle (-.9,.9);
		\draw[xshift=90pt,fill=black!20,color=black!20] (-.51,-.51) rectangle (-.9,-.9);
		\draw[xshift=90pt,fill=white] (0,1) -- (1,0) -- (0,-1) -- (-1,0) -- (0,1);
		\draw[color=blue,thick,xshift=90pt] plot [smooth,samples=30,domain=38.5:45,variable=\t] ({.5*\r*cos(\t)+.5},{.5*\r*sin(\t)+.5*tan(\t)});
		\draw[color=blue,thick,xshift=90pt] plot [smooth,samples=30,domain=38.5:45,variable=\t] ({-.5*\r*cos(\t)-.5},{.5*\r*sin(\t)+.5*tan(\t)});
		\draw[color=blue,thick,xshift=90pt] plot [smooth,samples=30,domain=38.5:45,variable=\t] ({-.5*\r*cos(\t)-.5},{-.5*\r*sin(\t)-.5*tan(\t)});
		\draw[color=blue,thick,xshift=90pt] plot [smooth,samples=30,domain=-45:-38.5,variable=\t] ({.5*\r*cos(\t)+.5},{.5*\r*sin(\t)+.5*tan(\t)});
		\draw[color=blue,thick,xshift=90pt] plot [smooth,samples=30,domain=38.5:45,variable=\t] ({.5*\r*sin(\t)+.5*tan(\t)},{.5*\r*cos(\t)+.5});
		\draw[color=blue,thick,xshift=90pt] plot [smooth,samples=30,domain=-45:-38.5,variable=\t] ({.5*\r*sin(\t)+.5*tan(\t)},{.5*\r*cos(\t)+.5});
		\draw[color=blue,thick,xshift=90pt] plot [smooth,samples=30,domain=-45:-38.5,variable=\t] ({.5*\r*sin(\t)+.5*tan(\t)},{-.5*\r*cos(\t)-.5});
		\draw[color=blue,thick,xshift=90pt] plot [smooth,samples=30,domain=-45:-38.5,variable=\t] ({-.5*\r*sin(\t)-.5*tan(\t)},{-.5*\r*cos(\t)-.5});
		\node at (-3.2,-1.2) {$\Omega$};
		\node at (0,-1.2) {$\boldsymbol{u}(\Omega \setminus \{ \boldsymbol{0}\})$}; 
		\node at (3.2,-1.2) {$\boldsymbol{y}(\Omega \setminus \{ \boldsymbol{0}\})$};
		\node[color=blue] at (-2.8,-.4) {\footnotesize $S(\boldsymbol{0},r)$};
		\node[color=blue] at (.4,-.8) {\footnotesize $\boldsymbol{u}(S(\boldsymbol{0},r))$};
		\node[color=blue] at (3.8,-.87) {\scriptsize $\boldsymbol{y}(S(\boldsymbol{0},r))$};
	\end{tikzpicture}
	\caption{The deformation in Example \ref{ex:superpos}.}
	\label{fig:superpos}
\end{figure}
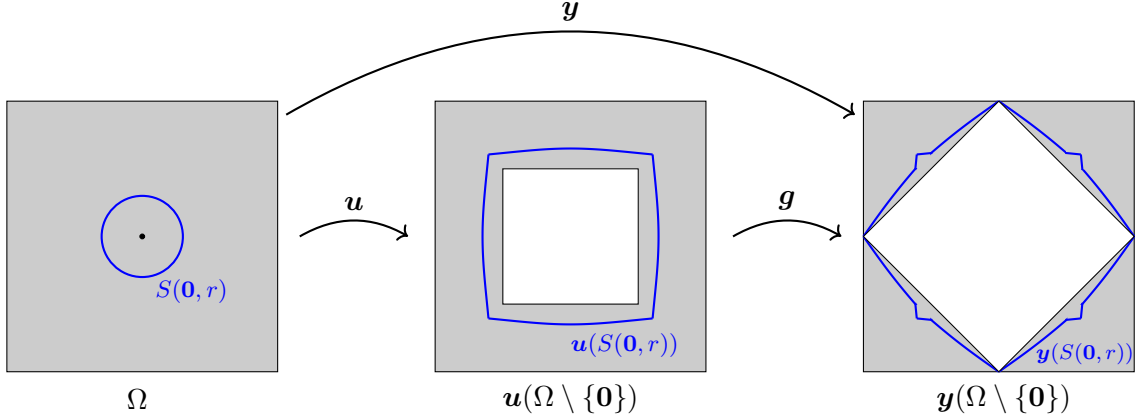

\begin{example}[Superposition]\label{ex:superpos}
Let $N=2$  and  $\Omega\coloneqq B_\infty(\boldsymbol{0},1)$. Define $\boldsymbol{u}\colon \Omega \to \R^2$ as 
\begin{equation*}
	\boldsymbol{u}(\boldsymbol{x})\coloneqq \frac{1}{2} (|\boldsymbol{x}|_\infty+1)\frac{\boldsymbol{x}}{|\boldsymbol{x}|_\infty} \quad \text{for all $\boldsymbol{x}\neq \boldsymbol{0}$.}
\end{equation*}
Thus, $\boldsymbol{u}(B_\infty(\boldsymbol{0},1)\setminus \{ \boldsymbol{0} \})=A_\infty(\boldsymbol{0},1/2,1)$. Define $\boldsymbol{g}\colon \closure{A}_\infty(\boldsymbol{0},1/2,1)\to \R^2$ as
\begin{equation*}
	\boldsymbol{g}(\boldsymbol{z})\coloneqq \begin{cases}
		\left( \mathrm{sgn}(z_1)(1-2|z_2|)+2z_1|z_2|,z_2 \right)^\top & \text{if $|z_1|>|z_2|$ and $|z_2|<\frac{1}{2}$,}\\
		\left( z_1, \mathrm{sgn}(z_2)(1-2|z_1|)+2|z_1|z_2 \right)^\top & \text{if $|z_2|>|z_1|$ and $|z_1 |<\frac{1}{2}$,}\\
		\boldsymbol{z} & \text{otherwise.} 
	\end{cases}
\end{equation*}
Set $\boldsymbol{y}\coloneqq \boldsymbol{g}\circ \boldsymbol{u}$. Thus, $\boldsymbol{y}\in\mathcal{Y}^A$  with $A\coloneqq \left\{ \boldsymbol{0}  \right\}$ for any $1<p<2$ \EEE and $\imt(\boldsymbol{y},\boldsymbol{0})=  \closure{B}_1(\boldsymbol{0},1)  $ in view of \cite[Appendix]{bresciani.friedrich.moracorral}.
Let $r\in (0,1)$.   Then,  $\per \left( \imt(\boldsymbol{y},B(\boldsymbol{0},r)) \right)=\H^1(\boldsymbol{y}(S(\boldsymbol{0},r)))$ in view of \cite[Lemma 2.28]{bresciani.friedrich.moracorral}.  By symmetry 
\begin{equation*}
	\H^1(\boldsymbol{y}(S(\boldsymbol{0},r)))=8\H^1(\boldsymbol{y}(S(\boldsymbol{0},r))\cap E),
\end{equation*}
where $E\coloneqq \left\{ \boldsymbol{\xi}\in\R^2:\:\xi_1>\xi_2>0  \right\}$. Also here we employ polar coordinates. Let $\boldsymbol{v}_r\colon [0,\pi/4] \to \R^2$ be defined as $\boldsymbol{v}_r(t)\coloneqq \boldsymbol{y}(r(\cos t, \sin t)^\top)$. Then
\begin{equation*}
	\H^1(\boldsymbol{y}(S(\boldsymbol{0},r))\cap E)=\int_0^{\pi/4} |\boldsymbol{v}'_r(t)|\,\d t.
\end{equation*} 
Consider the equation $r\sin \vartheta + \tan \vartheta=1$ in that interval $(0,\pi/4)$. It is easy to show that this equation has a unique solution $\vartheta_r\in (0,\pi/4)$ which satisfies $\vartheta_r \to \pi/4$, as $r\to 0^+$. By construction, we have 
\begin{equation*}
	\boldsymbol{v}_r(t)=\begin{cases}
		\left( 1-\frac{1}{2}\tan t+\frac{r^2}{2}\sin t \cos t,\frac{r}{2}\sin t+\frac{1}{2}\tan t   \right)^\top & \text{if $t\in [0,\vartheta_r]$,}\\
		\frac{1}{2} \left(r+\frac{1}{\cos t}\right) (\cos t,\sin t)^\top & \text{if $t\in [\vartheta_r,\pi/4]$.}
	\end{cases} 
\end{equation*}
Hence,
\begin{equation*}
	|\boldsymbol{v}_r'(t)|=\begin{cases}
		\frac{1}{\sqrt{2}\cos^2t}+\O(r) & \text{if $t\in [0,\vartheta_r]$,}\\
		\frac{1}{2\cos^2t}+\O(r) & \text{if $t\in [\vartheta_r,\pi/4]$.}
	\end{cases}
\end{equation*}
Letting $r\to 0^+$, with the aid of the dominated convergence theorem, we obtain
\begin{equation*}
	\begin{split}
		\lim_{r\to 0^+} \H^1(\boldsymbol{y}(S(\boldsymbol{0},r)))&=8 \lim_{r\to 0^+} \left( \int_0^{ \vartheta_r}|\boldsymbol{v}_r'(t)|\,\d t + \int_{\vartheta_r}^{\pi/4} |\boldsymbol{v}_t'(t)|\,\d t  \right)=\frac{8}{\sqrt{2}} \int_0^{\pi/4} \frac{1}{\cos^2 t}\,\d t= \frac{8}{\sqrt{2}} \EEE \\
		&=  \H^1(S_1(\boldsymbol{0},1))=\per \left(\imt(\boldsymbol{y},\boldsymbol{0})\right).  \EEE
	\end{split}
\end{equation*}
 This shows that $\boldsymbol{y}$ fulfills \eqref{eq:conv-perimeter}.  
\end{example}

  Eventually,  we come to an example where   \eqref{eq:conv-perimeter} is violated. 

\begin{figure}
	\begin{tikzpicture}[scale=2]
		\usetikzlibrary{math}
		\usetikzlibrary{calc}
		\def\r{.2};
		\def\R{.6}; %.5*(\r+1)
		%0
		\draw[->,thick] (-2,0) to[out=30,in=150] node[above,midway] {$\boldsymbol{u}$} (-1.2,0);
		\draw[->,thick] (-2.1,.8) to[out=30,in=150] node[above,midway] {$\boldsymbol{y}$} (2.1,.8);
		\draw[xshift=-90pt,fill=black!20] (0,0) circle (1);
		\draw[xshift=-90pt,fill=black] (0,0) circle (.5pt);
		\draw[xshift=-90pt,color=blue,thick] (0,0) circle (\r);
		\node[color=blue] at (-2.8,-.35) {\footnotesize $S(\boldsymbol{0},r)$};
		\node at (-3.2,-1.2) { $\Omega$};
		%1
		\draw[->,thick] (1.2,0) to[out=30,in=150] node[above,midway] {$\boldsymbol{g}$} (2,0);
		\draw[xshift=0pt,fill=black!20] (0,0) circle (1);
		\draw[xshift=0pt,fill=white] (0,0) circle (.5);
		\draw[xshift=0pt,dotted] (0,.5) -- (.51,.862);
		\draw[xshift=0pt,dotted] (0,.5) -- (-.51,.862);
		\draw[xshift=0pt,color=blue,thick] (0,0) circle (\R);
		\draw[color=blue,fill=blue] (.1285,.583) circle (.5pt);
		\draw[color=blue,fill=blue] (-.1285,.583) circle (.5pt);
		\node[color=blue,right] at (.2,.6) {\footnotesize $\boldsymbol{w}_r^+$};
		\node[color=blue,left] at (-.2,.6) {\footnotesize $\boldsymbol{w}_r^-$};
		\node[above] at (.33,.68) {\footnotesize $\ell^+$};
		\node[above,color=blue] at (.12,.57) {\footnotesize $\gamma_r^+$};
		\node[above,color=blue] at (-.12,.58) {\footnotesize $\gamma_r^-$};
		\node[above] at (-.27,.7) {\footnotesize $\ell^-$};
		\node[color=blue] at (.2,-.75) {\footnotesize $\boldsymbol{u}(S(\boldsymbol{0},r))$};
		\node at (0,-1.2) {$\boldsymbol{u}(\Omega \setminus \{ \boldsymbol{0} \})$};
		\draw[xshift=0pt,dotted] (0,.5)--(0,1);
		\node at (.05,.9) {\footnotesize $\ell$};
		%2
		\draw[xshift=90pt,fill=black!20] (0,0) circle (1);
		\draw[xshift=90pt,fill=white] (0,0) circle (.5);
		\draw[xshift=90pt,color=blue,thick] (0,0) circle (\R);
		\draw[xshift=90pt,fill=white,color=white] (.55,.863)--(0,1.1)--(-.55,.863);
		\draw[xshift=90pt,color=black!20,fill=black!20] (0,.5)--(.5,.86)--(0,1)--(-.5,.86)--(0,.5);
		\draw[xshift=90pt] (-.03,.5)--(.03,.5);
		\draw[xshift=90pt] (.51,.862)-- (0,1)-- (-.51,.862);
		\node[color=blue] at (3.4,.8) {\footnotesize $\sigma_r^+$};
		\node[color=blue] at (2.95,.8) {\footnotesize $\sigma_r^-$};
		\draw[xshift=90pt,dotted] (0,.5) -- (.51,.862);
		\draw[xshift=90pt,dotted] (0,.5) -- (-.51,.862);
		\draw[xshift=90pt,color=blue,thick] (0,1)--(.1285,.583);
		\draw[xshift=90pt,color=blue,thick] (0,1)--(-.1285,.583);
		\node[color=blue] at (3.35,-.75) {\footnotesize $\boldsymbol{y}(S(\boldsymbol{0},r))$};
		\node[color=blue] at (3.75,.45) {\footnotesize $C_r$};
		\node at (3.2,-1.2) { $\boldsymbol{y}(\Omega \setminus \{ \boldsymbol{0} \})$};
		\draw[xshift=90pt] (0,.5)--(0,1);
	\end{tikzpicture}
	\caption{The deformation in Example \ref{ex:counterexample}.}
	\label{fig:counterexample}
\end{figure}
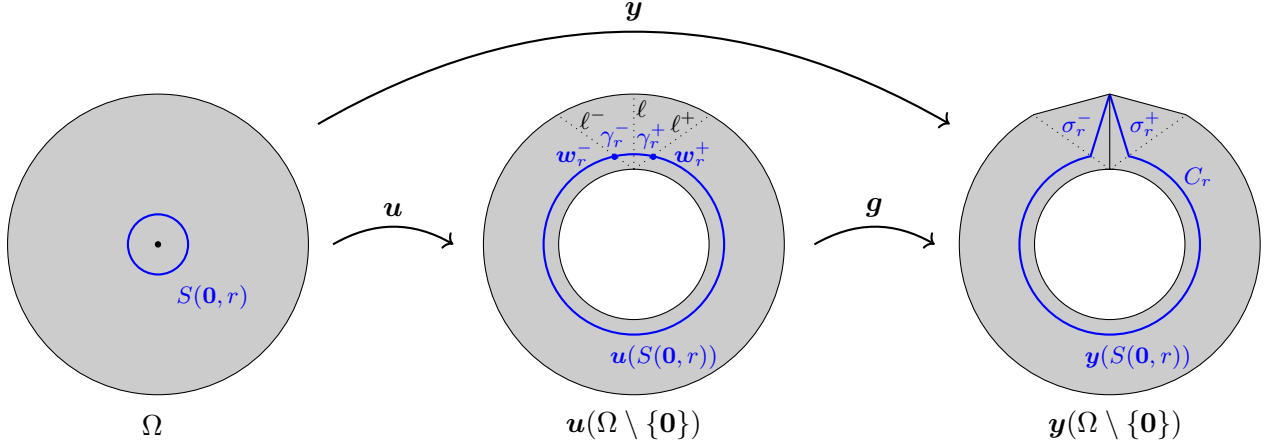

\begin{example}[Cavity with a spike]\label{ex:counterexample}
	Let $N=2$ and $\Omega\coloneqq B(\boldsymbol{0},1)$. Define $\boldsymbol{u}\colon \Omega \to \R^2$ by setting
	\begin{equation*}
		\boldsymbol{u}(\boldsymbol{x})\coloneqq \frac{1}{2} \left( |\boldsymbol{x}| +1  \right) \frac{\boldsymbol{x}}{|\boldsymbol{x}|} \quad \text{for all $\boldsymbol{x}\neq \boldsymbol{0}$.}
	\end{equation*}
	Thus, $\boldsymbol{u}(\Omega \setminus \{ \boldsymbol{0} \})=A(\boldsymbol{0},1/2,1)$. Set $V\coloneqq \left\{ \boldsymbol{z}\in A(\boldsymbol{0},1/2,1):\:z_2>(\sqrt{3}-1)|z_1|+\textstyle \frac{1}{2}  \right\}$.  The region $V$ is    bounded by the  two segments $\ell^\pm$ joining the points $(\pm 1/2,\sqrt{3}/2)\top$  with $(0,1/2)\top$ and the shorter of the two arcs on $S(\boldsymbol{0},1)$ connecting the two points $(\pm 1/2,\sqrt{3}/2)^\top$. We define $\boldsymbol{g}\colon \closure{A}(\boldsymbol{0},1/2,1)\to \R^2$ by setting  
	\begin{equation*}
		 \boldsymbol{g}(\boldsymbol{z})\coloneqq \begin{cases}
			\left (z_1, \frac{-9+4\sqrt{3}+(\sqrt{3}-1)\sqrt{1-4(5-\sqrt{3})|\boldsymbol{z}|^2}}{1-\sqrt{3}+\sqrt{1-4(5-\sqrt{3})|\boldsymbol{z}|^2}}
			 |z_1|+1  \right)^\top & \text{if $\boldsymbol{z}\in V$,}\\
			\boldsymbol{z}& \text{ if $\boldsymbol{z}\in \closure{A}(\boldsymbol{0},1/2,1)\setminus V$. }
		\end{cases}
	\end{equation*}
	We observe that $\boldsymbol{g}$ is a Lipschitz transformation which  maps the whole segment $\ell$ joining $(0,1/2)^\top$ and $(0,1)^\top$ to the point $(0,1)^\top$. In particular, $\boldsymbol{g}$ is injective on $\closure{A}(\boldsymbol{0},1/2,1)\setminus \ell$. 
	For every $r\in (0,1)$, denoting by $\boldsymbol{w}_r^\pm$ the intersection of $\boldsymbol{u}(S(\boldsymbol{0},r))=S(\boldsymbol{0},\frac{1}{2}(r+1))$ with $\ell^\pm$, the map $\boldsymbol{g}$ transforms  $\gamma_r^\pm$, the shorter of the two arcs   on $\boldsymbol{u}(S(\boldsymbol{0},r))$ connecting $\boldsymbol{w}_r^\pm$ to $(0,\frac{1}{2}(r+1))^\top$,   into the segment $\sigma_r^\pm$ connecting $\boldsymbol{w}_r^\pm$ with $(0,1)^\top$.

	We define $\boldsymbol{y}\coloneqq \boldsymbol{g}\circ \boldsymbol{u}$ and, from the results in \cite[Appendix]{bresciani.friedrich.moracorral}, we see that $\boldsymbol{y}\in\mathcal{Y}^A$ with $A\coloneqq \{ \boldsymbol{0} \}$ for any $1<p<2$. Moreover, $\imt(\boldsymbol{y},\boldsymbol{0})=\closure{B}(\boldsymbol{0},1/2) \cup \ell$. Let $r\in (0,1)$. If $C_r$ denotes the longer of the two arcs on  $S(\boldsymbol{0},\frac{1}{2}(r+1))$ connecting the points $\boldsymbol{w}_r^\pm$, we obtain
	\begin{equation*}
		\per \left( \imt(\boldsymbol{y},\boldsymbol{0})  \right)=\H^1(C_r)+\H^1(\sigma_r^+)+\H^1(\sigma_r^-).
	\end{equation*}
	 Observing that $\H^1(C_r)\to \pi$ and $\H^1(\sigma_r^\pm)\to 1/2$, as $r\to 0^+$, because $\boldsymbol{w}_r^\pm \to (0,1/2)^\top$, we find
	 \begin{equation*}
	 	\lim_{r\to 0^+} \per \left( \imt(\boldsymbol{y},\boldsymbol{0})  \right)=\pi + 1>\pi=\per \left(  \imt(\boldsymbol{y},\boldsymbol{0}) \right).
	 \end{equation*}
	Therefore, the deformation $\boldsymbol{y}$ violates \eqref{eq:conv-perimeter}.
\end{example}

\EEE 

\section*{Appendix: Differentiabiliy of maps on hypersurfaces}

\setcounter{section}{0}
\setcounter{theorem}{0}

\setcounter{equation}{0}

\renewcommand{\thetheorem}{A.\arabic{theorem}}

\renewcommand{\theequation}{A.\arabic{equation}}

In this appendix, we revise some facts about differentiability of maps defined on  hypersurfaces. Most of the results are certainly known to experts and we recall them just for the convenience of the reader.

We assume that $S\subset \RN$ is a  hypersurface of class $C^1$ and, for any $\boldsymbol{x}_0\in S$, we denote by $\Tan_{\boldsymbol{x}_0}S$ the tangent space to $S$ at $\boldsymbol{x}_0$. By a local chart of $S$ at $\boldsymbol{x}_0$  mean any injective map $\boldsymbol{\eta}\in C^1(Q;\RN)$ with $Q\subset \R^{N-1}$ open and $\boldsymbol{z}_0\coloneqq \boldsymbol{\eta}^{-1}(\boldsymbol{x}_0)\in Q$ such that the differential $\d_{\boldsymbol{z}_0}\boldsymbol{\eta}$ of $\boldsymbol{\eta}$ at $\boldsymbol{z}_0$  is also injective. We denote by $\boldsymbol{\nu}_S\colon S \to  S(\boldsymbol{0},1)$ a given unit normal field on $S$.  Below, $\mathscr{L}^{N-1}$ indicates the $(N-1)$-dimensional Lebesgue measure on $\R^{N-1}$. \EEE 

Below, we will make use of the notion of pseudoinverse for  matrices \cite{penrose}. 
Given  $\boldsymbol{H}\in \R^{N\times (N-1)}$ having full rank, we define its pseudoinverse as $\boldsymbol{H}^{\dag}  \coloneqq (\boldsymbol{H}^\top \boldsymbol{H} )^{-1}\boldsymbol{H}^\top$. Note that $\boldsymbol{H}^\top \boldsymbol{H} \in \R^{(N-1)\times (N-1)}$ is positive definite, so that $\boldsymbol{H}^\dag$ is well defined. Also, we observe that
\begin{equation}
	\label{eq:penrose1}
	\boldsymbol{H}^\dag \boldsymbol{H}=\boldsymbol{I},
\end{equation}
while
\begin{equation}
	\label{eq:penrose2}
	\text{$\boldsymbol{H}\boldsymbol{H}^\dag \boldsymbol{v}$ is the  orthogonal  projection of $\boldsymbol{v}$ onto the range of $\boldsymbol{H}$  for all $\boldsymbol{v} \in \RN$ .}
\end{equation}
The identity \eqref{eq:penrose1} is immediate. For the proof of \eqref{eq:penrose2}, see \cite[Lemma C.3]{skrepek}.

\subsection*{Approximate tangential differentiability} We refer to  \cite[Subsection 3.2.16]{federer} for the definition of approximate tangential differentiability. Precisely, according to the terminology in \cite{federer}, we are interested in the notion of  $(\haus\mres S,N-1)$-approximate differentiability. If $\boldsymbol{y}\colon S \to \RN$ is approximately tangentially differentiable at $\boldsymbol{x}_0\in S$, then the approximate tangential differential of $\boldsymbol{y}$ at $\boldsymbol{x}_0$ is a linear map $\d^{\rm t}_{\boldsymbol{x}_0}\boldsymbol{y}\colon \Tan_{\boldsymbol{x}_0}S \to \RN$ and we define the approximate tangential gradient as the unique matrix $\nabla^{\rm t}\boldsymbol{y}(\boldsymbol{x}_0)\in\RNN$ characterized by the two properties
\begin{equation}
	\label{eq:nabla-t}
	\d^{\rm t}_{\boldsymbol{x}_0}\boldsymbol{y}(\boldsymbol{v})=(\nabla^{\rm t}\boldsymbol{y}(\boldsymbol{x}_0))\boldsymbol{v} \quad \text{for all $\boldsymbol{v}\in\Tan_{\boldsymbol{x}_0}S$,} \qquad (\nabla^{\rm t}\boldsymbol{y}(\boldsymbol{x}_0))\boldsymbol{\nu}_S(\boldsymbol{x}_0)=\boldsymbol{0}.
\end{equation}

The approximate tangential gradient can be described using local charts as follows.

\begin{lemma}[Approximate tangential gradient by charts] \label{lem:diff-charts}
Let $\boldsymbol{y}\colon S \to \RN$ and $\boldsymbol{x}_0\in S$. Then, for any local chart $\boldsymbol{\eta}\in C^1(Q;\RN)$ of $S$ at $\boldsymbol{x}_0$, we have that $\boldsymbol{y}$ is approximately tangentially differentiable at $\boldsymbol{x}_0$  if and only if $\boldsymbol{u}\coloneqq \boldsymbol{y}\circ \boldsymbol{\eta}\colon Q\to \RN$ is approximately differentiable at  $\boldsymbol{z}_0\coloneqq \boldsymbol{\eta}^{-1}(\boldsymbol{x}_0)$. In this case,  we have 
 \begin{equation}
	\label{eq:tg-gc2}
	\nabla\boldsymbol{u}(\boldsymbol{z}_0)=\nabla^{\rm t}\boldsymbol{y} (\boldsymbol{x}_0)D\boldsymbol{\eta}(\boldsymbol{z}_0)
\end{equation}
or, equivalently, 
\begin{equation}
	\label{eq:tg-gc}
	\nabla^{\rm t}\boldsymbol{y}(\boldsymbol{x}_0)=\nabla\boldsymbol{u}(\boldsymbol{z}_0) (D\boldsymbol{\eta}(\boldsymbol{z}_0))^{\text{\rm \textdied}}.
\end{equation}  
\end{lemma}
Here, $\nabla \boldsymbol{u}$ denotes the approximate gradient of $\boldsymbol{u}$. Thus, if $\d_{\boldsymbol{z}_0} \boldsymbol{u}\colon \R^{N-1}\to \RN $ is the approximate differential of $\boldsymbol{u}$ at $\boldsymbol{z}_0$, then  $\d_{\boldsymbol{z}_0}\boldsymbol{u}=\d^{\rm t}_{\boldsymbol{x}_0}\boldsymbol{y}\circ \d_{\boldsymbol{z}_0} \boldsymbol{\eta}$ in view  of \eqref{eq:tg-gc2}.  
The equivalence of \eqref{eq:tg-gc2}--\eqref{eq:tg-gc} follows immediately from \eqref{eq:penrose1}--\eqref{eq:penrose2}. 
\begin{proof}
We prove only the sufficiency as the proof of necessity is analogous. Suppose that $\boldsymbol{u}$ is approximately differentiable at $\boldsymbol{z}_0$. By a standard property of approximate limits, we know that for some measurable set $R\subset Q$ with  
\begin{equation*}
	\Theta^{N-1}(R,\boldsymbol{z}_0)\coloneqq \lim_{ r  \to 0^+} \frac{\mathscr{L}^{N-1}(R\cap B(\boldsymbol{z}_0,r))}{\omega_{N-1}r^{N-1}}=0,
\end{equation*}
where  $\omega_{N-1}\coloneqq \mathscr{L}^{N-1}(B(\boldsymbol{0},1))$, \EEE we have
\begin{equation*}
	\lim_{\substack{ \boldsymbol{z}\to\boldsymbol{z}_0\\ \boldsymbol{z}\in Q\setminus R }} \frac{|\boldsymbol{u}(\boldsymbol{z})-\boldsymbol{u}(\boldsymbol{z}_0)-\nabla\boldsymbol{u}(\boldsymbol{z})(\boldsymbol{z}-\boldsymbol{z}_0)|}{|\boldsymbol{z}-\boldsymbol{z}_0|}=0.
\end{equation*} 
Setting $T\coloneqq \boldsymbol{\eta}(R)\subset S$, one checks that 
\begin{equation*}
	\Theta^{N-1}(T,\boldsymbol{x}_0)\coloneqq \lim_{r\to 0^+} \frac{\haus(T\cap B(\boldsymbol{x}_0,r))}{\alpha_{N-1}r^{N-1}}=0,	
\end{equation*}
where $\alpha_{N-1}\coloneqq \haus(S(\boldsymbol{0},1))$. Then,
for $\boldsymbol{z}\in Q\setminus R$ and $\boldsymbol{x}\coloneqq \boldsymbol{\eta}(\boldsymbol{z})$, we write
\begin{equation*}
	\begin{split}
		|\boldsymbol{y}(\boldsymbol{x})&-\boldsymbol{y}(\boldsymbol{x}_0)-\nabla\boldsymbol{u}(\boldsymbol{z}_0)(D\boldsymbol{\eta}(\boldsymbol{z}_0))^\dag(\boldsymbol{x}-\boldsymbol{x}_0)|\\
		&=|\boldsymbol{u}(\boldsymbol{z})-\boldsymbol{u}(\boldsymbol{z}_0)-\nabla\boldsymbol{u}(\boldsymbol{z}_0)(D\boldsymbol{\eta}(\boldsymbol{z}_0))^\dag \left( D\boldsymbol{\eta}(\boldsymbol{z}_0)(\boldsymbol{z}-\boldsymbol{z}_0)+\mathrm{o}(|\boldsymbol{z}-\boldsymbol{z}_0|)  \right)|\\
		&=|\boldsymbol{u}(\boldsymbol{z})-\boldsymbol{u}(\boldsymbol{z}_0)-\nabla\boldsymbol{u}(\boldsymbol{z}_0)(\boldsymbol{z}-\boldsymbol{z}_0)+\mathrm{o}(|\boldsymbol{z}-\boldsymbol{z}_0|)|,
	\end{split}
\end{equation*}
as $\boldsymbol{z}\to \boldsymbol{z}_0$. Here we used \eqref{eq:penrose1}.
Note that $\frac{|\boldsymbol{x}-\boldsymbol{x}_0|}{|\boldsymbol{z}-\boldsymbol{z}_0|}$ is controlled  both from above and below by the Lipschitz constants of $\boldsymbol{\eta}$ and $\boldsymbol{\eta}^{-1}$, so that $\boldsymbol{z}\to \boldsymbol{z}_0$ is equivalent to $\boldsymbol{x}\to \boldsymbol{x}_0$. Hence, we  obtain
\begin{equation*}
	\begin{split}
		\lim_{\substack{\boldsymbol{x}\to \boldsymbol{x}_0 \\ \boldsymbol{x}\in \boldsymbol{\eta}(Q)\setminus T }} \frac{|\boldsymbol{y}(\boldsymbol{x})-\boldsymbol{y}(\boldsymbol{x}_0)-\nabla\boldsymbol{u}(\boldsymbol{z}_0)(D\boldsymbol{\eta}(\boldsymbol{z}_0))^\dag(\boldsymbol{x}-\boldsymbol{x}_0)|}{|\boldsymbol{x}-\boldsymbol{x}_0|}=0,
	\end{split}
\end{equation*}  
which proves that $\boldsymbol{y}$ is approximately tangentially differentiable at $\boldsymbol{x}_0$ with  approximate tangential gradient as in \eqref{eq:tg-gc}  by same property of approximate limits as above.
\end{proof}

Henceforth, the adjoint of an operator is indicated with a star superscript. 
We define the approximate tangential Jacobian as in \cite{ambrosio.fusco.pallara,maggi,simon}.

\begin{definition}[Approximate tangential Jacobian]\label{def:atj}
Let  $\boldsymbol{y}\colon S \to \RN$ be approximately tangentially differentiable at $\boldsymbol{x}_0\in S$. Then, the approximate tangential Jacobian of $\boldsymbol{y}$ at $\boldsymbol{x}_0$ is defined as
\begin{equation*}
	\label{eq:atj}
	J^{\rm t}\boldsymbol{y}(\boldsymbol{x}_0)\coloneqq \sqrt{\det \left ((\d^{\rm t}_{\boldsymbol{x}_0}\boldsymbol{y})^*\circ (\d^{\rm t}_{\boldsymbol{x}_0}\boldsymbol{y})\right) }.
\end{equation*} 	
\end{definition}
 Up to identifying $\Tan_{\boldsymbol{x}_0}S$ and $\RN$ with their duals, the operator  $(\d^{\rm t}_{\boldsymbol{x}_0}\boldsymbol{y})^*\circ \d_{\boldsymbol{x}_0}^{\rm t}\boldsymbol{y}\colon \Tan_{\boldsymbol{x}_0}S \to \Tan_{\boldsymbol{x}_0}S$  acts as $\boldsymbol{v}\mapsto (\nabla^{\rm t}\boldsymbol{y}(\boldsymbol{x}_0))^\top(\nabla^{\rm t}\boldsymbol{y}(\boldsymbol{x}_0))\boldsymbol{v}$ for each $\boldsymbol{v}\in\Tan_{\boldsymbol{x}_0}S$. In this regard, \eqref{eq:nabla-t} ensures that $(\nabla^{\rm t}\boldsymbol{y}(\boldsymbol{x}_0))^\top(\nabla^{\rm t}\boldsymbol{y}(\boldsymbol{x}_0))\boldsymbol{v}\in \Tan_{\boldsymbol{x}_0}S$.
The next result provides a first way to compute the approximate tangential Jacobian using the approximate tangential gradient.

\begin{lemma}[Approximate tangential Jacobian via cofactor]
	\label{lem:j}
Let $\boldsymbol{y}\colon S \to \RN$ be approximately tangentially differentiable at $\boldsymbol{x}_0\in S$. Then,
\begin{equation}
	\label{eq:j}
	J^{\rm t}\boldsymbol{y}(\boldsymbol{x}_0)=\begin{cases}
		|(\cof \nabla^{\rm t}\boldsymbol{y}(\boldsymbol{x}_0))\boldsymbol{\nu}_S(\boldsymbol{x}_0)| & \text{if $\d_{\boldsymbol{x}_0}^{\rm t}\boldsymbol{y}$ is injective,}\\
		0 & \text{otherwise.}
	\end{cases}
\end{equation}	
\end{lemma}
\begin{proof}
For convenience, we set $\boldsymbol{L}\coloneqq \d^{\rm t}_{\boldsymbol{x}_0}\boldsymbol{y}$, $V\coloneqq \Tan_{\boldsymbol{x}_0}S$, and $\boldsymbol{F}\coloneqq \nabla^{\rm t}\boldsymbol{y}(\boldsymbol{x}_0)$. Let  $(\boldsymbol{w}_1,\dots,\boldsymbol{w}_d)$ be an orthonormal basis of $\boldsymbol{L}(V)$, where $d\leq N-1$ denotes the dimension of $\boldsymbol{L}(V)$. We consider  $N-1-d$ vectors $\boldsymbol{w}_{d+1},\dots,\boldsymbol{w}_{N-1}\in\RN$ such that the $(N-1)$-tuple $(\boldsymbol{w}_1,\dots,\boldsymbol{w}_{N-1})$ is orthonormal and, in turn, the space $W$ spanned by these vectors is a hyperplane. Define $\widetilde{\boldsymbol{L}}\colon V \to W$ by setting $\widetilde{\boldsymbol{L}}(\boldsymbol{v})\coloneqq \boldsymbol{F}\boldsymbol{v}$. Up to the identification of finite-dimensional spaces with their duals, we have  $\boldsymbol{L}^*\circ \boldsymbol{L}=\widetilde{\boldsymbol{L}}^*\circ \widetilde{\boldsymbol{L}}$. In particular, given an orthonormal basis $(\boldsymbol{v}_1,\dots,\boldsymbol{v}_{N-1})$ of $V$, basic linear algebra arguments show that 
$\det (\widetilde{\boldsymbol{L}}^*\circ \widetilde{\boldsymbol{L}})=(\det \boldsymbol{A})^2$, where $\boldsymbol{A}\in \R^{(N-1)\times(N-1)}$ is the matrix representing $\widetilde{\boldsymbol{L}}$ with respect to the basis $(\boldsymbol{v}_1,\dots,\boldsymbol{v}_{N-1})$ of $V$ and $(\boldsymbol{w}_1,\dots,\boldsymbol{w}_{N-1})$ of $W$, namely,
\begin{equation*}
	A^i_j\coloneqq \boldsymbol{F}\boldsymbol{v}_j\cdot \boldsymbol{w}_i \quad \text{for all $i,j=1,\dots,N-1$,}
\end{equation*}
 where the superscript $i$ indicates the row and the subscript $j$ indicates the column.   Note that $\boldsymbol{L}$ is injective if and only if the same holds for $\widetilde{\boldsymbol{L}}$. If $\boldsymbol{L}$ is not injective, then $\det \boldsymbol{A}=0$ which proves the second case in \eqref{eq:j}. Therefore, let us assume that $\boldsymbol{L}$ is injective. Setting $\widehat{\boldsymbol{F}}\coloneqq \boldsymbol{F}+\boldsymbol{\nu}\otimes \boldsymbol{\nu}$, where we write $\boldsymbol{\nu}$ for  $\boldsymbol{\nu}_S(\boldsymbol{x}_0)$, the map $\widehat{\boldsymbol{L}}\colon \RN \to \RN$ defined as $\widehat{\boldsymbol{L}}(\boldsymbol{v})\coloneqq \widehat{\boldsymbol{F}}\boldsymbol{v}$ provides an injective extension of $\boldsymbol{L}$ in view of \eqref{eq:tg-gc}.  Let  $\widetilde{\boldsymbol{\nu}}\coloneqq \frac{(\cof \boldsymbol{F})\boldsymbol{\nu} }{|(\cof \boldsymbol{F})\boldsymbol{\nu}|}$, and observe that $\widetilde{\boldsymbol{\nu}} \in W^\bot$ since for each $\boldsymbol{w} \in W$, with $\boldsymbol{v} \in V$ such that $\boldsymbol{w}  = \boldsymbol{F} \boldsymbol{v}$, we have
$\boldsymbol{w} \cdot  (\cof \boldsymbol{F})\boldsymbol{\nu}  = (\adj \boldsymbol{F}) \boldsymbol{F} \boldsymbol{v} \cdot \boldsymbol{\nu} = (\det \boldsymbol{F}) \boldsymbol{v} \cdot \boldsymbol{\nu} = 0$.   Let   $\boldsymbol{B}\in\RNN$ be  the matrix   representing $\widehat{\boldsymbol{L}}$ with respect to the two orthonormal basis $(\boldsymbol{v}_1,\dots,\boldsymbol{v}_{N-1},\boldsymbol{\nu})$ and $(\boldsymbol{w}_1,\dots,\boldsymbol{w}_{N-1},\widetilde{\boldsymbol{\nu}})$.  The matrix satisfies   $B^i_j=A^i_j$ for $i,j=1,\dots,N-1$, $B^N_j = 0$ for $j=1,\ldots, N-1$,   and $B^N_N=\frac{\det \widehat{\boldsymbol{F}}}{|(\cof\boldsymbol{F})\boldsymbol{\nu}|}$. Thus, taking into account that $\det \widehat{\boldsymbol{F}}\neq 0$, we find
\begin{equation*}
	|\det \widehat{\boldsymbol{F}}|=|\det \boldsymbol{B}|=\frac{|\det \boldsymbol{A}|\,|\det \widehat{\boldsymbol{F}}|}{|(\cof \boldsymbol{F})\boldsymbol{\nu}|}.
\end{equation*}
This yields $|\det \boldsymbol{A}|=|(\cof \boldsymbol{F})\boldsymbol{\nu}|$ which proves the first case in \eqref{eq:j}.
\end{proof}

In some literature, such as \cite[Section 3.2.1]{federer}, the approximate tangential Jacobian is often defined in terms of exterior algebra. As shown in the next lemma, that  choice  is actually equivalent to  Definition \ref{def:atj}.  We will employ customary notation of differential geometry for which we refer to \cite[Chapter 2]{federer}.

\begin{lemma}[Approximate tangential Jacobian via exterior algebra]
	\label{lem:j-ea}
Let $\boldsymbol{y}\colon S \to \RN$ be approximately tangentially differentiable at $\boldsymbol{x}_0\in S$. Then, we have
\begin{equation*}
	J^{\rm t}\boldsymbol{y}(\boldsymbol{x}_0)=\|\Lambda_{N-1}(\d^{\rm t}_{\boldsymbol{x}_0}\boldsymbol{y})\|.
\end{equation*}	
\end{lemma}
Here, $\|\Lambda_{N-1}(\d^{\rm t}_{\boldsymbol{x}_0}\boldsymbol{y})\|$ denotes the operator norm of $\Lambda_{N-1}(\d^{\rm t}_{\boldsymbol{x}_0}\boldsymbol{y})\colon \Lambda_{N-1}(\Tan_{\boldsymbol{x}_0}S)\to\Lambda_{N-1}\RN$.
\begin{proof}	
Also here, we set $\boldsymbol{L}\coloneqq \d^{\rm t}_{\boldsymbol{x}_0}\boldsymbol{y}$, $V\coloneqq \Tan_{\boldsymbol{x}_0}S$, and $\boldsymbol{F}\coloneqq \nabla^{\rm t}\boldsymbol{y}(\boldsymbol{x}_0)$.  The matrix $\boldsymbol{A}\in\R^{(N-1)\times (N-1)}$ representing $\boldsymbol{L}^*\circ \boldsymbol{L}$ with respect to an orthonormal basis $(\boldsymbol{v}_1,\dots,\boldsymbol{v}_{N-1})$ of $V$ is given by
\begin{equation*}
	A^i_j=\left( \boldsymbol{F}^\top \boldsymbol{F} \boldsymbol{v}_j \right) \cdot \boldsymbol{v}_i=(\boldsymbol{F}\boldsymbol{v}_j)\cdot (\boldsymbol{F}\boldsymbol{v}_i) \quad \text{for all $i,j=1,\dots,N-1$.}
\end{equation*}
Also, by definition, $J^{\rm t}\boldsymbol{y}(\boldsymbol{x})=\sqrt{\det \boldsymbol{A}}$.
The  one-dimensional  space $\Lambda_{N-1}V$ is spanned by $\boldsymbol{v}_1\wedge \dots \wedge \boldsymbol{v}_{N-1}$ and the induced map $\Lambda_{N-1}\boldsymbol{L}\colon \Lambda_{N-1}V \to \Lambda_{N-1}\RN$ uniquely determined by the identity
\begin{equation*}
	(\Lambda_{N-1}\boldsymbol{L})(\boldsymbol{v}_1 \wedge \dots \boldsymbol{v}_{N-1})=\boldsymbol{L}(\boldsymbol{v}_1)\wedge \dots \wedge \boldsymbol{L}(\boldsymbol{v}_{N-1})=(\boldsymbol{F}\boldsymbol{v}_1)\wedge \dots \wedge (\boldsymbol{F}\boldsymbol{v}_{N-1})
\end{equation*}
Recall that 
\begin{equation*}
	|\boldsymbol{w}_1\wedge \dots \wedge \boldsymbol{w}_{N-1}|^2=\det \boldsymbol{B} \quad \text{for all $\boldsymbol{w}_1,\dots,\boldsymbol{w}_{N-1}\in\RN$,}
\end{equation*}
where $\boldsymbol{B}\in \R^{(N-1)\times (N-1)}$ is defined by setting $B^i_j\coloneqq \boldsymbol{w}_j\cdot \boldsymbol{w}_i$ for all $i,j=1,\dots,N-1$.
Hence,  $|\boldsymbol{v}_1 \wedge \dots  \wedge\boldsymbol{v}_{N-1}|=1$ and, in turn, 
\begin{equation*}
	\|\Lambda_{N-1}\boldsymbol{L}\|^2=|(\boldsymbol{F}\boldsymbol{v}_1)\wedge \dots \wedge (\boldsymbol{F}\boldsymbol{v}_{N-1})|^2=\det \boldsymbol{A},
\end{equation*}
which proves the claim.
\end{proof}

As a consequence of Lemma \ref{lem:diff-charts}, the approximate tangential Jacobian can be also computed by using local charts.

\begin{lemma}[Approximate tangential Jacobian by charts]
	\label{lem:jac-charts}
Let $\boldsymbol{y}\colon S \to \RN$ be approximately tangentially differentiable at $\boldsymbol{x}_0\in S$. Then, for any local chart $\boldsymbol{\eta}\in C^1(Q;\RN)$ of $S$ at $\boldsymbol{x}_0$, setting $\boldsymbol{z}_0\coloneqq \boldsymbol{\eta}^{-1}(\boldsymbol{x}_0)$ and  $\boldsymbol{u}\coloneqq \boldsymbol{y}\circ \boldsymbol{\eta}\colon Q\to \RN$, we have 
\begin{equation*}
	J\boldsymbol{u}(\boldsymbol{z}_0)=J^{\rm t}\boldsymbol{y}(\boldsymbol{x}_0)\,J\boldsymbol{\eta}(\boldsymbol{z}_0).
\end{equation*}	
\end{lemma}
Here, $J\boldsymbol{u}(\boldsymbol{z}_0)\coloneqq \sqrt{\det \left( (\d_{\boldsymbol{z}_0}\boldsymbol{u})^*\circ \d_{\boldsymbol{z}_0}\boldsymbol{u}   \right)}$ is the approximate Jacobian of $\boldsymbol{u}$ at $\boldsymbol{z}_0$ and, analogously, $J\boldsymbol{\eta}(\boldsymbol{z}_0)\coloneqq \sqrt{\det \left( (\d_{\boldsymbol{z}_0}\boldsymbol{\eta})^*\circ \d_{\boldsymbol{z}_0}\boldsymbol{\eta}   \right)}$ is the tangential Jacobian of $\boldsymbol{\eta}$ at the same point.
\begin{proof}
We argue as in \cite[Section 2, Chapter 4]{simon}. Note that $\boldsymbol{u}$ is approximately differentiable at $\boldsymbol{z}_0$ by Lemma~\ref{lem:diff-charts}. Let  $\boldsymbol{A}\in\R^{(N-1)\times (N-1)}$ be defined as
\begin{equation*}
	A^i_j=\partial_i \boldsymbol{u}(\boldsymbol{z}_0) \cdot \partial_j \boldsymbol{u}(\boldsymbol{z}_0) \quad \text{for all $i,j=1,\dots,N-1$.}
\end{equation*}
By definition, $J\boldsymbol{u}(\boldsymbol{z}_0)=\sqrt{\det \boldsymbol{A}}$. 
 Given an orthonormal basis $(\boldsymbol{v}_1,\dots,\boldsymbol{v}_{N-1})$ of $\Tan_{\boldsymbol{x}_0}S$, using Lemma \ref{lem:diff-charts}, we write
\begin{equation*}
	\partial_i \boldsymbol{u}(\boldsymbol{z}_0)=D^{\rm t}\boldsymbol{y}(\boldsymbol{x}_0)\partial_i\boldsymbol{\eta}(\boldsymbol{z}_0)=\sum_{k=1}^{N-1} (\partial_i \boldsymbol{\eta}(\boldsymbol{z}_0)\cdot \boldsymbol{v}_k) \nabla^{\rm t}\boldsymbol{y}(\boldsymbol{x}_0)\boldsymbol{v}_k.
\end{equation*} 
Thus,
\begin{equation}
	\label{eq:abc}
	A^i_j=\sum_{k,l=1}^{N-1} (\partial_i \boldsymbol{\eta}(\boldsymbol{z}_0)\cdot \boldsymbol{v}_k) \, (\partial_j \boldsymbol{\eta}(\boldsymbol{z}_0)\cdot \boldsymbol{v}_l) \left( \nabla^{\rm t}\boldsymbol{y}(\boldsymbol{x}_0)\boldsymbol{v}_k \right) \cdot  \left ( \nabla^{\rm t}\boldsymbol{y}(\boldsymbol{x}_0)\boldsymbol{v}_l \right).
\end{equation}
Now, if we define $\boldsymbol{B},\boldsymbol{C}\in\R^{(N-1)\times (N-1)}$ by setting
\begin{equation*}
	B^i_j\coloneqq \partial_i \boldsymbol{\eta}(\boldsymbol{z}_0)\cdot \boldsymbol{v}_j, \quad C^i_j\coloneqq \left(\nabla^{\rm t}\boldsymbol{y}(\boldsymbol{x}_0)\boldsymbol{v}_i \right) \cdot \left( \nabla^{\rm t}\boldsymbol{y}(\boldsymbol{x}_0)\boldsymbol{v}_j \right) \qquad \text{for all $i,j=1,\dots,N-1$,}
\end{equation*}
then, the identity \eqref{eq:abc} yields $\boldsymbol{A}=\boldsymbol{B}\boldsymbol{C}\boldsymbol{B}^\top$ and, in turn, $\det \boldsymbol{A}=(\det \boldsymbol{B})^2(\det \boldsymbol{C})$. Therefore, by observing that $J\boldsymbol{\eta}(\boldsymbol{z}_0)=|\det \boldsymbol{B}|$ and $J^{\rm t}\boldsymbol{y}(\boldsymbol{x}_0)=\sqrt{\det \boldsymbol{C}}$, we recover the desired formula.
\end{proof}

With Lemma \ref{lem:diff-charts} and Lemma \ref{lem:jac-charts} at hand, the area formula for maps on hypersurfaces follows easily from the classical area formula. In the case of Lipschitz functions, the result is given in \cite[Corollary 3.2.20]{federer} for arbitrary codimension taking into account Lemma \ref{lem:j-ea}.

\begin{theorem}[Area formula for maps on hypersurfaces]
	\label{thm:area-formula surface}
Let  $\boldsymbol{y}\colon S \to \RN$ be $\haus$-almost everywhere approximately tangentially differentiable. Denote by $X_{\boldsymbol{y}}\subset S$  the set of approximate tangential differentiability points of $\boldsymbol{y}$.  Then, for any $\haus$-measurable set $E\subset S$, the function $\boldsymbol{\xi}\mapsto m(\boldsymbol{y},E,\boldsymbol{\xi})\coloneqq \H^0( \left\{ \boldsymbol{x}\in E \cap X_{\boldsymbol{y}}: \:\boldsymbol{y}(\boldsymbol{x})=\boldsymbol{\xi} \right\})$ is $\haus$-measurable. Moreover, for any $\haus$-measurable map $\psi\colon \boldsymbol{y}(E\cap X_{\boldsymbol{y}})\to \R$, we have 
\begin{equation}\label{eq:af-s2}
	\int_{\RN} \psi(\boldsymbol{\xi}) m(\boldsymbol{y},E,\boldsymbol{\xi})\,\d\haus(\boldsymbol{\xi})=\int_{E} \psi(\boldsymbol{y}(\boldsymbol{x})) J^{\rm t}\boldsymbol{y}(\boldsymbol{x}) \,\d\haus(\boldsymbol{x})
\end{equation}
whenever one of the two integrals is well defined.
\end{theorem}
 Note that  on the left-hand side of \eqref{eq:af-s2} the domain of integration can be replaced by $\boldsymbol{y}(E\cap X_{\boldsymbol{y}})$. 
\begin{proof}
By splitting $E$ in a suitably way, we can reduce to the case $E\subset \boldsymbol{\eta}(Q)$ with  $\boldsymbol{\eta}\in C^1(Q;\RN)$ being a local chart of $S$. Let $\boldsymbol{u}\coloneqq \boldsymbol{y}\circ \boldsymbol{\eta}$ and recall Lemma \ref{lem:diff-charts}. If $Z_{\boldsymbol{u}}\subset Q$ denotes the set of  approximate differentiability points of $\boldsymbol{u}$, then  $X_{\boldsymbol{y}}=\boldsymbol{\eta}(Z_{\boldsymbol{u}})$ and, setting $F\coloneqq \boldsymbol{\eta}^{-1}(E)$, we have
\begin{equation*}
	m(\boldsymbol{u},F,\boldsymbol{\xi})\coloneqq \H^0(\{ \boldsymbol{z}\in F\cap Z_{\boldsymbol{u}}:\:\boldsymbol{u}(\boldsymbol{z})=\boldsymbol{\xi}  \})=\H^0(\{  \boldsymbol{x}\in E \cap X_{\boldsymbol{y}}:\:\boldsymbol{y}(\boldsymbol{x})=\boldsymbol{\xi}  \})=	m(\boldsymbol{y},E,\boldsymbol{\xi})
\end{equation*}
for all $\boldsymbol{\xi}\in\RN$. Using Federer's area formula for approximately differentiable maps \cite[Theorem 1, p. 220]{gms.cc1}  and Lemma \ref{lem:jac-charts}, we compute
\begin{equation*}
	\begin{split}
			\int_{\RN}\psi(\boldsymbol{\xi})m(\boldsymbol{y},E,\boldsymbol{\xi})\,\d\haus(\boldsymbol{\xi})&=\int_{\RN} \psi(\boldsymbol{\xi}) m(\boldsymbol{u},F,\boldsymbol{\xi})\,\d\haus(\boldsymbol{\xi})=\int_F \psi(\boldsymbol{u}(\boldsymbol{z})) J\boldsymbol{u}(\boldsymbol{z})\,\d \boldsymbol{z}\\
			&=\int_F \psi(\boldsymbol{y}(\boldsymbol{\eta}(\boldsymbol{z}))) J^{\rm t}\boldsymbol{y}(\boldsymbol{\eta}(\boldsymbol{z}))\,J\boldsymbol{\eta}(\boldsymbol{z})\,\d\boldsymbol{z}=\int_E \psi(\boldsymbol{y}(\boldsymbol{x}))J^{\rm t}\boldsymbol{y}(\boldsymbol{x})\,\d\boldsymbol{x}.
	\end{split}
\end{equation*}
\end{proof}

\subsection*{Weak tangential differentiability}

Henceforth, we fix $p\in [1,\infty]$. The Sobolev space $W^{1,p}(S;\RN)$ is customarily defined as the class of maps $\boldsymbol{y}\in L^p(S;\RN)$ such that, for any local chart $\boldsymbol{\eta}\in C^1(Q;\RN)$ of $S$, we have $\boldsymbol{u}\coloneqq \boldsymbol{y}\circ \boldsymbol{\eta}\in W^{1,p}(Q;\RN)$.  Similarly to \eqref{eq:tg-gc},  the weak tangential gradient is defined as follows.  

\begin{definition}[Weak tangential gradient]\label{def:wtg}
Let $\boldsymbol{y}\in W^{1,p}(S;\RN)$. Then, for any local chart $\boldsymbol{\eta}\in C^1(Q;\RN)$ of $S$, the weak tangential gradient of $\boldsymbol{y}$ in $\boldsymbol{\eta}(Q)$ is defined as 
\begin{equation*}
	D^{\rm t}\boldsymbol{y}(\boldsymbol{x})\coloneqq D\boldsymbol{u}(\boldsymbol{z})(D\boldsymbol{\eta}(\boldsymbol{z}))^{\text{\rm \textdied}} \quad \text{for $\haus$-almost every $\boldsymbol{x}\in \boldsymbol{\eta}(Q)$,}
\end{equation*}
where $\boldsymbol{u}\coloneqq \boldsymbol{y}\circ \boldsymbol{\eta}$ and $\boldsymbol{z}\coloneqq \boldsymbol{\eta}^{-1}(\boldsymbol{x})$.
\end{definition}

The next result shows that the previous definition yields $D^{\rm t}\boldsymbol{y}$ as a well defined field on $S$, at least up to equivalence $\haus$-almost everywhere.

\begin{lemma}[Independence on local charts]
Let $\boldsymbol{y}\in W^{1,p}(S;\RN)$. Then, for any two local charts $\boldsymbol{\eta}\in C^1(Q;\RN)$ and $\widetilde{\boldsymbol{\eta}}\in C^1(\widetilde{Q};\RN)$ of $S$, we have
\begin{equation*}
	D\boldsymbol{u}(\boldsymbol{z})\left (D\boldsymbol{\eta}(\boldsymbol{z})\right)^{\text{\rm \textdied}} = D\widetilde{\boldsymbol{u}}(\widetilde{\boldsymbol{z}})\left(  D\widetilde{\boldsymbol{\eta}}(\widetilde{\boldsymbol{z}})\right)^{\text{\rm \textdied}} \quad \text{for $\haus$-almost every $\boldsymbol{x}\in\boldsymbol{\eta}(Q)\cap \widetilde{\boldsymbol{\eta}}(\widetilde{Q})$,}
\end{equation*}
where we set $\boldsymbol{u}\coloneqq \boldsymbol{y}\circ \boldsymbol{\eta}$, $\widetilde{\boldsymbol{u}}\coloneqq \boldsymbol{y}\circ \widetilde{\boldsymbol{\eta}}$, $\boldsymbol{z}\coloneqq \boldsymbol{\eta}^{-1}(\boldsymbol{x})$, and $\widetilde{\boldsymbol{z}}\coloneqq \widetilde{\boldsymbol{\eta}}^{-1}(\boldsymbol{x})$.	
\end{lemma} 
\begin{proof}
We argue along the lines of \cite[Appendix B]{skrepek}.	
First, we write ${\boldsymbol{\eta}}=\widetilde{\boldsymbol{\eta}}\circ \widetilde{\boldsymbol{\eta}}^{-1}\circ {\boldsymbol{\eta}}$ and we apply the chain rule to we find
\begin{equation*}
	D{\boldsymbol{\eta}}({\boldsymbol{z}})=D\widetilde{\boldsymbol{\eta}}(\widetilde{\boldsymbol{z}})D(\widetilde{\boldsymbol{\eta}}^{-1}\circ {\boldsymbol{\eta}})({\boldsymbol{z}}).
\end{equation*}	
By multiplying both sides of the previous equation on the left by $(D\widetilde{\boldsymbol{\eta}}(\widetilde{\boldsymbol{z}}))^\dag$, we obtain
\begin{equation}
	\label{eq:Dee}
	D(\widetilde{\boldsymbol{\eta}}^{-1}\circ {\boldsymbol{\eta}})({\boldsymbol{z}})=(D\widetilde{\boldsymbol{\eta}}(\widetilde{\boldsymbol{z}}))^\dag 	D{\boldsymbol{\eta}}({\boldsymbol{z}})
\end{equation}
thanks to \eqref{eq:penrose1}. Next, we write $\boldsymbol{u}=\widetilde{\boldsymbol{u}}\circ \widetilde{\boldsymbol{\eta}}\circ \boldsymbol{\eta}$ and we apply again the chain rule. We get
\begin{equation*}
	D\boldsymbol{u}(\boldsymbol{z})=D\widetilde{\boldsymbol{u}}(\widetilde{\boldsymbol{z}})D(\widetilde{\boldsymbol{\eta}}^{-1}\circ \boldsymbol{\eta})(\boldsymbol{z}),
\end{equation*}
which, together with \eqref{eq:Dee},  gives
\begin{equation*}
	D\boldsymbol{u}(\boldsymbol{z})=D\widetilde{\boldsymbol{u}}(\widetilde{\boldsymbol{z}})(D\widetilde{\boldsymbol{\eta}}(\widetilde{\boldsymbol{z}}))^\dag 	D{\boldsymbol{\eta}}({\boldsymbol{z}}). 
\end{equation*} 
 By multiplying  the both sides of the previous equation on the right by $(D\boldsymbol{\eta}(\boldsymbol{z}))^\dag$, we deduce 
\begin{equation*}
	D\boldsymbol{u}(\boldsymbol{z})\left(D{\boldsymbol{\eta}}({\boldsymbol{z}})\right)^\dag=D\widetilde{\boldsymbol{u}}(\widetilde{\boldsymbol{z}})(D\widetilde{\boldsymbol{\eta}}(\widetilde{\boldsymbol{z}}))^\dag D{\boldsymbol{\eta}}({\boldsymbol{z}})(D\boldsymbol{\eta}(\boldsymbol{z}))^\dag. 
\end{equation*}
As $D\boldsymbol{\eta}(\boldsymbol{z})$ and $D\widetilde{\boldsymbol{\eta}}(\widetilde{\boldsymbol{z}})$ have the same range $\Tan_{\boldsymbol{x}}S$, the conclusion follows from \eqref{eq:penrose2}.
\end{proof}

It is clear from Definition \ref{def:wtg} that $D^{\rm t}\boldsymbol{y}\in L^p(S;\RNN)$. Thus,  $W^{1,p}(S;\RN)$ can be equipped with the norm 
\begin{equation*}
	\|\boldsymbol{y}\|_{W^{1,p}(S;\RN)}\coloneqq \left( \|\boldsymbol{y}\|^p_{L^p(S;\RN)}+\|D^{\rm t}\boldsymbol{y}\|^p_{L^p(S;\RNN)}  \right)^{1/p}.
\end{equation*}  
It is routine to show that this  provides   a Banach space structure on  $W^{1,p}(S;\RN)$, as well as that   the usual characterization of weak convergence  holds true. 

 Next,    we discuss the  approximate tangential differentiability of Sobolev maps on hypersurfaces.

\begin{proposition}[Weak and approximate tangential gradients]
	\label{prop:watg}
Let $\boldsymbol{y}\in W^{1,p}(S;\RN)$. Then, $\boldsymbol{y}$ is approximately tangentially differentiable at $\boldsymbol{x}$ with approximate tangential gradient $D^{\rm t}\boldsymbol{y}(\boldsymbol{x})$ for $\haus$-almost every $\boldsymbol{x}\in S$. 
\end{proposition}
\begin{proof}
Let $\boldsymbol{\eta}\in C^1(Q;\RN)$ be a local chart of $S$. By definition, $\boldsymbol{u}\coloneqq \boldsymbol{y}\circ \boldsymbol{\eta}\in W^{1,p}(Q;\RN)$, so that  $\boldsymbol{u}$ is approximately differentiable at $\mathscr{L}^{N-1}$-almost every $\boldsymbol{z}\in Q$ with approximate gradient $D\boldsymbol{u}(\boldsymbol{z})$. For the previous claim, see, e.g., \cite[Theorem 2, p. 216]{gms.cc1}. Therefore, the conclusion follows in view of Lemma \ref{lem:diff-charts} and Definition \ref{def:wtg}.
\end{proof}

At this point, we realize that Proposition \ref{prop:cov-surface} is just a special case of  Theorem \ref{thm:area-formula surface} taking into account  Lemma \ref{lem:j} and  Proposition \ref{prop:watg}.

 Finally, we  address   the existence of continuous representatives for Sobolev maps on hypersurfaces. For the notion of tangential differentiability, we adopt the definition in \cite[Subsection 3.1.22]{federer}.

\begin{proposition}[Continuous representative]\label{prop:cont-repr}
Let $\boldsymbol{y}\in W^{1,p}(S;\RN)$ with $p>N-1$. Then, there exists a representative $\overline{\boldsymbol{y}}$ of $\boldsymbol{y}$ which is continuous, tangentially differentiable with tangential gradient $D^{\rm t}\boldsymbol{y}(\boldsymbol{x})$ at  $\mathscr{H}^{N-1}$-almost  every $\boldsymbol{x}\in S$, and satisfies Lusin's condition (N), namely
\begin{equation}\label{eq:LusinN}
	\haus(\overline{\boldsymbol{y}}(E))=0 \quad \text{whenever} \quad \text{$E\subset S$ satisfies $\haus(E)=0$.}
\end{equation}
\end{proposition}
In analogy with the classical case, one can actually prove that $\overline{\boldsymbol{y}}$ is locally (globally when $S$ is compact) H\"{o}lder-continuous with exponent $1-\frac{N-1}{p}$.
\begin{proof}
It is enough to show that all the desired properties hold on a coordinate patch. 	
Let $\boldsymbol{\eta}\in C^1(Q;\RN)$ be a local chart of $S$ and set $\boldsymbol{u}\coloneqq \boldsymbol{y}\circ \boldsymbol{\eta}\in W^{1,p}(Q;\RN)$. By \cite[Theorem 3, p.~223]{gms.cc1}, the map $\overline{\boldsymbol{u}}$ admits a representative $\overline{\boldsymbol{u}}\in C^0(\closure{Q};\RN)$ which is differentiable with gradient $D\boldsymbol{u}(\boldsymbol{z})$ for $\mathscr{L}^{N-1}$-almost every $\boldsymbol{z}\in Q$, and satisfies Lusin's condition (N) in the sense that 
 \begin{equation}
 	\label{eq:Nchart}
 	\haus(\overline{\boldsymbol{u}}(F))=0 \quad \text{whenever} \quad \text{$F\subset Q$ satisfies $\mathscr{L}^{N-1}(F)=0$.}
 \end{equation}
Define $\overline{\boldsymbol{y}}\coloneqq \overline{\boldsymbol{u}}\circ \boldsymbol{\eta}^{-1}$. Then, $\overline{\boldsymbol{y}}$ is continuous, and arguments analogous to those  of Lemma \ref{lem:diff-charts} show that $\overline{\boldsymbol{y}}$ is tangentially differentiable with tangential gradient $D^{\rm t}\boldsymbol{y}(\boldsymbol{x})$ for $\haus$-almost every $\boldsymbol{x}\in\boldsymbol{\eta}(Q)$.  Let $E\subset \boldsymbol{\eta}(Q)$ satisfy $\haus(E)=0$. Then, $\mathscr{L}^{N-1}(\boldsymbol{\eta}^{-1}(E))=0$, so that \eqref{eq:Nchart} gives $\haus(\overline{\boldsymbol{y}}(E))=\haus(\overline{\boldsymbol{u}}(\boldsymbol{\eta}^{-1}(E)))=0$.
\end{proof}

\section*{Acknowledgements}

M.~Bresciani acknowledges the support of the Alexander von Humboldt Foundation. He  thanks the Institute of Information Theory and Automation of the Czech Academy of Sciences, where part of this work has been developed, for the hospitality.   M.~Friedrich's research was funded by the Deutsche Forschungsgemeinschaft (DFG, German Research Foundation) - 377472739/GRK 2423/2-2023.  \EEE

\vskip 10pt
{
	%\textbf{Data availability statement.} Data sharing not applicable to this article as no datasets were generated or analyzed during the current study.
	}

\end{document}